\documentclass[11pt,a4paper,english, twoside]{book}
\usepackage[latin1]{inputenc}
\usepackage[T1]{fontenc}
\usepackage{babel}
\usepackage{amsmath,amscd,amssymb,latexsym}
\usepackage{epsfig,graphics,array,float}
\usepackage{makeidx}
\usepackage[latin1]{inputenc}
\usepackage[T1]{fontenc,url}
\usepackage{babel}
\usepackage{amsmath,amsfonts,amssymb,amsthm}
\usepackage{picins,graphicx,graphics,subfigure, subfig}
\usepackage{epsfig}
\usepackage{fancyhdr}
\usepackage{epic,eepic}
\usepackage{bm,color}
\usepackage{afterpage}
\usepackage{multirow}
\usepackage{wrapfig}

\theoremstyle{plain}
\newtheorem{thm}{Theorem}[section]
\newtheorem{lem}[thm]{Lemma}
\newtheorem{prp}[thm]{Proposition}
\newtheorem{cor}[thm]{Corollary}
\newtheorem{conj}[thm]{Conjecture}

\theoremstyle{definition}
\newtheorem{defi}[thm]{Definition}
\newtheorem{ex}[thm]{Example}

\theoremstyle{remark}
\newtheorem{rem}[thm]{Remark}

\newtheorem*{pf}{Proof}

\newcommand{\Po}{\ensuremath{\mathbb{P}^2\,}}

\newcommand{\V}{\ensuremath{\mathcal{V}}}
\newcommand{\noi}{\noindent}
\newcommand{\ol}{\overline}
\newcommand{\nom}{\noalign{\medskip}}
\newcommand{\beq}{\begin{equation*}}
\newcommand{\eeq}{\end{equation*}}
\newcommand{\bsp}{\begin{split}}
\newcommand{\esp}{\end{split}}
\newcommand{\scr}{\scriptsize}
\newcommand{\f}{\varphi}
\newcommand{\je}{\ell}

\makeindex

\author{Torgunn Karoline Moe}

\title{Rational Cuspidal Curves}

\begin{document}
\newpage
\thispagestyle{empty}
\pagenumbering{roman}



\begin{center}
  \qquad \\
  \qquad \\
  \qquad \\
  \Huge
  \textbf{Rational Cuspidal Curves} \\
  \normalsize
  \vspace{5mm}
  \textsl{by} \\
  \vspace{5mm}
  {\large
  \textbf{Torgunn Karoline Moe} \\}
  \vspace{35mm}
  {\normalsize
  {{\textsl{Thesis for the degree of} \\
  \vspace{2mm}}}}
  {\large
  {\bf{{Master in Mathematics}}} \\}
  \vspace{2mm}
  {\normalsize \textsl {(Master of Science)}}\\
  \vspace{5mm}
  \centerline{\includegraphics[width=7cm]{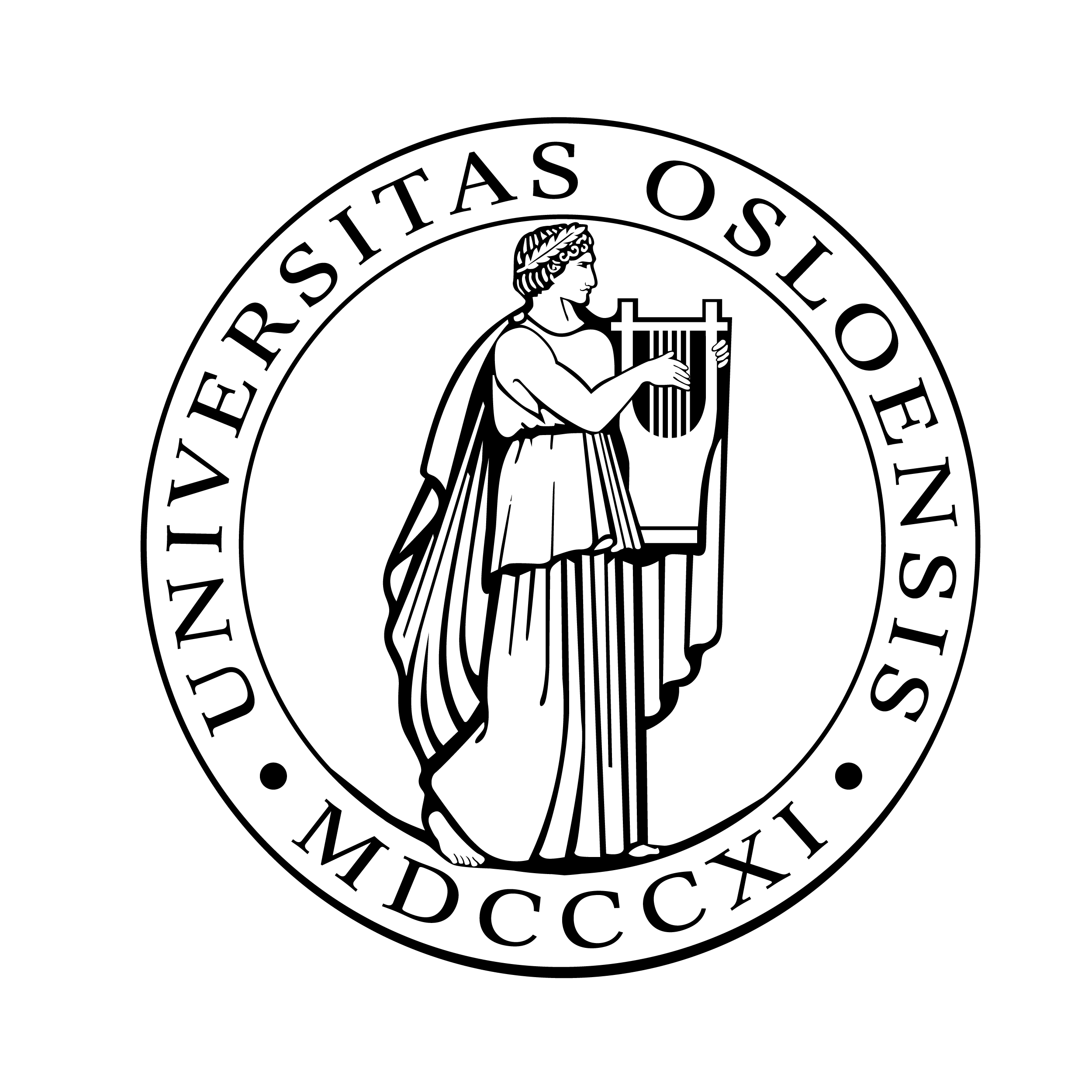}}
  \vspace{5mm}
  \textsl{Department of Mathematics} \\
  \textsl{Faculty of Mathematics and Natural Sciences} \\
  \textsl{University of Oslo} \\
  \vspace{5mm}
  \normalsize
  \textsl{May 2008} 
  \normalsize
\end{center}

\newpage
\thispagestyle{empty}
\cleardoublepage
\newpage
\thispagestyle{empty}
\cleardoublepage
\begin{center}
\qquad \\
\qquad \\
\qquad \\
  \vspace{25mm}
  {\Large
  \textbf{Rational Cuspidal Curves} \\}
  \vspace{5mm}
  \textsl{by} \\
  \vspace{5mm}
  \textbf{Torgunn Karoline Moe} \\
  \vspace{25mm}
  \textsl{Supervised by}\\
  \vspace{5mm}
 \textbf{Professor Ragni Piene}\\
   \vspace{50mm}
  \vspace{18mm}
    \textsl{Department of Mathematics} \\
    \textsl{Faculty of Mathematics and Natural Sciences} \\
  \textsl{University of Oslo} \\
  \vspace{5mm}
  \textsl{May 2008} \\
  \vspace{5mm}
  \end{center}
  
  \newpage
  \thispagestyle{empty}
  \cleardoublepage


 \newpage
  \thispagestyle{empty}
\cleardoublepage

\newpage
\thispagestyle{plain}
\include{PRE}

\newpage
\thispagestyle{empty}
\cleardoublepage
\thispagestyle{empty}

\newpage
\thispagestyle{plain}
\tableofcontents
\cleardoublepage

\newpage
\pagestyle{fancy}
\pagenumbering{arabic}
\chapter{Introduction}

A classical problem in algebraic geometry is the question of how many and what kind of singularities a plane curve of a given degree can have. This problem is interesting in itself. Additionally, the problem is interesting because it appears in other contexts, for example in the classification of open surfaces.

A curve in the projective plane is called rational if it is birational to a projective line. Furthermore, if all its singularities are cusps, we call the curve cuspidal. In this thesis we will investigate the above problem for rational cuspidal curves.

\begin{center}
{\em How many and what kind of cusps can a rational cuspidal curve have?}
\end{center}


\noi This problem has been boldly attacked with a variety of methods by a number of mathematicians. Some fundamental properties of rational cuspidal curves can be deduced from well known results in algebraic geometry. Additionally, very powerful results have been discovered recently. Rational cuspidal curves of low degree have been classified by Namba in \cite{Namba} and Fenske in \cite{Fen99b}. Series of rational cuspidal curves have been discovered and constructed by Fenske in \cite{Fen99b} and \cite{Fen99a}, Orevkov in \cite{Ore}, Tono in \cite{Tono00}, and Flenner and Zaidenberg in \cite{FlZa95} and \cite{FlZa97}. New properties of rational cuspidal curves have been found by Flenner and Zaidenberg in \cite{FlZa95}, Matsuoka and Sakai in \cite{MatsuokaSakai}, Orevkov in \cite{Ore}, Fernandez de Bobadilla et al. in \cite{Bobadillains}, and Tono in \cite{Tono05}.


Although a lot of technical tools have been developed, a definite answer to the above question has not been found. However, a vague contour of a partial, mysterious and intriguing answer has appeared. 

\begin{conj}\label{c:1}
A rational cuspidal curve can not have more than four cusps.
\end{conj}

\pagebreak
\noi In this thesis we present some of the results given in the mentioned works and give an overview of most known rational cuspidal curves. One very important tool in the mentioned works is Cremona transformations. We will therefore give a thorough definition of Cremona transformations and use them to construct some rational cuspidal curves of low degree. Moreover, a rational cuspidal curve in the plane can be viewed as a resulting curve of a projection of a curve in a higher-dimensional projective space. This represents a new and interesting way to approach such curves. In this thesis we will therefore also investigate the rational cuspidal curves from this point of view.\\

\noi In Chapter \ref{TB} we set notation and give an overview of the theoretical tools used in this thesis in the analyzation of rational cuspidal curves.

In Chapter \ref{rccq} we use some of the theoretical background to argue for the existence of the rational cuspidal cubic and quartic curves. We briefly introduce these curves by giving some essential properties of each curve.

In Chapter \ref{proj} we give a general description of how rational cuspidal curves can be constructed from the rational normal curve in a projective space. We will also analyze the cubic and quartic curves and the particular projections by which they can be constructed.

In Chapter \ref{ccc} we give a thorough definition of Cremona transformations. We use these transformations to construct and also investigate the construction of rational cuspidal cubic and quartic curves. In this process we encounter some issues concerning inflection points, which will be briefly discussed. Last in this chapter we present a conjecture linked to both rational cuspidal curves and Cremona transformations.

In Chapter \ref{rcq} we construct all rational cuspidal quintic curves with Cremona transformations and prove that they are the only rational cuspidal curves of this degree.

In Chapter \ref{mcc} we present a few series of rational cuspidal curves, some of which are just recently discovered.

In Chapter \ref{onc} we address the question of how many cusps a rational cuspidal curve can have, and we present the most recent discoveries on the problem. Two particular curves draw our attention, and these curves will be investigated in great detail. We additionally view the question from the perspective of projections.

In Chapter \ref{mrr} we present miscellaneous results which are closely related to rational cuspidal curves. First, we discuss whether all cusps on a cuspidal curve can have real coordinates. Second, we propose and investigate a conjecture concerning the intersection multiplicity of a curve and its Hessian curve. Third, we present an example of a reducible toric polar Cremona transformation.

\pagebreak
\noi The work in this thesis has led to neither a confirmation nor a contradiction of Conjecture \ref{c:1}. The thesis presents an overview of rational cuspial curves of low degree and explains how they can be constructed by Cremona transformations. Nothing new concerning cusps of a curve has been discovered in this work, but questions concerning the construction of inflection points have arisen. We have additionally shown that viewing rational cuspidal curves from the perspective of projection might introduce some new possibilities, but there are great obstacles blocking the way of new results, which we have not been able to step over. 

A possibly interesting subject for further investigations is how Cremona transformations can restrict the number of cusps of a rational cuspidal curve. Although there is no apparent way of attacking this problem generally, it seems to be strongly dependent of properties of rational cuspidal curves of low degree and the Coolidge--Nagata problem.\\   

\noi All explicit information concerning the rational cuspidal curves presented in this thesis have been found using the computer programs {\em Maple} \cite{Maple} and {\em Singular} \cite{Singular}. For examples of code and calculations, see Appendix \ref{calculationsandcode}.\\

\noi The figures in this thesis are made in {\em Maple} or drawn in {\em GIMP} \cite{Gimp}. Note that the illustrations only represent the real images of the curves and that there sometimes are properties of the curves which we can not see. 

\chapter{Theoretical background}\label{TB}

Quite a lot of definitions, notations and results concerning algebraic curves are needed in order to explain what a {\em rational cuspidal curve} actually is. Not surprisingly, explaining known and finding new properties of such curves demand even more of the above. This chapter is devoted to the mentioned tasks and presents most of the theoretical background material upon which this thesis is based.

\section{Rational cuspidal curves}
Let $\mathbb{P}^2$ be the projective plane over $\mathbb{C}$, and let $(x:y:z)$ denote the coordinates of a {\em point} in $\mathbb{P}^2$. Furthermore, let $\mathbb{C}[x,y,z]$ be the ring of polynomials in $x$, $y$ and $z$ over $\mathbb{C}$. Let $F(x,y,z) \in \mathbb{C}[x,y,z]$ be a homogeneous irreducible polynomial, and let $\V(F)$ denote the zero set of $F$. Then $C = \mathcal{V}(F) \subset \Po$ is called a {\em plane algebraic curve}. By convention, when $F$ is a polynomial of degree $d$, we say that $C$ has degree $d$. Furthermore, if $F=F_1 \cdot \ldots \cdot F_{\nu}$ is a reducible polynomial and all $F_i$ are distinct, then the zero set of $F$ defines a union of curves $\V(F) = \V(F_1) \cup \ldots \cup \V(F_{\nu})$. If $F$ is a reducible polynomial and some of the factors $F_i$ are multiple, i.e., $F=F_1^{w_1} \cdot \ldots \cdot F_{\nu}^{w_{\nu}}$, then we define $\V(F)$ to be the zero set of the reduced polynomial $F=F_1 \cdot \ldots \cdot F_{\nu}$.

A curve $C$ is {\em rational} if it is birationally equivalent to $\mathbb{P}^1$ and hence admits a parametrization.

{\samepage
A point $p=(p_0:p_1:p_2)$ of $C$ is a called a {\em singularity} or, equivalently, a {\em singular point} if the partial derivatives $F_x$, $F_y$ and $F_z$ satisfy \begin{equation}
F_x(p)=F_y(p)=F_z(p)=0.
\label{sing}
\end{equation}
Otherwise, we call $p$ a {\em smooth point}. The set of singularities of a curve $C$ is usually referred to as $\mathrm{Sing}\,C$, and this is a finite set of points \cite[Cor. 3.6., pp.45--46]{Fisc:01}.}\\

\noi Given $C$, to each point $p \in \mathbb{P}^2$  we assign an integer value $m_p$, called the {\em multiplicity} of $p$ on $C$. If $p \notin C$, we define $m_p=0$. If $p \in C$, we move $p$ to $(0:0:1)$ using a linear change of coordinates. We write 
\begin{equation*}
\begin{split}
F(x,y,1)&=f(x,y)\\
&= f_{(m)}(x,y) + f_{(m+1)}(x,y) + \ldots + f_{(d)}(x,y),
\end{split}
\end{equation*}
where each $f_{(i)}(x,y)$ denotes a homogeneous polynomial in $x$ and $y$ of degree $i$. We define $m_p=m$. For a point $p$ of an irreducible algebraic curve $C$, $0 < m_p < d$, since $F(x,y,1)=f_{(d)}(x,y)$ contradictorily implies that $F$ is reducible. Additionally, it follows from the definition (\ref{sing}) that $p$ is a singularity if and only if $m_p > 1$ and that $p$ is a smooth point if and only if $m_p=1$.\\


\noi The {\it tangent} to a curve $C$ at a point $p=(p_0:p_1:p_2)$ is denoted by $T_pC$, or simply $T_p$ if there is no ambiguity. If $p$ is a smooth point, then there exists a unique tangent $T_p$ to $C$ at $p$, given by \cite[Prop. 3.6., pp.45--46]{Fisc:01} $$T_p=p_0F_x+p_1F_y+p_2F_z.$$ 

\noi If $p$ is a singularity, this definition fails. Relocating $p$ to $(0:0:1)$, we have that $$f_{(m)}(x,y)=\prod_{i=1}^m L_i(x,y),$$ where $L_i(x,y)$ are linear polynomials, not necessarily distinct. For the reduced polynomial $$f_{(m)}(x,y)=\prod_{i=1}^k L_i(x,y),$$ where the $k$, $1 \leq k \leq m$, polynomials $L_i(x,y)$ are distinct, let $T_i=\V(L_i(x,y))$. Then $\V(f_{(m)}(x,y))$ is a union of $k$ lines $T_i$ through $p$, $$\V(f_{(m)}(x,y))=\bigcup^k_{i=1} T_i. $$ The $k$ lines $T_i$ are called the {\em tangents} to $C$ at $p$ \cite[pp.41--42]{Fisc:01}. In the particular case that $k=1$ and $C$ only has one branch through $p$, $p$ is called a {\em cusp}.\\

\noi If the set of singular points of $C$ only consists of cusps, we call the curve {\em cuspidal}.

\begin{defi}[Rational cuspidal curve]
A {\em rational cuspidal curve} is a plane algebraic curve which is birational to $\mathbb{P}^1$ and is such that all its singularities are cusps.
\end{defi}

\noi Note that since all curves in this thesis are rational, we often refer to these curves as merely {\em cuspidal curves}.

\section{Invariants and conditions}
Now that we have defined a rational cuspidal curve, we add new, and further investigate the previously defined, properties of particular points on a curve.

\subsubsection{Linear change of coordinates}\label{pcc}
A {\em linear change of coordinates} in $\mathbb{P}^2$, given by a map $\tau$, will in the following be represented by an invertible $3 \times 3$ matrix $\mathcal{T} \in PGL_3(\mathbb{C})$.
$$\begin{array}{rccc}
\tau: & \mathbb{P}^2 &\longrightarrow &\mathbb{P}^2\\
&\rotatebox[origin=c]{90}{$\in$}&&\rotatebox[origin=c]{90}{$\in$}\\
&(x:y:z) &\longmapsto &(x:y:z) \cdot \mathcal{T}^{-1}.
\end{array}$$

\noi Observe that we may easily trace points under the transformation. The rows in $\mathcal{T}$, representing points in $\mathbb{P}^2$, are moved to the respective coordinate points. The first row is moved to the point $(1:0:0)$, the second row to $(0:1:0)$ and the third row to $(0:0:1)$.\\

\noi Two curves $C$ and $D$ are called {\em projectively equivalent} if there exists a linear change of coordinates such that $C$ is mapped onto $D$.

\subsubsection{Monoidal transformations}
Let $Y$ be a nonsingular surface and $p$ a point of $Y$. A {\em monoidal transformations} is the operation of {\em blowing-up}  $Y$ at $p$  \cite[p.386]{Hart:1977}. We denote this by $\sigma:\bar{Y}\longrightarrow Y$. The transformation $\sigma$ induces an isomorphism of $\bar{Y} \setminus \sigma^{-1}(p)$ onto $Y \setminus p$. The inverse image of $p$ is a curve $E$, which is isomorphic to $\mathbb{P}^1$ and is called the {\em exceptional line}.

If $C$ is a curve in $Y$, we define the {\em strict transform} $\bar{C}$ of $C$ as the closure in $\bar{Y}$ of $\sigma^{-1}(C \cap (Y \setminus p))$.

We will refer to a monoidal transformation as a blowing-up of a point, and the inverse operation will be referred to as a blowing-down of an exceptional line.

\subsubsection{Multiplicity sequence}
Let $(C,p)$ denote an irreducible analytic plane curve germ $(C,p) \subset (\mathbb{C}^2,0)$. Furthermore, let
\begin{center}
$\begin{array}{rlll}
\mathbb{C}^2 = &Y \stackrel{\sigma_1}{\longleftarrow} &Y_1 \stackrel{\sigma_2}{\longleftarrow} \ldots \stackrel{\sigma_n}{\longleftarrow} &Y_n\\
&\cup &\cup &\cup\\
(C,p)=&C \longleftarrow &C_1 \longleftarrow \ldots \longleftarrow &C_n,
\end{array}$
\end{center}
be a sequence of blowing-ups over $p$, where $C=C_0$, and $C_{i+1}$ is the strict transform of $C_i$. Let the sequence of blowing-ups be such that it resolves the singularity $p$ on $C$. Moreover, let the sequence be such that it additionally ensures that the reduced total inverse image $D=\sigma_n^{-1} \circ \ldots \circ \sigma_1^{-1}(C)$ is a simple normal crossing divisor, but $\sigma_{n-1}^{-1} \circ \ldots \circ \sigma_1^{-1}(C)$ is not. Then this sequence of blowing-ups is called the {\em minimal embedded resolution} of the cusp.\\

\begin{figure}[H]
{\includegraphics[width=0.25\textwidth]{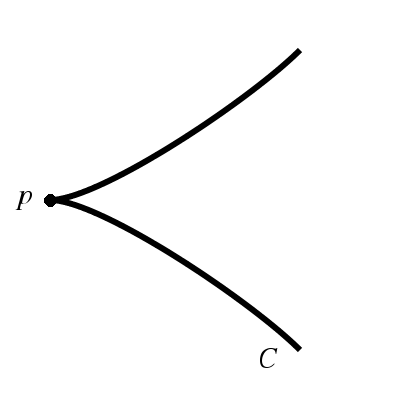}}{\includegraphics[width=0.25\textwidth]{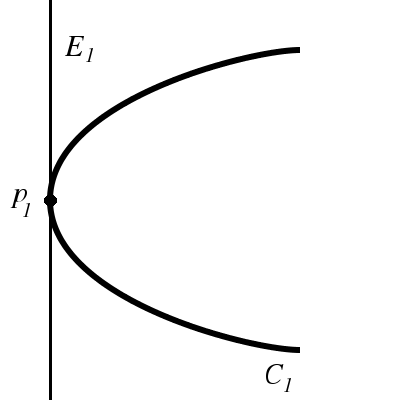}}{\includegraphics[width=0.25\textwidth]{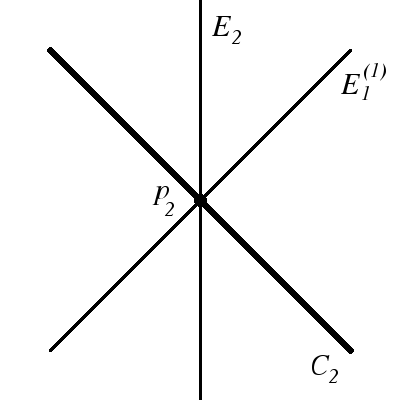}}{\includegraphics[width=0.25\textwidth]{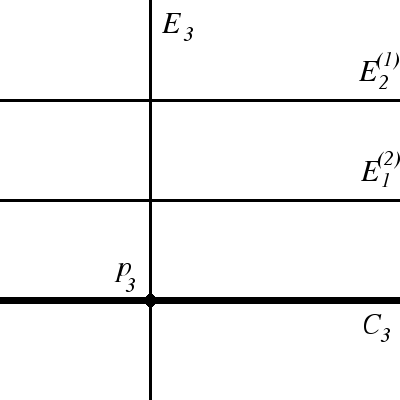}}
\caption{Minimal embedded resolution of a cusp with multiplicity sequence $(2)$.}
\label{minemres}
\end{figure}

\noi For every $i$ denote by $p_i$ the point corresponding to $p \in C$ on the curve $C_i$. The points $p_i$ are {\em infinitely near points} of $p$ on $C$, and they are referred to as the {\em strict transforms} of $p$ on $C$. Furthermore, let $m_{p.i}$ denote the multiplicity of the point $p_i \in C_i$. Then we define the {\em multiplicity sequence} of $p$ as
$$\overline{m}_p=(m_{p.0},m_{p.1},\ldots,m_{p.n}).$$ 
The index $p$ will be omitted whenever the reference point is clear from the context, and we write $m_{p.i}=m_{i}$. Note that $m_0=m_p$, which by the previous convention often is written merely $m$.

There are many important results concerning the multiplicity sequence of a point. First of all, the multiplicity sequence of a cusp $p$ has the property that \cite[p.440]{FlZa95} $$m_0 \geq m_1 \geq \ldots \geq m_n=1.$$


\noi We also have the following important result \cite[Prop. 1.2., p.440]{FlZa95}. 
\begin{prp}[On multiplicity sequences] Let $\ol{m}$ be the multiplicity sequence of a cusp $p$.
\begin{itemize}
\item[--] For each $i=1,\ldots,n$ there exists $k \geq 0$ such that $$m_{i-1}=m_i+ \ldots + m_{i+k},$$ where $$m_i=m_{i+1}=\ldots=m_{i+k-1}.$$
\item[--] The number of ending $1$'s in the multiplicity sequence equals the smallest $m_i>1$.
\end{itemize}
\label{multiseq}
\end{prp}

\noi In order to simplify notation, we introduce two conventions. First, whenever there are $k_i$ subsequent identical terms $m_i$ in the sequence, we compress the notation by writing $\ol{m}_p=(m,m_1,\ldots,(m_i)_{k_i},\ldots,1)$.  We usually also omit the ending $1$'s in the sequence. For example, if a cusp has multiplicity sequence $(4,2,2,2,1,1)$, we write merely $(4,2_3)$.\\

\noi We define the {\em delta invariant} $\delta_p$ of any point $p$ of $C$ by $$\delta_p=\sum \frac{m_{q}(m_q-1)}{2},$$ where the sum is taken over all infinitely near points $q$ lying over $p$, including $p$ \cite[Ex. 3.9.3., p.393]{Hart:1977}.

For a cusp $p$ with multiplicity sequence $\ol{m}_p$ we have \cite[p.440]{FlZa95},$$\delta_p = \sum_{i=0}^n \frac{m_i(m_i-1)}{2}.$$\medskip

\noi Let $C$ be a rational cuspidal curve with cusps $p$, $q$, $r$, \ldots. Then the curve can be described by the multiplicity sequences of the cusps. We write $[(\ol{m}_p),(\ol{m}_q),(\ol{m}_r),\ldots]$ and call this the cuspidal configuration of the curve.\\

\noi We define the {\em genus} $g$ of a curve \cite[Thm. 9.9, p.180]{Fisc:01}, $$g=\frac{(d-1)(d-2)}{2} - \sum_{p \in \mathrm{Sing}\,C} \delta_p.$$
Furthermore, a rational curve has genus $g=0$. From the above definition we derive a formula which is valid for rational cuspidal curves. 

\begin{thm}[Genus formula for a rational cuspidal curve] \label{thm:genus} 
Let $d$ be the degree of a rational cuspidal curve $C$ with singularities $p_j,\; j=1,...,s$, and let $m_{j.i}$ be the multiplicity of $p_j$ after $i$ blowing-ups. Let $n_{j}$ be the number of blowing-ups required to resolve the singularity $p_j$. Let $\delta_j$ be the delta invariant of $p_j$.  Then $$\frac{(d-1)(d-2)}{2}=\sum_{j=1}^s \delta_{j}= \sum_{j=1}^s \sum_{i=0}^{n_{j}} \frac{m_{j.i}(m_{j.i}-1)}{2}.$$
\end{thm}
\medskip

\noi The multiplicity sequence is often used to describe a cusp. Sometimes, however, it is convenient to use a different notation. In this thesis we will inconsistently refer to a cusp by either the multiplicity sequence, the Arnold classification or simply a common name. Customary notations for some of the more frequently encountered cusps are given in Table \ref{tab:notationofcusp}.\\ 
\begin{table}[H]
\centering
\setlength{\extrarowheight}{2pt}
\begin{tabular}{lcc}
\hline
{\small {\bf Common name}} & {\small {\bf Multiplicity sequence}} & {\small {\bf Arnold type}}\\
\hline
{\small Simple cusp of multiplicity 2} & $(2)$ & $A_2$\\
{\small Double cusp} & $(2_2)$ & $A_4$\\
{\small Ramphoid cusp} ($(k-1)$th type) & $(2_k)$ & $A_{2k}$ \\
{\small Simple cusp of multiplicity 3} & $(3)$ & $E_6$\\
{\small Fibonacci cusp ($k$th type)} & $(\f_k, \f_{k-1},\ldots,1,1)$\footnote{$\f_k$ is the kth Fibonacci number} & \\
\hline
\multicolumn{3}{l}{\footnotesize 1. $\f_k$ is the $k$th Fibonacci number, see Chapter \ref{mcc}.}
\end{tabular}
 \caption{What will you call a beautiful cusp?}
 \label{tab:notationofcusp}
 \end{table}

\noi The multiplicity and the multiplicity sequence serve as two very important invariants of a cusp. If two cusps have the same multiplicity sequence, then they are called {\em topologically equivalent}. This classification is, most of the time, sufficient to give a good description of a cuspidal curve. We sometimes do, however, need a finer classification of singularities. The {\em intersection multiplicity} of a cusp with its tangent appears to be an essential invariant in this context. 

\subsubsection{Intersection multiplicity}
Let $C=\V(F)$ and $D=\V(G)$ be algebraic curves which do not have any common components. If a point $p$ is such that $p \in C$ and $p \in D$, we say that $C$ and $D$ intersect at $p$. The point $p$ is called an {\em intersection point}. For an intersection point $p=(0:0:1)$, the {\em intersection multiplicity} $(C \cdot D)_p$ is defined as
\begin{equation*}
(C \cdot D)_p=\dim_{\mathbb{C}}\mathbb{C}[x,y]_{(x,y)}/(f,g),
\end{equation*}
where $f=F(x,y,1)$ and $g=G(x,y,1)$ \cite[pp.75--76]{Fulton}.\\

\noi The intersection multiplicity can be calculated directly by $$(C \cdot D)_p=\sum m^C_{p_i}m^D_{p_i},$$ where $p_i \in C_i \cap D_i$ are infinitely near points of $p$, and $m_{p_i}^C$ and $m_{p_i}^D$ denote the multiplicities of the points $p_i$ with respect to the curves $C_i$ and $D_i$ respectively.\\

\noi When working implicitly with curves, we are not able to calculate $(C \cdot D)_p$ directly. We can, however, estimate $(C \cdot D)_p$. \medskip \\
\noi First of all, we have Bezout's theorem \cite[Thm 2.7., p.31]{Fisc:01}. It provides a powerful global result on the intersection of two curves and hence an upper bound for an intersection multiplicity of two curves at an intersection point.

\begin{thm}[B\'{e}zout's theorem]\label{thm:B\'{e}zout}\label{Bez}
For plane algebraic curves $C$ and $D$ of degree $\deg C$ and $\deg D$ which do not have any common component, we have
\begin{equation*}
\sum_{p \in C \cap D} (C \cdot D)_p= \deg C \cdot \deg D.
\end{equation*}
\end{thm}

\noi In particular, for the intersection between a curve $C$ of degree $d$ and a line $L$, we have
\begin{equation*}
\sum_{p \in C \cap L} (C \cap L)_{p}= d.
\end{equation*}

\noi By B\'{e}zout's theorem, the set of intersection points of two curves $C$ and $D$ with no common component is finite. Let $p_j$, $j=1,\ldots, s$, denote the intersection points of $C$ and $D$. Then we write $$C \cdot D = (C \cdot D)_{p_1} \cdot {p_1} +  \ldots +(C \cdot D)_{p_s} \cdot p_s.$$\smallskip

\noi Second, if $L$ is a line and $p \in C \cap L$, then \cite[Prop. 3.4, p.41]{Fisc:01}\begin{equation*}\label{mintnotstrict} m_p \leq (C \cdot L)_p.\end{equation*} Furthermore, for the tangent line $T_p$, the inequality is strict, \begin{equation}\label{mint} m_p < (C \cdot T_p)_p.\end{equation} Hence, we have the inequality $$\sum_{p \in C \cap L} m_{p} \leq d.$$ Note that the inequality is strict if and only if $L$ is tangent to $C$ at one or more of the intersection points.

Moreover, if $C$ is smooth at $p$, then $(C \cdot T_p)_p \geq 2$. If $(C \cdot T_p)_p = 2$, we call $T_p$ a {\em simple tangent}. If $(C \cdot T_p)_p \geq 3$, we call $T_p$ an {\em inflectional tangent}. In the latter case we call the smooth point $p$ an {\em inflection point}. Note that we refine the definition of inflection points by calling $p$ an inflection point of type $t=(C \cdot T_p)_p - 2$. \medskip \\
\noi Third, we have a lemma linking multiplicity sequences and intersection multiplicities \cite[Lemma 1.4., p.442]{FlZa95}. For this lemma we change the notation and define the multiplicity sequence to be infinite, setting $m_{\nu}=1$ for all $\nu \geq n$. Note that in this notation a smooth point has multiplicity sequence $(1,1,\ldots)$.

\begin{lem}\label{intmult}\label{FlZa14}
Let $(C,p)$ be an irreducible germ of a curve, and let $p$ have multiplicity sequence $\ol{m}_p$. Then there exists a germ of a smooth curve $(\Gamma,p)$ through $p$ with $(\Gamma \cdot C)_p = k$ if and only if $k$ satisfies the condition
$$k = m_0+m_1+ \ldots +m_a\; \text{ for some } a>0\; \text{ with } m_0=m_1= \ldots =m_{a-1}.$$\smallskip
\end{lem}

\noi All the above results can be used to estimate $(C\cdot T_p)_p$ for a cusp. We will frequently use the letter $r$ for this invariant, i.e., $r_p=(C \cdot T_p)_p$. B\'{e}zout's theorem (\ref{thm:B\'{e}zout}) provides an upper bound for $(C\cdot T_p)_p$, while Lemma \ref{intmult} combined with inequality (\ref{mint}) provides a lower bound. \begin{equation} m_0+m_1 \leq (C \cdot T_p)_p \leq d.\label{upperandlower} \end{equation} Lemma \ref{intmult} additionally provides information about the possible values between the upper and lower bound.
\begin{equation*}
\begin{split}
(C \cdot T_p)_p =& \sum_{i=0}^a m_i\\
=&\;a \cdot m_0 + m_a \text { for some } a \geq 1.
\end{split}
\label{tangentintersection}
\end{equation*}

\subsubsection{Puiseux parametrization}\label{puiseux}
In order to investigate a point on a curve in more detail, we will occasionally parametrize the curve locally. Since smooth points and cusps are unibranched, each point on a cuspidal curve can be given a local parametrization by power series, a {\em Puiseux parametrization}. Let $(C,p)$ be the germ of a cuspidal curve $C$ at the point $p=(0:0:1)$, and let $\V(y)$ be the tangent to $C$ at $p$. With $m=m_p \geq 1$ and $r=(C \cdot T_p)_p > m$, the germ $(C,p)$ can be parametrized by \cite[Cor. 7.7, p.135]{Fisc:01}
\begin{equation}
\begin{split}
&x=t^m,\\
&y=c_rt^r + \ldots,\\
&z=1,
\label{pp}
\end{split}
\end{equation}
where $\ldots$ denotes higher powers of $t$, the coefficients of $t^i$ in the power series expansion of $y$ are $c_i \in \mathbb{C}$, and $c_r \neq 0$.


Observe that, in this form, the Puiseux parametrization reveals both the multiplicity of $p$ and the intersection multiplicity of the curve and the tangent at the point. So far, the Puiseux parametrization seems like a straightforward matter. There are, however, some subtleties involved. 

\begin{ex}\label{normform}
Cusps of type $A_{2k}$ can topologically be represented by the normal form \cite[Table 2.2., p.219]{Piene99} $$y^2+x^{2k+1}.$$ The normal form implies the parametrization $$(t^2:t^{2k+1}:1).$$ 

\noi We frequently need to describe the $A_{2k}$-cusps in more detail. For example, if the curve has degree $d=4$, then the tangent intersects the curve at the $A_{2k}$-cusp with multiplicity $4$ for $k>1$. The cusp can then be parametrized by $$(t^2:c_4t^4+\text{(even powers of t)}+ c_{2k+1}t^{2k+1}+\ldots:1),\;c_4,c_{2k+1} \neq 0.$$ 

\begin{table}[H]
\centering
\setlength{\extrarowheight}{2pt}
\begin{tabular}{cl}
\hline
{\bf Type} & {\bf Puiseux parametrization}\\
\hline
$A_2$ & $(t^2:c_3t^3+c_4t^4+\ldots:1)$\\
$A_4$ & $(t^2:c_4t^4+c_5t^5+\ldots:1)$\\
$A_6$ & $(t^2:c_4t^4+c_6t^6+c_7t^7+\ldots:1)$\\
\hline
\end{tabular}
\caption{Puiseux parametrization for cusps of type $A_{2k},\; k=1,2,3$, on a curve of degree $d=4$.}
\label{tab:pui}
\end{table}

\noi If the curve has degree $d=5$, the picture gets even more complicated. For example, if $p$ is an $A_4$-cusp of a quintic curve, then the tangent may intersect the curve with multiplicity $4$ or $5$. The value of $r$ must be determined by other methods.\medskip
\end{ex}

\noi The example reveals that the multiplicity sequence does not determine the full complexity of the Puiseux parametrization. We are, however, able to to calculate the multiplicity sequence from the Puiseux parametrization \cite[Thm. 12., p.516]{Brieskorn}\cite[p.234]{MatsuokaSakai}.

Given $(C,p)$ and a Puiseux parametrization on the form (\ref{pp}), let the {\em characteristic terms} of the Puiseux parametrization be the terms $c_{\beta_{\je}}t^{\beta_{\je}}$ of the power series expansion of $y$ defined by
\begin{itemize}
\item[--] $m>\gcd(m,\beta_1)>\ldots>\gcd(m,\beta_1,\ldots,\beta_g)=1,$
\item[--] $c_{\beta_{\je}} \neq 0 \text{ for } \je=1,\ldots,g,$
\item[--] if $\beta_1,\ldots,\beta_{\je-1}$ have been defined and if $\gcd(m,\beta_1,\ldots,\beta_{\je-1})>1$, then $\beta_{\je}$ is the smallest $\beta$ such that $c_{\beta_{\je}}\neq 0$ and $\gcd(m,\beta_1,\ldots,\beta_{\je-1})>\gcd(m,\beta_1,\ldots,\beta_{\je-1},\beta_{\je}).$
\end{itemize}


\noi Let $(D,q)$ be a germ given by the Puiseux parametrization of $(C,p)$ in such a way that the power series expansion of $y$ only consists of characteristic terms,
\begin{equation*}
\begin{split}
x&=t^m\\
y&=c_{\beta_1}t^{\beta_1}+c_{\beta_2}t^{\beta_2}+ \ldots +c_{\beta_g}t^{\beta_g}\\
z&=1.
\end{split}
\end{equation*}

\noi Although Example \ref{normform} reveals that there potentially are many differences between $(C,p)$ and $(D,q)$, the point $p$ of the germ $(C,p)$ has the same multiplicity sequence as the point $q$ of the germ $(D,q)$. Furthermore, we can calculate the multiplicity sequence.

\begin{thm}\label{thm:mp}
Let $q$ be a point of an irreducible germ $(D,q)$ where the Puiseux parametrization only consists of characteristic terms. Then the multiplicity sequence of $q$ is determined by a chain of Euclidian algorithms. Let $\gamma_{\je}=\beta_{\je}-\beta_{\je-1},$ and $\beta_0=0$. For each $\je$, let
$$\begin{array}{rclc}
\gamma_{\je}&=&a_{\je,1}m_{\je,1}+m_{\je,2} &\qquad (0<m_{\je,2}<m_{\je,1})\\
m_{\je,1}&=&a_{\je,2}m_{\je,2}+m_{\je,3} &\qquad (0<m_{\je,3}<m_{\je,2})\\
&\ldots& &\qquad \ldots\\
m_{\je,q_{\je}-1}&=&a_{\je,q_{\je}}m_{\je,q_{\je}},&
\end{array}$$
where $m_{1,1}=m$, $m_{\je+1,1}=m_{\je,q_{\je}}$, and $m_{g,q_g}=1$. 
The multiplicity sequence of the point $q$ on $D$ is given by 
$$\ol{m}_q = (\overbrace{m_{1,1},\ldots,m_{1,1}}^{a_{1,1}},\ldots,\overbrace{m_{\je,k},\ldots,m_{\je,k}}^{a_{\je,k}},\ldots,\overbrace{1,\ldots,1}^{a_{g,q_g}}).$$
\end{thm}


\subsubsection{Properties of the blowing-up process}\label{properelm}
The blowing-up process has certain elementary properties that will be invaluable in the later study of curves. \medskip \\
\noi First of all, we have the {\em self-intersection} of the exceptional line $E$ on $\bar{Y}$. We will use, but not define, self-intersection here, see Hartshorne \cite[pp.360--361]{Hart:1977} for a formal definition. For any monoidal transformation we have that the self-intersection of $E$ on $\bar{Y}$ is $E^2=-1$ \cite[p.386]{Hart:1977}.\medskip \\
\noi Second, we have the following important lemma from Flenner and Zaidenberg \cite[Lemma 1.3., pp.440--441]{FlZa95}.
\begin{lem}
Let $\ol{m}_p$ be the multiplicity sequence of a point $p$ on a curve $C$ as defined prior to Lemma \ref{FlZa14}. Let $\sigma_i$ be a sequence of blowing-ups and let $Y_i$ be the corresponding surfaces. Denote by $E_i^{(k)}$ the strict transform of the exceptional divisor $E_i$ of $\sigma_i$ at the surface $Y_{i+k}$. Then
\begin{align*}
(E_i \cdot C_i)_{p_i}&=m_{i-1},\\
(E_i^{(k)} \cdot C_{i+k})_{p_{i+k}}&=\max \{0,m_{i-1}-m_i-\ldots-m_{i+k-1}\},\;k>0,\\
(E_i^{(1)} \cdot C_{i+1})_{p_{i+1}}&=m_{i-1}-m_i.
\end{align*} 
\end{lem} 

\noi Third, note that we may calculate intersection multiplicities of strict transforms of curves. Since $$(C \cdot D)_p=\sum m^C_{p_i}m^D_{p_i},$$ for points $p_i \in C_i \cap D_i$, we see that for a fixed $k\geq 0$,
\begin{align*}
\centering
(C_k \cdot D_k)_{p_k}&=\sum_{i \geq k} m^{C}_{p_i}m^{D}_{p_i}\\
&=(C \cdot D)_p-\sum_{i<k}m^C_{p_i}m^D_{p_i}.
\end{align*}

\noi We will frequently use the fact that $$(C_1 \cdot D_1)_{p_1}=(C \cdot D)_p-m^C_{p}m^D_{p}.$$

\section{Derived curves}
There are a number of associated curves which are useful in the analysis of a curve $C$. 

\subsubsection{The polar curves}
The {\em polar} $P_pC$ of a curve $C$ with respect to a point $p=(p_0:p_1:p_2) \in \mathbb{P}^2$ is defined as $$P_pC=\mathcal{V}(p_0 F_x+p_1 F_y+p_2 F_z).$$ This curve has degree $d-1$. 

The points in the intersection $P_pC \cap C$ are the points $p_j$, $j=1,\ldots,s$, for which the tangents $T_{p_j}$ to $C$ at $p_j$ go through $p$, and additionally the singularities of $C$ \cite[Thm. 4.3., p.64]{Fisc:01}.

\subsubsection{The dual curve}
The dual space $\mathbb{P}^{2*}$ consists of all lines in $\mathbb{P}^2$. Since smooth points of a curve $C=\mathcal{V}(f)$ have a unique tangent, we define the rational map $\zeta$,
$$\begin{array}{cccc}
\zeta: & C \setminus \mathrm{Sing}\,C \subset \mathbb{P}^2 &\longrightarrow &\mathrm{Im} (\zeta) \subset \mathbb{P}^{2*}\\
& \rotatebox[origin=c]{90}{$\in$}&&\rotatebox[origin=c]{90}{$\in$}\\
&p=(p_0:p_1:p_2) &\longmapsto &(F_x(p):F_y(p):F_z(p))
\end{array}$$

\noi We define the {\em dual curve} $C^*$ as the closure of $\mathrm{Im}(\zeta)$,
$$C^*=\mathrm{Cl}(\mathrm{Im}(\zeta)).$$ Furthermore, $C$ and $C^*$ has the same genus \cite[p.179]{Fisc:01}.\\

\noi We may describe the dual germ $(C^*,p^*)$ of a germ $(C,p)$.
Let $(C,p)$ be given by its Puiseux parametrization,
$$(C,p)=(t^m: c_rt^r+c_{\alpha}t^{\alpha}+ \ldots:1), \qquad \; c_i \neq 0, \; \alpha > r > m.$$
Then $(C^*,p^*)$ can be found by calculating the minors of the matrix \cite[pp.73--94]{Fisc:01},
$$\begin{bmatrix}x(t)&y(t)&1 \\ \nom x'(t)&y'(t)&0 \end{bmatrix}.$$\medskip
$$(C^*,p^*)=(a^*t^{r-m} + \ldots:1:c_r^*t^r+c^*_{\alpha}t^{\alpha}+ \ldots).$$
We have that $$a^*=-\frac{c_rr}{m}, \qquad c_r^*=c_r\left(\frac{r}{m}-1\right), \qquad c_{\alpha}^*=c_{\alpha}\left(\frac{\alpha}{m}-1\right).$$ Since $c_i \neq 0$ and $\alpha>r>m$, the constants $a^*,\; c_r^*,\;c_{\alpha}^* \neq 0$. As a consequence of the calculation, the power series $c_r^*t^r+c^*_{\alpha}t^{\alpha}+ \ldots$ contains precisely the same powers of $t$ as the power series $c_rt^r+c_{\alpha}t^{\alpha}+ \ldots$. 

Using properties of the Puiseux parametrization, we may determine important invariants, like the multiplicity sequence, of the dual point $p^*$ on $C^*$. In particular, observe that we can find the multiplicity $m^*$ of the dual point $p^*$ on $C^*$,
$$m^*=r-m.$$    

\noi Additionally, a classical Plücker formula gives the degree $d^*$ of the dual curve $C^*$ \cite[p.316]{Fen99b}.
\begin{thm}\label{plucker} Let $C$ be a curve of genus $g$ and degree $d$ having $j$ singularities $p_j$ with multiplicities $m_{p_j}=m_j$. Let $b_j$ denote the number of branches of the curve at $p_j$. Then the degree $d^*$ of the dual curve is given by
$$d^*=2d+2g-2-\sum_{p_j \in \mathrm{Sing}\,C}(m_{j}-b_{j}).$$
\end{thm}

\begin{cor}
For rational cuspidal curves we have
$$d^*=2d-2-\sum_{p_j \in \mathrm{Sing}\,C}(m_{j}-1).$$
\end{cor}

\subsubsection{The Hessian curve}


Let $\mathcal{H}$ be the matrix given by
$$\mathcal{H}=\begin{bmatrix} F_{xx} & F_{xy} & F_{xz} \\ \nom F_{yx} & F_{yy} & F_{yz} \\ \nom F_{zx} & F_{zy} & F_{zz}\end{bmatrix},$$ 
where $F_{ij}$ denote the double derivatives of $F$ with respect to $i$ and $j$ for $i,j \in \{x,y,z\}$.\\

\noi Define a polynomial $H_F$, $$H_F=\det \mathcal{H}.$$  Then the {\em Hessian curve}, $H_C$ of degree $3(d-2)$, is given by $$H_C=\mathcal{V}(H_F).\label{Hessian}$$ 
By B\'{e}zout's theorem, \begin{equation}\sum_{p \in C \cap H_C} (C \cdot H_C)_p = 3d(d-2).\label{HESS}\end{equation} 
Moreover, $H_C$ and $C$ intersect at the singular points and the inflection points of $C$ \cite[p.67]{Fisc:01}.\\

\noi We have an interesting formula relating several invariants regarding the cuspidal configuration of a curve $C$ to the total intersection number between $C$ and its Hessian curve $H_C$. The below formula is given for rational cuspidal curves, but a similar result is valid for more general curves ~\cite[Thm. 2., pp.586--597]{Brieskorn}.

\begin{thm}[Inflection point formula]\label{Ghessian}
Let $C$ be a rational cuspidal curve. Let $s$ be the number of inflection points on $C$, counted such that an inflection point  of type $t$ accounts for $t$ inflection points. Let $p_j$ be the cusps of $C$ with multiplicity sequences $\ol{m}_j$, delta invariants $\delta_j$ and tangent intersection multiplicities $r_j$ at $p_j$. Let $m_j^*$ denote the multiplicities of the dual points $p_j^*$ on the dual curve $C^*$. Then the number of inflection points, counted properly, is given by
\begin{equation*}
\begin{split}
s=&\;3d(d-2)-6 \sum_{p_j \in \mathrm{Sing}\,C} \delta_{j}- \sum_{p_j \in \mathrm{Sing}\,C} (2m_{j} +m_{j}^*-3)\\
=&\;3d(d-2)-6 \sum_{p_j \in \mathrm{Sing}\,C} \delta_{j}- \sum_{p_j \in \mathrm{Sing}\,C} (m_{j} +r_{j}-3).
\end{split}
\end{equation*}
\end{thm}

\noi Using a few identities, we can rewrite this formula. For an inflection point $q$, we have that $m_q=1$, which means that $\delta_q=0$. Additionally, the type $t$ is a function of $m_q$ and $m_q^*$,
\begin{equation*}
\begin{split}
t&=(C\cdot T_q)_q-2\\
&=m_q+m_q^*-2\\
&=2m_q+m_q^*-3.
\end{split}
\end{equation*}
Moreover, if $q_i$ are the inflection points of $C$, then $s=\sum t_{i}$.\\

\noi We substitute for $s$ and use identity (\ref{HESS}) in the inflection point formula. Then we obtain the following corollary.

\pagebreak
\begin{cor}\label{CHC}
Let $C$ be a rational cuspidal curve. Let $p_j$ denote the set of both inflection points and cusps on $C$. Let $\ol{m}_{j}$ be their respective multiplicity sequences, let $r_j=(C \cdot T_{p_j})_{p_j}$, and let $\delta_{j}$ be the delta invariant of the points. Let $m_j^*$ denote the multiplicities of the dual points $p_j^*$ on the dual curve $C^*$. Then
\begin{equation*}
\begin{split}
\sum_{p_j \in C \cap H_C} (C \cdot H_C)_{p_j} &= \sum_{p_j \in C \cap H_C} (6\delta_{j}+ 2 m_{j} +m_{j}^*-3)\\
&= \sum_{p_j \in C \cap H_C} (6\delta_{j}+  m_{j} +r_{j}-3).
\end{split}
\end{equation*}
\end{cor}

\section{Other useful results}

\subsubsection{Euler's identity}
There is a fundamental dependency between a homogeneous polynomial $F$ and its partial derivatives \cite[p.45]{Fisc:01}.

\begin{thm}[Euler's identity]\label{thm:Euler}
If $F\in \mathbb{C}[x,y,z]$ is homogeneous of degree $d$, then
$$xF_x+yF_y+zF_z=d \cdot F.$$
\end{thm}

\subsubsection{The ramification condition}
We have another condition on the multiplicities of points on a rational cuspidal curve $C$, which is based on the Riemann--Hurwitz formula \cite[Lemma 3.1., p.446]{FlZa95}.
\begin{lem}[From Riemann--Hurwitz]  \label{rhof}
Let $C \subset \mathbb{P}^2$ be a rational cuspidal curve of degree d with a cusp $p \in C$ of multiplicity $m_p$ with multiplicity sequence $\ol{m}_p=(m_{p},m_{p.1},...,m_{p.n})$. Then the rational projection map $\pi_p:C \longrightarrow \mathbb{P}^1$ from $p$ has at most $2(d-m_p-1)$ ramification points. Furthermore, if $p_1,...,p_s$ are the other cusps of C and $m_{p_j}=m_j$, then
$$\sum_{j=1}^s (m_j-1)+(m_{p.1}-1) \leq 2(d-m_p-1).$$
\end{lem}

\subsubsection{On the maximal multiplicity}
Let $C$ be a rational cuspidal curve with cusps $p_j, \; j=1,\ldots,s$. Let $m_{p_j}$ denote the multiplicities of the cusps. Let $\mu$ denote the largest multiplicity of any cusp on the curve, $$\mu=\mathrm{max}_{p_j}\{m_{p_j}\}.$$\smallskip
For every rational cuspidal curve there has to be at least one cusp with a multiplicity that is quite large \cite[Thm., p.233]{MatsuokaSakai}.

\begin{thm}[Matsuoka--Sakai]
Let $C$ be a rational cuspidal plane curve of degree $d$. Let $\mu$ denote the maximum of the multiplicities of the cusps. Then 
\begin{equation*}
d < 3 \mu.
\label{MatsuokaSakai}
\end{equation*}
\end{thm}\medskip

\noi For $\mu \geq 9$ we have a better estimate \cite[Thm. A., p.657]{Ore}.
\begin{thm}[Orevkov]
Let $C$ be a rational cuspidal plane curve of degree $d$. Let $\alpha=\frac{3+\sqrt{5}}{2}$. Then
\begin{equation*}
d < \alpha (\mu+1)+\textstyle{\frac{1}{\sqrt{5}}}.
\label{MatsuokaSakaiulikhet}
\end{equation*}
\end{thm}

\section{Getting an overview}
The theoretical background in this chapter provides powerful tools for the study of rational cuspidal curves. In the next chapters we will explore and apply this theory to cuspidal curves of low degree. Before we go on with this analysis, we will give an overview of the invariants directly involved in the study and description of a particular rational cuspidal curve.

Starting out with either a parametrization or a homogeneous defining polynomial, we may investigate a rational cuspidal curve in depth. The first thing we are interested in is finding the number cusps of the curve. Next we want to study each cusp in detail. We first find its multiplicity and its multiplicity sequence, which gives us the cuspidal configuration of the curve. We then find the tangent of each cusp and the intersection multiplicity of the tangent and the curve at the point. This enables us to distinguish cusps with identical multiplicity sequences. 

The above gives us the necessary overview of the cusps of a cuspidal curve. There is, however, more to a rational cuspidal curve than its cusps. For example, two curves with identical cuspidal configurations are not necessarily projectively equivalent. They may have different number and types of inflection points. In some of the descriptions of rational cuspidal curves in this thesis, we will therefore include a discussion of the inflection points of the curve.

Since we have a restriction on the total intersection multiplicity, and because we discuss the local intersection multiplicity of a curve and its Hessian curve in Section \ref{ICHC}, we also provide the intersection multiplicity of the Hessian curve and the curve at cusps and inflection points when we present the curves.


\chapter{Rational cuspidal cubics and quartics}\label{rccq}
In this chapter we will use the results of Chapter \ref{TB} to obtain a list of possible rational cuspidal cubics and quartics. Furthermore, in order to get an overview of the curves, we briefly describe all cuspidal curves of mentioned degrees up to projective equivalence.

\section{Rational cuspidal cubics}
Let $C$ be a rational cuspidal cubic. Substituting $d=3$ in Theorem \ref{thm:genus} gives $$\frac{(3-1)(3-2)}{2}=1=\sum_{p\in \mathrm{Sing}\,C}\sum_{i=0}^{n_p}\frac{m_i(m_{i}-1)}{2}.$$ We see from this formula that $C$ can only have one cusp. In particular, the cusp must have multiplicity sequence $\ol{m}=(2)$. Hence, we have only one possible cuspidal configuration for a cubic curve, $[(2)]$.


\subsubsection{The cuspidal cubic -- $[(2)]$}
The cuspidal cubic can be given by the parametrization $$(s^3:st^2:t^3).$$ An illustration of the cuspidal cubic and a brief summary of its properties are given in Table \ref{fig:cuspidalcubic}.\\ 

\begin{table}[htb!]
\centering
\setlength{\extrarowheight}{2pt}
	\begin{tabular}{c}{\includegraphics[width=0.3\textwidth]{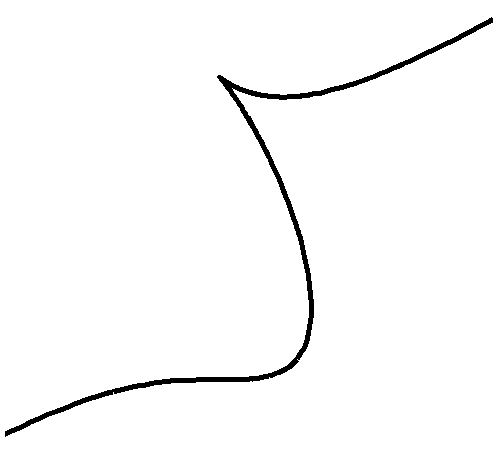}} \end{tabular}
	{\begin{tabular}{ccc}
	\multicolumn{3}{c}{} \\ 
		\multicolumn{3}{c}{$(s^3:st^2:t^3)$} \\ 
		\multicolumn{3}{c}{} \\ 
	\hline
	\multicolumn{3}{c}{{\bf \# Cusps = 	1 }}\\
	{\bf Cusp $p_j$} & {\bf $(C \cdot T_{p_j})_{p_j}$}& {\bf $(C \cdot H_C)_{p_j}$ }\\
	\hline 
	$(2)$& 3 & 8\\
	\hline
		&&\\
		\hline
			\multicolumn{3}{c}{{\bf \# Inflection points =	1 }}\\
{\bf Inflection point $q_j$} & {\bf $(C \cdot T_{q_j})_{q_j}$}& {\bf $(C \cdot H_C)_{q_j}$ }\\	
	\hline
	$q_1$ & 3 & 1 \\
	\hline
	\end{tabular}}
	\caption{Cuspidal cubic -- $[(2)]$}
	\label{fig:cuspidalcubic}
\end{table}

\noi Using {\em Singular} and the code given in Appendix \ref{calculationsandcode}, we find that the defining polynomial of this curve is $F=y^3-xz^2$. The partial derivatives of $F$ vanish at $p=(1:0:0)$, hence this point is the cusp. $C$ has tangent $T_p=\V(z)$ at $p$, and $T_p$ intersects $C$ at $p$ with multiplicity $(C \cdot T_p)_p=3$. 

The Hessian curve $H_C$ is given by $H_F=24yz^2$. Since $(0:0:1)$ is a smooth point and $$H_C \cap C=\{(1:0:0),(0:0:1)\},$$ $C$ has an inflection point at $q=(0:0:1)$. Indeed, we have the tangent at $q$ given by $T_q=\V(x)$, and this line intersects $C$ at $q$ with multiplicity $(C \cdot T_q)_q=3$. 

The parametrization of $C$ can be studied locally. Setting $s=1$, we find the germ of the curve at the cusp $p$, $$(C,p)=(1:t^2:t^3).$$ Similarly, setting $t=1$, we find the germ of the curve at the inflection point $q$, $$(C,q)=(s^3:s:1).$$

\section{Rational cuspidal quartics}
Let $C$ be a rational cuspidal quartic. Since $d=4 > m_p \geq 2$, any cusp on $C$ must have multiplicity $m=3$ or $m=2$. Additionally, substituting $d=4$ in Theorem \ref{thm:genus} gives \begin{equation}\frac{(4-1)(4-2)}{2}=3=\sum_{p\in \mathrm{Sing}\,C}\sum_{i=0}^{n_p}\frac{m_i(m_{i}-1)}{2}.\label{4:genus}\end{equation}\medskip 

\noi Assume that $C$ has a cusp with $m=3$. By (\ref{4:genus}), $C$ can not have any other cusps. Moreover, the cusp must have multiplicity sequence $(3)$. \medskip \\
\noi Assume that $C$ has a cusp with $m=2$. By (\ref{4:genus}), $C$ can not have more than three cusps. If there are three cusps on $C$, each cusp must have multiplicity sequence $(2)$. If there are two cusps on $C$, then one cusp must have multiplicity sequence $(2_2)$, while the other cusp must have multiplicity sequence $(2)$. If there is just one cusp on $C$ and $m=2$, then this cusp must have multiplicity sequence $(2_3)$.\\

\noi For each of the possible cuspidal configurations there exists at least one quartic curve, up to projective equivalence. The classification of rational cuspidal quartic curves up to projective equivalence is given by Namba in \cite[pp.135,146]{Namba}. The cuspidal quartic curves with maximal multiplicity $m=2$ are unique up to projective equivalence. For the curve with a cusp with multiplicity $m=3$, however, there are two possibilities. An overview of all existing rational cuspidal quartic curves up to projective equivalence is given in Table \ref{tab:degree4}. 

\begin{table}[htb]
  \renewcommand\thesubtable{}
  \setlength{\extrarowheight}{2pt}
\centering
	{\begin{tabular}{ccll}
	\hline
	{\bf \# Cusps}&	{\bf Curve} &{\bf Cuspidal configuration} & {\bf \# Curves}\\
	\hline 
	\multirow{1}{14mm}{3}&$C_1$&$\qquad(2),(2),(2)$&1\\
	\hline
	\multirow{1}{14mm}{2}&$C_2$&$\qquad (2_2),(2)$&1\\
		\hline
		\multirow{2}{14mm}{1}&$C_3$&$\qquad (2_3)$&1\\
	&$C_4$&$\qquad (3)$&2 -- AB\\
	\hline
	\end{tabular}}
	\caption {Rational cuspidal quartic curves.}
	\label{tab:degree4}
	\end{table}


\subsubsection{$C_1$ -- Tricuspidal quartic -- $[(2),(2),(2)]$}

\begin{table}[H]
  \renewcommand\thesubtable{}
  \setlength{\extrarowheight}{2pt}
\centering
	\begin{tabular}{c}{\includegraphics[width=0.3\textwidth]{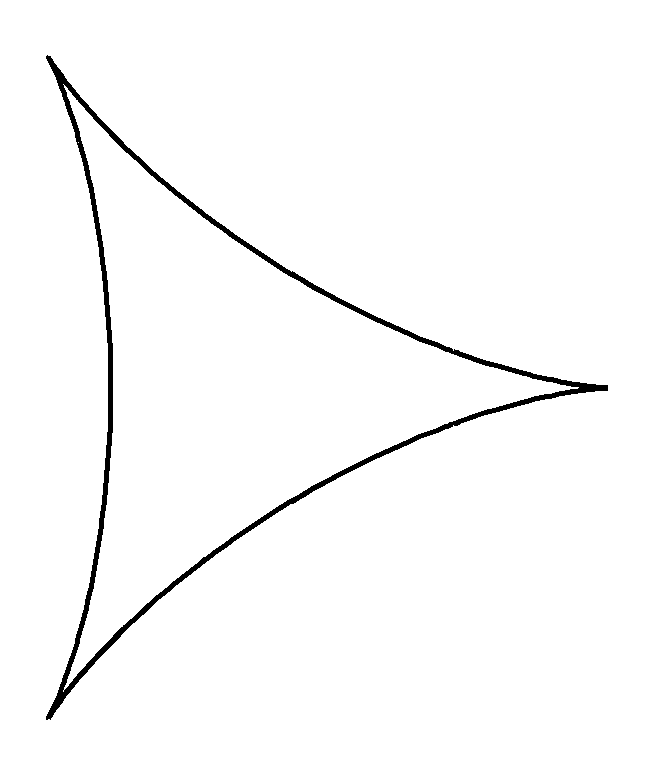}}\end{tabular}
	{\begin{tabular}{ccc}
			\multicolumn{3}{c}{} \\ 
	\multicolumn{3}{c}{$(s^3t-\textstyle{\frac{1}{2}}s^4:s^2t^2:t^4-2st^3)$}\\
			\multicolumn{3}{c}{} \\ 
	\hline
	\multicolumn{3}{c}{{\bf \# Cusps = 	3 }}\\
{\bf Cusp $p_j$} & {\bf $(C \cdot T_{p_j})_{p_j}$}& {\bf $(C \cdot H_C)_{p_j}$ }\\
	\hline 
	(2)& 3 & 8\\
	(2)& 3 & 8\\
	(2)& 3 & 8\\
	\hline
		&&\\
		\hline
			\multicolumn{3}{c}{{\bf \# Inflection points =	0}}\\
	\hline
	\end{tabular}}
	\end{table}

\vspace{25mm}
\subsubsection{$C_2$ -- Bicuspidal quartic -- $[(2_2),(2)]$}
\begin{table}[H]
  \renewcommand\thesubtable{}
  \setlength{\extrarowheight}{2pt}
\centering
	\begin{tabular}{c}{\includegraphics[width=0.3\textwidth]{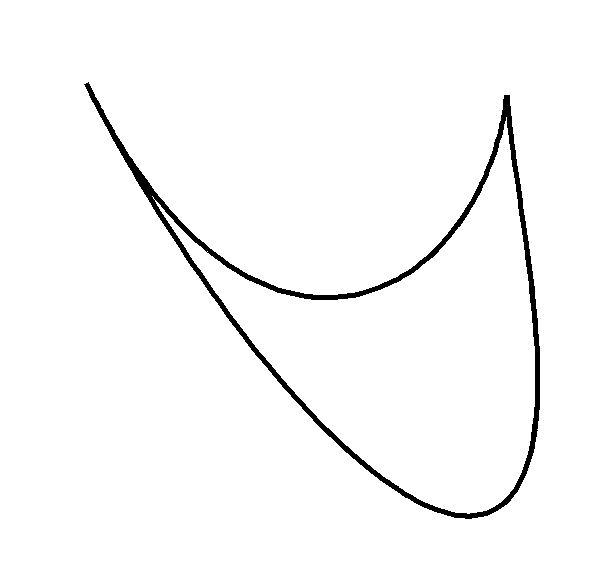}}\end{tabular}
	{\begin{tabular}{ccc}
				\multicolumn{3}{c}{} \\ 
	\multicolumn{3}{c}{$(s^4+s^3t:s^2t^2:t^4)$}\\
			\multicolumn{3}{c}{} \\ 
	\hline
	\multicolumn{3}{c}{{\bf \# Cusps = 	2 }}\\
	{\bf Cusp $p_j$} & {\bf $(C \cdot T_{p_j})_{p_j}$}& {\bf $(C \cdot H_C)_{p_j}$ }\\
	\hline 
	$(2_2)$& 4 & 15\\
	$(2)$& 3 & 8\\
	\hline
		&&\\
		\hline
			\multicolumn{3}{c}{{\bf \# Inflection points =	1 }}\\
{\bf Inflection point $q_j$} & {\bf $(C \cdot T_{q_j})_{q_j}$}& {\bf $(C \cdot H_C)_{q_j}$ }\\	
	\hline
	$q_1$ & 3 & 1 \\
	\hline
	\end{tabular}}
	\end{table}

\vspace{25mm}
\subsubsection{$C_3$ -- Unicuspidal ramphoid quartic -- $[(2_3)]$}
\begin{table}[H]
  \renewcommand\thesubtable{}
  \setlength{\extrarowheight}{2pt}
\centering
	\begin{tabular}{c}{\includegraphics[width=0.3\textwidth]{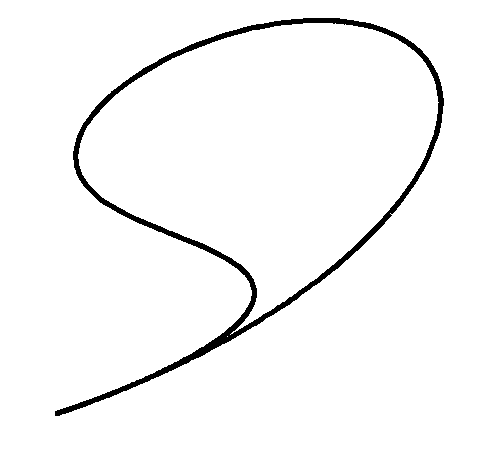}}\end{tabular}
	{\begin{tabular}{ccc}
		\multicolumn{3}{c}{} \\ 
	\multicolumn{3}{c}{$(s^4+st^3:s^2t^2:t^4)$}\\
			\multicolumn{3}{c}{} \\
	\hline
	\multicolumn{3}{c}{{\bf \# Cusps = 	1 }}\\
	{\bf Cusp $p_j$} & {\bf $(C \cdot T_{p_j})_{p_j}$}& {\bf $(C \cdot H_C)_{p_j}$ }\\
	\hline 
	$(2_3)$& 4 & 21\\
	\hline
		&&\\
		\hline
			\multicolumn{3}{c}{{\bf \# Inflection points =	3 }}\\
{\bf Inflection point $q_j$} & {\bf $(C \cdot T_{q_j})_{q_j}$}& {\bf $(C \cdot H_C)_{q_j}$ }\\		
	\hline
	$q_1$ & 3 & 1 \\
	$q_2$ & 3 & 1 \\
	$q_3$ & 3 & 1 \\
	\hline
	\end{tabular}}
	\end{table}
	
	\vspace{25mm}
	\subsubsection{$C_{4A}$ -- Ovoid quartic A -- $[(3)]$}
	\vspace{15mm}
	\begin{table}[H]
  \renewcommand\thesubtable{}
  \setlength{\extrarowheight}{2pt}
\centering
	\begin{tabular}{c}{\includegraphics[width=0.3\textwidth]{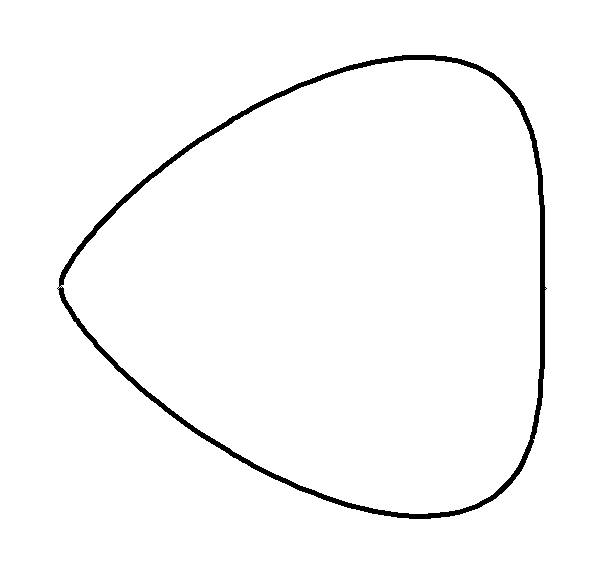}}\end{tabular}
	{\begin{tabular}{ccc}
		\multicolumn{3}{c}{} \\ 
	\multicolumn{3}{c}{$(s^4:st^3:t^4)$}\\
			\multicolumn{3}{c}{} \\
	\hline
	\multicolumn{3}{c}{{\bf \# Cusps = 	1 }}\\
	{\bf Cusp $p_j$} & {\bf $(C \cdot T_{p_j})_{p_j}$}& {\bf $(C \cdot H_C)_{p_j}$ }\\
	\hline 
	$(3)$& 4 & 22\\
	\hline
		&&\\
		\hline
			\multicolumn{3}{c}{{\bf \# Inflection points =	1 }}\\
{\bf Inflection point $q_j$} & {\bf $(C \cdot T_{q_j})_{q_j}$}& {\bf $(C \cdot H_C)_{q_j}$ }\\	
	\hline
	$q_1$ & 4 & 2 \\
	\hline
	\end{tabular}}
	\end{table}
	
\vspace{35mm}	
	\subsubsection{$C_{4B}$ -- Ovoid quartic B -- $[(3)]$}
	\begin{table}[b!]
  \renewcommand\thesubtable{}
  \setlength{\extrarowheight}{2pt}
\centering
	\begin{tabular}{c}{\includegraphics[width=0.3\textwidth]{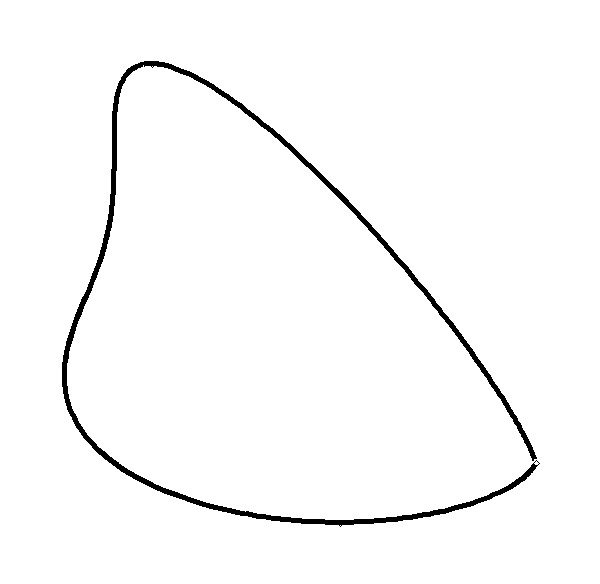}}\end{tabular}
	{\begin{tabular}{ccc}
		\multicolumn{3}{c}{} \\ 
	\multicolumn{3}{c}{$(s^3t-s^4:st^3:t^4)$}\\
			\multicolumn{3}{c}{} \\
	\hline
	\multicolumn{3}{c}{{\bf \# Cusps = 	1 }}\\
	{\bf Cusp ${p_j}$} & ${(C \cdot T_{p_j})_{p_j}}$& ${(C \cdot H_C)_{p_j} }$\\
	\hline 
	$(3)$& 4 & 22\\
	\hline
		&&\\
		\hline
			\multicolumn{3}{c}{{\bf \# Inflection points =	2 }}\\
{\bf Inflection point $q_j$} & {\bf $(C \cdot T_{q_j})_{q_j}$}& {\bf $(C \cdot H_C)_{q_j}$ }\\	
	\hline
$q_1$ & 3 & 1 \\
$q_2$ & 3 & 1 \\
	\hline
	\end{tabular}}
	\end{table}


\chapter{Projections}\label{proj}
Projection is a method by which it is possible to construct curves in general and, particularly, cuspidal curves \cite{Heidi} \cite{Namba} \cite{Johnsen} \cite{Piene81}. In this thesis we will not use projections to construct cuspidal curves. Rather, we will use known properties of a particular cuspidal curve and the projection map to a posteriori analyze how this curve was constructed. 

In this chapter we will first give an outline of the method of projection in general. Then we will define the necessary tools to analyze a curve. Last, we will take a closer look at the construction of the cuspidal cubics and quartics.

\section{The projection map}
Let $(x_0:x_1:\ldots:x_n)$ denote the coordinates of a point in the n-dimensional projective space $\mathbb{P}^n$. Let $X$ be a projective variety of dimension $r-1$ in $\mathbb{P}^n$. Furthermore, let $V \subset \mathbb{P}^n$ be a linear subspace of dimension $n-r-1$. $V$ is called the {\em projection center}, and it can be given by the zero set 
\begin{equation*}
V=\V(H_0,\ldots,H_r),
\end{equation*}
where $H_i \in \mathbb{C}[x_0,\ldots,x_n]$, $i=0,\ldots,r$, are linearly independent linear polynomials, $$H_i=\sum_{k=0}^na_{ik}x_k.$$

\noi Let $A_V$ be the $(r+1) \times (n+1)$ coefficient matrix of the linear polynomials $H_i$,
$$A_V=\begin{bmatrix} a_{00}& \ldots &a_{0n}\\\noalign{\medskip} \vdots & \ddots & \vdots \\\noalign{\medskip} a_{r0}& \ldots &a_{rn}\end{bmatrix}.$$ \\

\noi With a variety $X$ and a projection center $V$ in $\mathbb{P}^n$ we define the {\em projection map} $\rho_V$\label{def:projection},
\begin{center}
$\begin{array}{cccc}
\rho_V:&X& \longrightarrow & \mathbb{P}^r \\
&\rotatebox[origin=c]{90}{$\in$}&&\rotatebox[origin=c]{90}{$\in$}\\
&(p_0: \ldots:p_n) & \longmapsto & (H_0: \ldots : H_r) \\
&&& \rotatebox[origin=c]{90}{$=$}\\
& & & (\sum a_{0k}p_k: \ldots:\sum a_{rk}p_k).
\end{array}$
\end{center}
In the language of matrices, this is nothing more than the matrix product
$$ \rho_V:\; \begin{bmatrix} p_0 \\\noalign{\medskip} \vdots \\\noalign{\medskip} p_n \end{bmatrix} \longmapsto \begin{bmatrix} a_{00}& \ldots &a_{0n}\\\noalign{\medskip} \vdots & \ddots & \vdots \\\noalign{\medskip}  a_{r0}& \ldots &a_{rn}\end{bmatrix} \cdot \begin{bmatrix} p_0 \\\noalign{\medskip} \vdots \\\noalign{\medskip} p_n \end{bmatrix}.$$

\noi Since $H_i$ are linearly independent, we have that the kernel $K_V$ of the map $\rho_V$, $K_V= \ker(A_V)$, is a linear subspace of $A_V$. $K_V$ can be given by $n-r$ linearly independent basis vectors, $$\vec{b}_i=\begin{bmatrix} b_{i0} & \ldots & b_{in} \end{bmatrix}, \qquad i=1, \ldots, n-r.$$ Furthermore, $K_V$ will frequently be given by a $(n-r) \times (n+1)$ matrix where the rows are given by the basis vectors, $$K_V=\begin{bmatrix} b_{10} & \ldots & b_{1n} \\ \nom \vdots & \ddots & \vdots \\ \nom b_{(n-r)0} & \ldots & b_{(n-r)n} \end{bmatrix}.\label{projectioncenter}$$ Moreover, the rows of the matrix $K_V$ span the projection center $V$, and we will therefore often describe $V$ by $K_V$. Note that we have the relations $K_V=\ker (A_V)$ and, conversely, $A_V= \ker(K_V)$.

\section{The rational normal curve}
All rational cuspidal curves in $\mathbb{P}^2$ are the resulting curves of different projections from a particular curve in $\mathbb{P}^n$. In this section we define the {\em rational normal curve} $C_n$ and some associated varieties of this curve.

\subsubsection{The rational normal curve}\label{def:rationalnormalcurve}
Let $\gamma$ be the map
\begin{center}
$\begin{array}{cccc}
\gamma: & \mathbb{P}^1 & \longrightarrow & \mathbb{P}^n\\
& \rotatebox[origin=c]{90}{$\in$}&&\rotatebox[origin=c]{90}{$\in$}\\
&(s:t) & \longmapsto & (s^n:s^{n-1}t:\ldots:st^{n-1}:t^n).\\
\end{array}$
\end{center}
The {\em rational normal curve} $C_n$ is a 1-dimensional variety in $\mathbb{P}^n$, given by $\mathrm{Im}(\gamma(s,t))$. It can be described in vector notation by $$\vec{\gamma}=\begin{bmatrix}s^n&s^{n-1}t&\ldots&st^{n-1}&t^n\end{bmatrix}.$$
\noi Additionally, the rational normal curve is given by the common zero set of $$x_ix_j-x_{i-1}x_{j+1},\qquad 1 \leq i \leq j \leq n-1.$$

\subsubsection{The tangent and the tangent developable}\label{def:tangent}
For every point $\gamma(s,t)$ of $C_n$ we define the {\em tangent} $T(s,t)$, 
$$T(s,t)=a_{00} \gamma(s,t)+a_{10} \frac{\partial}{\partial s} (\gamma(s,t))+a_{11}\frac{\partial}{\partial t}(\gamma(s,t)), \qquad a_{ij} \in \mathbb{C}.$$
By Euler's identity, the three terms above are linearly dependent. Hence, the tangent can be rewritten in matrix form as the row space of the matrix $T_M$,
\begin{equation*}
\begin{split}
T_M&=\begin{bmatrix}&\vec{\gamma}_s&\\\noalign{\medskip}& \vec{\gamma}_t&\end{bmatrix}\\ \nom 
&=\begin{bmatrix}ns^{n-1}&(n-1)s^{n-2}t&\ldots&t^{n-1}&0\\\noalign{\medskip} 0&s^{n-1}&\ldots&(n-1)st^{n-2}&nt^{n-1}\end{bmatrix}.
\end{split}
\end{equation*}\medskip


\label{def:tangentdevelopable}
\noi We define the {\em tangent developable} $T_n$ of $C_n$ as the union of all the tangents $T(s,t)$. It is a 2-dimensional surface in $\mathbb{P}^n$ which, by the homogeneity of the rational normal curve, has similar properties for all values $(s:t)$. We observe that $C_n \subset T_n$. The tangent developable $T_n$ is smooth outside $C_n$, but the rational normal curve constitutes a cuspidal edge on $T_n$. 

The tangent developable in $\mathbb{P}^n$ can be given by defining polynomials in $\mathbb{C}[x_0,x_1,\ldots,x_n]$ by elimination of $s$ and $t$, see Appendix \ref{calculationsandcode}. \medskip

\begin{ex} For degree $d=4$, the tangent developable $T_4$ is given by
{\footnotesize
\begin{equation*}
\begin{split}
[1]&=3x_2^2-4x_1x_3+x_0x_4,\\
[2]&=2x_1x_2x_3-3x_0x_3^2-3x_1^2x_4+4x_0x_2x_4,\\
[3]&=8x_1^2x_3^2-9x_0x_2x_3^2-9x_1^2x_2x_4+14x_0x_1x_3x_4-4x_0^2x_4^2.
\end{split}
\end{equation*}
}
\end{ex}

{\samepage{
\subsubsection{Osculating $k$-planes}\label{def:osculating}
For every point $\gamma(s,t)$ of the rational normal curve $C_n$ we define the {\em osculating $k$-plane} $O^k(s,t)$,
\begin{equation*}
\begin{split}
O^k(s,t)=&\;a_{00}\gamma+a_{10}\frac{\partial \gamma}{\partial s}+a_{11}\frac{\partial \gamma}{\partial t}+ \ldots \\ &+a_{k0}\frac{\partial^k \gamma}{\partial s^k}+a_{k1}\frac{\partial^k \gamma}{\partial s^{(k-1)} \partial t}+\ldots+a_{k(k-1)}\frac{\partial^k \gamma}{\partial s \partial t^{(k-1)}}+a_{kk}\frac{\partial^k \gamma}{\partial t^k},
\end{split}
\end{equation*} where $a_{ij} \in \mathbb{C},\; i=0,\ldots,k, \; j=0,\ldots,i$. Note that $T(s,t)=O^1(s,t)$.}}\\

\noi The terms of $O^k(s,t)$ are linearly dependent, hence the $k$-dimensional osculating $k$-plane can be rewritten in matrix form as the row space of the $(k+1)\times (n+1)$ matrix $O^k_M$,
\begin{equation*}
\begin{split}
O^k_M&=\begin{bmatrix}&\vec{\gamma}_{s^k}&\\\noalign{\medskip} &\vec{\gamma}_{s^{k-1}t}&\\\noalign{\medskip}  &\vdots& \\\noalign{\medskip} &\vec{\gamma}_{st^{k-1}}&\\\noalign{\medskip} &\vec{\gamma}_{t^{k}}& \end{bmatrix}.
\end{split}
\end{equation*}


\noi Observe that we have obvious relations between the rational normal curve, the tangents and the osculating $k$-planes. For every value of $(s:t)$, hence for every point $\gamma(s,t) \subset \mathbb{P}^n$, we have the chain $$\gamma(s,t) \subset T(s,t) \subset O^2(s,t) \subset O^3(s,t) \subset \ldots \subset O^{n-1}(s,t) \subset O^n(s,t)=\mathbb{P}^n.$$

\subsubsection{Secant variety}\label{def:secantvariety}
The {\em secant variety} $S_n$ of the rational normal curve can be given as the ideal generated by all $2\times 2$ minors of the matrix $\mathcal{S}_\alpha$ for any $\alpha$ such that $n-\alpha,\;\alpha \geq 2$ \cite[Prop. 9.7., p.103]{Harris},
$$\mathcal{S}_\alpha=\begin{bmatrix} x_0&x_1&x_2&\ldots&x_{n-\alpha}\\ \nom x_1&x_2&x_3&\ldots&x_{n-\alpha+1}\\ \nom
\vdots&&&&\\ \nom x_{\alpha}&x_{\alpha+1}&x_{\alpha+2}&\ldots&x_n\end{bmatrix}.$$ 

\noi The secant variety is a subspace of $\mathbb{P}^n$ with the property $$C_n \subset T_n \subset S_n.$$

\section{Cuspidal projections from $C_n$}
With a few exceptions there are so far not known {\em sufficient} conditions which can be imposed on the projection center $V$, such that the resulting curve $C'$ of a projection from $C_n$ is rational of degree $n$ and cuspidal, with cusps of a particular type. However, we do have some {\em necessary} conditions on the projection center so that the resulting curve is cuspidal.\\

\noi Let $\rho_V$ be the projection map mapping $C_n \subset \mathbb{P}^n$ to a curve $C' \subset \mathbb{P}^2$. Counting dimensions, the projection center $V$ of $\rho_V$ must be a $n-3$-dimensional linear subspace of $\mathbb{P}^n$. Furthermore, $V$ can be given by the intersection of the zero sets of three linearly independent linear polynomials, $$V=\V(H_0,H_1,H_2),$$ $$A_V=\begin{bmatrix} a_{00}& \ldots &a_{0n}\\\noalign{\medskip}  a_{10}& \ldots &a_{1n}\\\noalign{\medskip} a_{20}& \ldots &a_{2n}\end{bmatrix}.$$
The projection center $V$ can also be given by a $(n-2)\times(n+1)$ matrix $K_V$, where the rows consist of the $(n-2)$ linearly independent basis vectors $\vec{b}_{i}$ for the kernel of $A_V$.\\

\noi In the language of matrices, the projection map $\rho_V$ can be given by

$$\begin{array}{cccc}
\centering 
\rho_V:& \mathbb{P}^n & \longrightarrow &\Po\\ 
	& \cup & & \cup\\
	& C_n & \longrightarrow & C'\\
	\nom
	& \rotatebox[origin=c]{90}{$\in$}&&\rotatebox[origin=c]{90}{$\in$}\\
	\nom
&\begin{bmatrix} s^n\\\noalign{\medskip}s^{n-1}t\\\noalign{\medskip} \vdots \\\noalign{\medskip}st^{n-1}\\\noalign{\medskip}t^n \end{bmatrix} &\longmapsto& \begin{bmatrix} a_{00}& \ldots &a_{0n}\\\noalign{\medskip}  a_{10}& \ldots &a_{1n}\\\noalign{\medskip} a_{20}& \ldots &a_{2n}\end{bmatrix} \cdot \begin{bmatrix} s^n\\\noalign{\medskip}s^{n-1}t\\\noalign{\medskip} \vdots \\\noalign{\medskip}st^{n-1}\\\noalign{\medskip}t^n \end{bmatrix}. 
\end{array}$$\\

\noi If the resulting curve $C'$ of a projection map $\rho_V$ of $C_n$ is rational and cuspidal of degree $n$, we say that $\rho_V$ is a {\em cuspidal projection}. If $V$ fulfills the following criteria, then $\rho_V$ is a cuspidal projection \cite[pp.89--90]{Johnsen} \cite[pp.95--97]{Piene81}.

\begin{itemize}
\item[--] $V$ can not intersect $C_n$. If it did, $C'$ would not be of degree $n$. Using the matrices above, we find that this is equivalent to the criterion $$\mathrm{rank}\,C_nV=n-1,$$ where $C_nV$ is the $(n-1) \times (n+1)$ matrix
$$C_nV=\begin{bmatrix} s^n&s^{n-1}t& \ldots &st^{n-1}&t^n\\\noalign{\medskip} &&K_V&& \\\noalign{\medskip} \end{bmatrix}.$$

\item[--] $V$ must intersect the tangent developable $T_n$. If $V$ intersects $T_n$ at $T(s_0,t_0)$, then the image of the point $\gamma(s_0,t_0)$ will be a cusp on $C'$. This is easily seen by looking at the Puiseux parametrization of a cusp, which is on the form $$(t^m:c_rt^r+\ldots:1), \qquad r>m>1.$$ Getting $j$ cusps on a curve requires that $V$ intersects the tangent developable $T_n$ in $j$ points. This is equivalent to $$\mathrm{rank}\,T_nV=n-1$$ for $j$ values $(s:t)$, where $T_nV$ is the $n \times (n+1)$ matrix
$$T_nV=\begin{bmatrix} ns^{n-1}&(n-1)s^{n-2}t&\ldots&t^{n-1}&0\\\noalign{\medskip} 0&s^{n-1}&\ldots&(n-1)st^{n-2}&nt^{n-1}\\ \noalign{\medskip} &&K_V&& \\\noalign{\medskip} \end{bmatrix}.$$

\item[--] $V$ can not intersect the secant variety $S_n$ of $C_n$ outside $T_n$. If $V$ did intersect $S_n \setminus T_n$, $C'$ would not be purely cuspidal. How to impose this restriction is unknown. However, given $V$ and using the matrix representation of $S_n$ on page \pageref{def:secantvariety}, it can be checked that $V$ does not intersect $S_n \setminus T_n$.

\end{itemize}		 

\begin{rem}
Observe that since $T(s_0,t_0) \subset O^k(s_0,t_0)$ for all $k > 1$, any projection center $V$ which intersects $T(s_0,t_0)$ will intersect $O^k(s_0,t_0)$ as well.\\
\end{rem}

\noi If $V$ fulfills the above criteria, then we get a rational cuspidal curve with $j$ cusps from the projection map $\rho_V$ of $C_n$. However, we do not know how different choices of $V$ give different kinds of cusps on $C'$. Although it is possible to give qualified suggestions for $V$ in order to get cusps with relatively simple multiplicity sequences on $C'$, finding general patterns for more complex cases seems difficult. To illustrate this problem, we will briefly explore the subject for the quartic curves later in this chapter.\\

\noi Although not directly involved in the discussion of the number of or types of cusps on a curve, inflection points represent important properties of a curve. We have criteria on $V$ so that an inflection point is produced by the projection map. 

If $V$ intersects $O^2(s_0,t_0)$, but it intersects neither the curve $C_n$ nor the tangent $T(s_0,t_0)$, then the image of the point $\gamma(s_0,t_0)$ will be an inflection point on $C'$. This follows from the Puiseux parametrization of an inflection point, which generally is on the form $$(t:c_rt^r+\ldots:1), \qquad r\geq3.$$ Furthermore, getting $j$ inflection points on a curve requires that, for $j$ values $(s:t)$, $$\det O^2_nV=0 \qquad \text{and} \qquad V \cap T(s,t)=\emptyset,$$ where $O^2_nV$ is the $(n+1) \times (n+1)$ matrix
$$O^2_nV=\begin{bmatrix} &\vec{\gamma}_{s^2}& \\\nom &\vec{\gamma}_{st}& \\\nom &\vec{\gamma}_{t^2}& \\ \noalign{\medskip}&K_V& \\\noalign{\medskip} \end{bmatrix}.$$\medskip


\noi Imposing the above restrictions and finding the appropriate projection centers so that we may construct curves by projection is quite hard. For maximally inflected curves of degree $d=4$, this was done by Mork in \cite[pp.45--62]{Heidi}. But generally there are many restrictions and many unknown parameters which have to be determined. Therefore, even with the help of computer programs, it is difficult to find suitable projection centers. 

Since we, by for example \cite[pp.135,146,179--182]{Namba} and \cite[pp.327--328]{Fen99b}, explicitly know some cuspidal curves, we will use the parametrization of the curves and essentially read off the associated projection centers $V$ instead of determining it based on the restrictions. For some curves the projection centers will be investigated closely for the purpose of finding out more about the geometry of cuspidal projections. This will involve intersecting the projection center with different osculating $k$-planes.

\section{Cuspidal projections from $C_3$}
A projection of the rational normal curve $C_3$ in $\mathbb{P}^3$ to a curve $C'$ in \Po must have as projection center a 0-dimensional linear variety, a point $P$, $$K_P=\begin{bmatrix}b_0&b_1&b_2&b_3\end{bmatrix}.$$ Because $T_3$ is smooth outside $C_3$, there is no point on the surface $T_3 \setminus C_3$ where two or more tangents intersect. Hence, we can maximally have one cusp on the curve $C'$. 

The cuspidal cubic curve can be represented by the parametrization $$(s^3:st^2:t^3).$$
We read off the parametrization that $P$ is given by
$$P=\V(x_0,x_2,x_3),$$
$$A_P=\begin{bmatrix}1&0&0&0\\\noalign{\medskip}0&0&1&0\\\noalign{\medskip}0&0&0&1\end{bmatrix},\qquad K_P=\begin{bmatrix}0&1&0&0\end{bmatrix}.$$

\noi With this information on the projection center we observe, with the help of {\em Singular}, that
\begin{itemize}
\item[--] $P \cap C_3=\emptyset$. \\
\noi For all $(s:t)$ we have $$\mathrm{rank}\,\begin{bmatrix}s^3&s^2t&st^2&t^3\\ 0&1&0&0 \end{bmatrix}=2.$$
{\scriptsize
\begin{verbatim}
ring r=0, (s,t),dp;
matrix C[2][4]=s3,s2t,st2,t3,0,1,0,0;
ideal I=(minor(C,2));
solve(std(I));

[1]:
   [1]:
      0
   [2]:
      0
\end{verbatim}
}

\item[--] $P \cap T_3=\{p_1\}$.\\
\noi For $(s:t)=\{(1:0)\}$, we have $$\mathrm{rank}\,\begin{bmatrix}3s^2&2st&t^2&0\\\noalign{\medskip} 0&s^2&2st&3t^2 \\\noalign{\medskip}0&1&0&0 \end{bmatrix}=2.$$
{\scriptsize
\begin{verbatim}
matrix T[3][4]=3s2,2st,t2,0,0,s2,2st,3t2,0,1,0,0;
ideal I=(minor(T,3));
ideal Is=I,s-1;
ideal It=I,t-1;
solve(std(Is));

[1]:
   [1]:
      1
   [2]:
      0

solve(std(It));

? ideal not zero-dimensional
\end{verbatim}
}
This corresponds to precisely one cusp on the curve.

\item[--] $P \cap O^2(s,t)=\{p_1,p_2\}$.\\
\noi For $(s:t)=\{(1:0),(0:1)\}$, we have $$\det \begin{bmatrix}3s&t&0&0\\ 0&2s&2t&0 \\0&0&s&3t \\0&1&0&0 \end{bmatrix}=0.$$
{\scriptsize
\begin{verbatim}
matrix O_2[4][4]=3s,t,0,0,0,2s,2t,0,0,0,s,3t,0,1,0,0;
ideal I=(det(O_2));
ideal Is=I,s-1;
ideal It=I,t-1;
solve(std(Is));

[1]:
   [1]:
      1
   [2]:
      0

solve(std(It));

[1]:
   [1]:
      0
   [2]:
      1

\end{verbatim}
}
Only one of these values, $(s:t)=(0:1)$, additionally fulfill the restriction $P \cap T(s,t)=\emptyset$. Hence, we have one inflection point on the curve.

\end{itemize}

\noi We conclude that the cuspidal cubic has one cusp and one inflection point.





%

\section{Cuspidal projections from $C_4$}\label{cuspidalprojections4}
The projection center of a projection of the rational normal curve $C_4$ in $\mathbb{P}^3$ to a curve $C'$ in $\mathbb{P}^2$ must be a 1-dimensional linear variety, a line $L$, given by
$$K_L=\begin{bmatrix}b_{10}&b_{11}&b_{12}&b_{13}&b_{14}\\ \noalign{\medskip}b_{20}&b_{21}&b_{22}&b_{23}&b_{24}\end{bmatrix}.$$

In order to give a cuspidal projection, the following conditions must always be fulfilled by $L$.
\begin{itemize}
\item[--] $L \cap C_4=\emptyset$.
$$\mathrm{rank}\,\begin{bmatrix}s^4&s^3t&s^2t^2&st^3&t^4\\\noalign{\medskip}b_{10}&b_{11}&b_{12}&b_{13}&b_{14}\\\noalign{\medskip} b_{20}&b_{21}&b_{22}&b_{23}&b_{24}\end{bmatrix}=3 \; \text{ for all } (s:t).$$
\item[--] $L \cap T_4=\{p_j\}$.
$$ \mathrm{rank}\,\begin{bmatrix}4s^3&3s^2t&2st^2&t^3&0\\\noalign{\medskip}0&s^3&2s^2t&3st^2&4t^3\\\noalign{\medskip}b_{10}&b_{11}&b_{12}&b_{13}&b_{14}\\\noalign{\medskip} b_{20}&b_{21}&b_{22}&b_{23}&b_{24}\end{bmatrix}=3 \; \text{ for $j$ values of } (s:t).$$ 
\end{itemize}

\noi In \cite[pp.55,65]{Telling} it is explained that after a projection of $C_4$ to a curve $C' \subset \mathbb{P}^3$ from a point in $\mathbb{P}^4$, there exist points in $\mathbb{P}^3$ where tangents of $C'$ meet. The maximal number of tangents that intersect in a point in $\mathbb{P}^3$, regardless of the projection center of the initial projection, is three. This translates to the fact that any line $L \subset \mathbb{P}^4$ can intersect the tangent developable $T_4$ in maximally three points in $\mathbb{P}^4$. Hence, the maximal number of cusps of a rational quartic curve is three, which we have already seen is true.\\ 


\noi We have have $j$ inflection points on the curve if, for $j$ pairs $(s:t)$,$$L \cap T(s,t)=\emptyset \qquad \text{and} \qquad L \cap O^2(s,t)=\{p_j\}.$$ The last requirement is equivalent to $$\det \begin{bmatrix}6s^2&3st&t^2&0&0\\\nom 0&3s^2&4st&3t^2&0 \\\nom 0&0&s^2&3st&6t^2 \\\noalign{\medskip}b_{10}&b_{11}&b_{12}&b_{13}&b_{14}\\\noalign{\medskip} b_{20}&b_{21}&b_{22}&b_{23}&b_{24} \end{bmatrix}=0.$$

\noi We may investigate the position of $L$ in further detail. We know that $L \subset O^2(s,t)$ if  
$$\mathrm{rank}\,\begin{bmatrix}6s^2&3st&t^2&0&0\\\nom 0&3s^2&4st&3t^2&0 \\\nom 0&0&s^2&3st&6t^2 \\\noalign{\medskip}b_{10}&b_{11}&b_{12}&b_{13}&b_{14}\\ \nom b_{20}&b_{21}&b_{22}&b_{23}&b_{24} \end{bmatrix}=3.$$

\noi Additionally, we know that $L \subset O^3(s,t)$ if
$$\mathrm{rank}\,\begin{bmatrix}4s&t&0&0&0\\\nom 0&3s&2t&0&0 \\\nom 0&0&2s&3t&0 \\ \nom 0&0&0&s&4t \\\noalign{\medskip}b_{10}&b_{11}&b_{12}&b_{13}&b_{14} \\ \nom b_{20}&b_{21}&b_{22}&b_{23}&b_{24}\end{bmatrix}=4.$$

\noi The importance of this will become apparent when we next discuss the projection centers of the rational cuspidal quartic curves. See Appendix \ref{calculationsandcode} for the code used in Singular to produce these results.

\subsubsection{\bf Tricuspidal quartic -- $[(2),(2),(2)]$} 
The cuspidal quartic curve with three $A_2$ cusps is given by the parametrization $$(s^3t-\textstyle{\frac{1}{2}}s^4:s^2t^2:t^4-2st^3).$$ This parametrization corresponds to the projection center $L$, described by
$$L=\V(x_1-\textstyle{\frac{1}{2}}x_0,x_2,x_4-2x_3),$$
$$A_L=\left[ \begin {array}{ccccc} -\textstyle{\frac{1}{2}}&1&0&0&0\\\noalign{\medskip}0&0&1&0&0
\\\noalign{\medskip}0&0&0&-2&1\end {array} \right], \qquad K_L=\begin{bmatrix}2&1&0&0&0\\\noalign{\medskip}0&0&0&1&2\end{bmatrix}.$$

\noi With this information on the projection center we observe that
\begin{itemize}
\item[--] $L \cap C_4=\emptyset.$
\item[--] $L \cap T_4=\{p_1,p_2,p_3\}$ for $(s:t)=\{(1:0),(1:1),(0:1)\}$.
\item[--] $L \cap O^2(s,t)=\{p_1,p_2,p_3\}$ for $(s:t)=\{(1:0),(1:1),(0:1)\}$.
\item[--] $L \nsubseteq O^2(s,t)$ for any $(s:t)$.
\item[--] $L \nsubseteq O^3(s,t)$ for any $(s:t)$.
\end{itemize}

\noi The results imply well known properties of this curve. It has three cusps and no inflection points. We furthermore observe that we get an $A_2$-cusp when $L$ intersects $T_4$, but is not contained in any other osculating $k$-plane. This is consistent with the standard Puiseux parametrization of an $A_2$-cusp, i.e., $(t^2:t^3+\ldots:1)$, and the results of Mork in \cite[p.47]{Heidi}.

\subsubsection{\bf Bicuspidal quartic -- $[(2_2),(2)]$} 
The cuspidal quartic curve with two cusps, one $A_4$-cusp and one $A_2$-cusp, and one inflection point of type $1$ is given by the parametrization $$(s^4+s^3t:s^2t^2:t^4).$$ This parametrization corresponds to the projection center $L$, described by $$L=\V(x_0+x_1,x_2,x_4),$$
$$A_L=\left[ \begin {array}{ccccc} 1&1&0&0&0\\\noalign{\medskip}0&0&1&0&0
\\\noalign{\medskip}0&0&0&0&1\end {array} \right], \qquad K_L=\begin{bmatrix}0&0&0&1&0\\\noalign{\medskip}-1&1&0&0&0\end{bmatrix}.$$

\noi With this information on the projection center we observe that
\begin{itemize}
\item[--] $L \cap C_4=\emptyset.$
\item[--] $L \cap T_4=\{p_1,p_2\}$ for $(s:t)=\{(1:0),(0:1)\}$.
\item[--] $L \cap O^2(s,t)=\{p_1,p_2,p_3\}$ for $(s:t)=\{(1:0),(0:1),(1:-\frac{8}{3})\}$. 
\item[--] $L \nsubseteq O^2(s,t)$ for any $(s:t)$.
\item[--] $L \subset O^3(s,t)$ for $(s:t)=\{(1:0)\}$.
\end{itemize}

\noi The first three observations are consistent with the fact that we have two cusps and one inflection point on this quartic. Interestingly, the last observation reveals that the two cusps are different. We see that $L \subset O^3(s,t)$ for $(s:t)=(1:0)$. Although the below analysis of the cuspidal quartic with one ramphoid cusp of type $A_6$ reveals that this is not a sufficient condition, it apparently accounts for the $A_4$-type of the cusp corresponding to this value of $(s:t)$. This is consistent with the parametrization of the $A_4$-cusp given in Table \ref{tab:pui} on page \pageref{tab:pui} and the results of Mork in \cite[p.50]{Heidi}.

\subsubsection{\bf Unicuspidal ramphoid quartic -- $[(2_3)]$} 
The cuspidal quartic curve with one ramphoid cusp, an $A_6$-cusp, and three inflection points of type $1$ is given by the parametrization $$(s^4+st^3:s^2t^2:t^4).$$ This parametrization corresponds to the projection center $L$, described by
$$L=\V(x_0+x_3,x_2,x_4),$$
$$A_L=\left[ \begin {array}{ccccc} 1&0&0&1&0\\\noalign{\medskip}0&0&1&0&0
\\\noalign{\medskip}0&0&0&0&1\end {array} \right], \qquad K_L=\begin{bmatrix}0&1&0&0&0\\\noalign{\medskip}-1&0&0&1&0\end{bmatrix}.$$

\noi With this information on the projection center we observe that
\begin{itemize}
\item[--] $L \cap C_4=\emptyset.$
\item[--] $L \cap T_4=\{p_1\}$ for $(s:t)=\{(1:0)\}$.
\item[--] $L \cap O^2(s,t)=\{p_1,p_2,p_3,p_4\}$ for $$(s:t)=\{(1:0),(1:-2),(1:1-\mathrm{i}\sqrt{3}),(1:1+\mathrm{i}\sqrt{3})\}.$$ 
\item[--] $L \nsubseteq O^2(s,t)$.
\item[--] $L \subset O^3(s,t)$ for $(s:t)=\{(1:0)\}$.
\end{itemize}

\noi The first three observations confirm that we have one cusp and three inflection points on the quartic curve. Furthermore, we observe that $L \subset O^3(s,t)$ for $(s:t)=(1:0)$, which seems to be the reason why the cusp is of type $A_6$. This is consistent with the parametrization of the $A_6$-cusp given in Table \ref{tab:pui} on page \pageref{tab:pui}, where there is no term with $t^3$.\\

\begin{rem} Note that the above discussion reveals that we do not have sufficient criteria for the formation of $A_4$- and $A_6$-cusps on curves of degree $d=4$ under projection. We need further restrictions, and several attempts have been made to find these. David and Wall prove in \cite[Lemma 4.2., p.558]{David} and \cite[p.363]{Wall} that the production of an $A_6$-cusp is a special case of the production of an $A_4$-cusp. Translating these results to the language of matrices has so far not been successful.
\end{rem}

\subsubsection{\bf Ovoid quartic A-- $[(3)]$} 
The cuspidal quartic curve with one ovoid cusp, an $E_6$-cusp, and one inflection point of type $2$ is given by the parametrization $$(s^4:st^3:t^4).$$ This parametrization corresponds to the projection center $L$, described by
$$L=\V(x_0,x_3,x_4),$$
$$A_L=\left[ \begin {array}{ccccc} 1&0&0&0&0\\\noalign{\medskip}0&0&0&1&0
\\\noalign{\medskip}0&0&0&0&1\end {array} \right], \qquad K_L=\begin{bmatrix}0&1&0&0&0\\\noalign{\medskip}0&0&1&0&0\end{bmatrix}.$$

\noi With this information on the projection center we observe that
\begin{itemize}
\item[--] $L \cap C_4=\emptyset.$
\item[--] $L \cap T_4=\{p_1\}$ for $(s:t)=\{(1:0)\}$.
\item[--] $L \cap O^2(s,t)=\{p_1,p_2\}$ for $(s:t)=\{(1:0),(0:1)\}$.
\item[--] $L \subset O^2(s,t)$ for $(s:t)=\{(1:0)\}$.
\item[--] $L \subset O^3(s,t)$ for $(s:t)=\{(1:0),(0:1)\}$.
\end{itemize}

\noi We observe that $L$ intersects $T_4$ for the value $(s:t)=(1:0)$, and that the image of this point will be a cusp. We also note that $L$ additionally is contained in $O^2(s,t)$ for $(s:t)=(1:0)$, which accounts for the multiplicity $m=3$ of the cusp. This is consistent with the Puiseux parametrization of the cusp, $(t^3:t^4:1)$, and Mork in \cite[p.50]{Heidi}. 

Additionally, $L$ intersects $O^2(s,t)$ for $(s:t)=(0:1)$, hence the image of this point is an inflection point. We also note that $L$ is contained in $O^3(s,t)$ for $(s:t)=(0:1)$, which ensures that the inflection point is of type 2. 


\subsubsection{\bf Ovoid quartic B -- $[(3)]$} 
The cuspidal quartic curve with one ovoid cusp, an $E_6$-cusp, and two inflection points of type $1$ is given by the parametrization $$(s^3t-s^4:st^3:t^4).$$ This parametrization corresponds to the projection center $L$, described by
$$L=\V(x_1-x_0,x_3,x_4),$$
$$A_L=\left[ \begin {array}{ccccc} -1&1&0&0&0\\\noalign{\medskip}0&0&0&1&0
\\\noalign{\medskip}0&0&0&0&1\end {array} \right], \qquad K_L=\begin{bmatrix}1&1&0&0&0\\\noalign{\medskip}0&0&1&0&0\end{bmatrix}.$$

\noi With this information on the projection center we observe that
\begin{itemize}
\item[--] $L \cap C_4=\emptyset.$
\item[--] $L \cap T_4=\{p_1\}$ for $(s:t)=\{(1:0)\}$.
\item[--] $L \cap O^2(s,t)=\{p_1,p_2,p_3\}$ for $(s:t)=\{(1:0),(1:2),(0:1)\}$.
\item[--] $L \subset O^2(s,t)$ for $(s:t)=\{(1:0)\}$.
\item[--] $L \subset O^3(s,t)$ for $(s:t)=\{(1:0)\}$.
\end{itemize}

\noi As above, the first two and the fourth observation accounts for the cusp and its multiplicity $m=3$. The third observation is consistent with the fact that we have two inflection points on this curve, and these are of type 1 by the last observation.

\chapter{Cremona transformations}\label{ccc}
The concept of {\em Cremona transformations} provides a powerful tool that makes construction of curves in general, and cuspidal curves in particular, quite simple. Constructing cuspidal curves with Cremona transformations can be approached in two radically different, although entwined, ways. We may regard Cremona transformations algebraically and apply an explicitly given transformation to a polynomial, or we may use geometrical properties of Cremona transformations and implicitly prove the existence of curves. In this chapter we will first define and describe Cremona transformations, with particular focus on quadratic Cremona transformations. Then we will use both given approaches to construct and prove the existence of the cuspidal cubics and quartics.

\section{Quadratic Cremona transformations}\label{def:cremona}
Let $\psi$ be a {\em birational transformation},
$$\begin{array}{rccc}
\psi: & \mathbb{P}^2 &\dashrightarrow &\mathbb{P}^2\\
&\rotatebox[origin=c]{90}{$\in$}&&\rotatebox[origin=c]{90}{$\in$}\\
& (x:y:z) &\longmapsto &(G^x(x,y,z):G^y(x,y,z):G^z(x,y,z)),
\end{array}$$
where $G^x$, $G^y$ and $G^z$ are linearly independent homogeneous polynomials of the same degree $d$ with no common factor. Then $\psi$ is called a plane {\em Cremona transformation} of {\em order $d$}. In particular, if $\psi$ is a birational transformation of order $2$, then $\psi$ is called a {\em quadratic Cremona transformation}.\\

\noi The birational transformations are precisely the maps for which the set $\V(G^x) \cap \V(G^y) \cap \V(G^z)$ consists of exactly $d^2-1$ points, counted with multiplicity \cite{Moe}. These points are called the {\em base points} of $\psi$. This immediately implies that any linear change of coordinates is a Cremona transformation. A quadratic Cremona transformation must have $3$ base points, counted with multiplicity.

By Hartshorne \cite[Theorem 5.5, p.412]{Hart:1977} it is possible to factor a birational transformation into a finite sequence of monoidal transformations and their inverses. More specifically, a quadratic Cremona transformation acts on $\mathbb{P}^2$ by blowing up the three base points and blowing down three associated lines. This process will transform any curve $C=\mathcal{V}(F) \subset \mathbb{P}^2$, and we will go through the process in detail in the next sections. 

Before we discuss how Cremona transformations can transform curves, we observe that by allowing the base points to be not only proper points in $\mathbb{P}^2$, but also infinitely near points of any other base point, we get three apparently different kinds of quadratic Cremona transformations \cite[Section 2.8., pp.63--66]{Albe}.\medskip \\
\noi {\bf Three proper base points -- $\psi_3$}\\
\noi $\psi_3$ is a Cremona transformation which has three proper base points, $p$, $q$ and $r$, in $\mathbb{P}^2.$ We will write $\psi_3(p,q,r)$ for this transformation. Note that since we have three base points, orienting curves in $\mathbb{P}^2$ such that explicit applications of this kind of Cremona transformations have the desired effect, is easily done.\medskip \\
\noi {\bf Two proper base points -- $\psi_2$}\\
\noi $\psi_2$ is a Cremona transformation which has two proper base points, $p$ and $q$, in $\mathbb{P}^2$, and one infinitely near base point, $\hat{q}$. The latter point is here defined to be the infinitely near point of $q$ lying in the intersection of the exceptional line $E$ of $q$ and the strict transform of a specified line $L \subset \Po$ through $q$. We write $\psi_2(p,q,L)$ for this transformation. Note that although the point $\hat{q}$ is not in $\mathbb{P}^2$, we are still able to orient a curve such that this kind of Cremona transformation gives the sought after effect. This is because $q$ together with any point $r \in L \setminus{q} \subset \Po$ determine $L$.\medskip \\ 
\noi {\bf One proper base point -- $\psi_1$}\\
\noi $\psi_1$ is a Cremona transformation which has one proper base point, $p$, in $\mathbb{P}^2$. The last two base points, called $\hat{p}$ and $\hat{\hat{p}}$, are infinitely near points of $p$. The point $\hat{p}$ is an infinitely near point of $p$ in the intersection of the exceptional line $E$ of $p$ and the strict transform of a line $L$ through $p$. The point $\hat{\hat{p}}$ is an infinitely near point of both $p$ and $\hat{p}$, lying somewhere on the exceptional line $E_{\hat{p}}$ of the blowing-up of $\hat{p}$. We write $\psi_1(p,L,-)$ for this transformation. Note that there is no apparent representative for $\hat{\hat{p}}$ in $\mathbb{P}^2$, hence it may be difficult to appropriately orient curves to get the desired effect from  explicit applications of Cremona transformations of this kind.\\

\noi These quadratic Cremona transformations are not that different after all. All of them can be written as a product of linear transformations and quadratic Cremona transformations with three proper base points. The transformation with two proper base points is a product of two, and the transformation with one proper base point is a product of four quadratic Cremona transformations with three proper base points \cite[pp.246--247]{Albe}.

\section{Explicit Cremona transformations}
Cremona transformations can be applied to $\mathbb{P}^2$, thereby transforming a curve $C=\mathcal{V}(F)$ \cite[pp.\,170--178]{Fulton}. For each kind of quadratic Cremona transformation, we choose an explicit representation of $\psi$ such that $\psi \circ \psi$ is the identity on $\mathbb{P}^2$ outside the base points. We may apply $\psi$ to the curve $C$ by making the substitutions $x=G^x$, $y=G^y$ and $z=G^z$ in the defining polynomial $F$. The polynomial $F^Q(x,y,z)=F(G^x,G^y,G^z)$ may be reducible, and $C^Q=\V(F^Q)$ is called the {\em total transform} of $C$. Carefully removing linear factors of $F^Q$, we get a polynomial $F'(x,y,z)$, which is the defining polynomial of the {\em strict transform} $C'$ of $C$. Note that the removal of linear factors in $F^Q$ is depending on the particular Cremona transformation. Moreover, if $C$ is irreducible, so is $C'$. We also have that $(C')'=C$. \\

\noi The three different types of quadratic Cremona transformations presented in Section \ref{def:cremona} can be given explicitly by standard maps. By combining these maps with linear changes of coordinates, we can produce all quadratic Cremona transformations. \medskip \\ 
\noi{\bf Three proper base points -- $\psi_3$}\\
\noi After a linear change of coordinates which moves $p$ to $(1:0:0)$, $q$ to $(0:1:0)$, and $r$ to $(0:0:1)$, $\psi_3(p,q,r)$ can be written on standard form. This map is commonly referred to as the standard Cremona transformation.
\begin{equation*}
\psi_3: (x:y:z) \longmapsto (yz:xz:xy).
\end{equation*}
$\psi_3$ applied to $F$ gives $F^Q$, the defining polynomial of the total transform, $C^Q$. To get $F'$, the defining polynomial of the strict transform $C'$, remove the factors $x$, $y$ and $z$ with multiplicity $m_p$, $m_q$ and $m_r$ respectively \cite[p.173]{Fulton}. \medskip \\
\noi {\bf Two proper base points -- $\psi_2$}\\
\noi After a linear change of coordinates which moves $p$ to $(1:0:0)$, $q$ to $(0:1:0)$, and lets $L$ be the line $\mathcal{V}(x)$, $\psi_2(p,q,L)$ can be written on standard form. 
\begin{equation*}
\psi_2: (x:y:z) \longmapsto (z^2:xy:zx).
\end{equation*}
$\psi_2$ applied to $F$ gives $F^Q$. To get $F'$, remove the factors $x$ with multiplicity $m_p$ and $z$ with multiplicity $m_q+m_{\hat{q}}$. Note that $m_{\hat{q}}=0$ if and only if $L$ is {\em not} the tangent line to $C$ at $q$. If $L=T_q$, then $m_{\hat{q}}=m_{q.1}$. \medskip \\
\noi {\bf One proper base point -- $\psi_1$}\\
\noi Although we are not able to control the orientation of this kind of Cremona transformation completely, the following explicit map is appropriate in later examples. The base point of this transformation is $(1:0:0)$.
\begin{equation*}
\psi_1: (x:y:z) \longmapsto (y^2-xz:yz:z^2).
\end{equation*}
$\psi_1$ applied to $F$ gives the defining polynomial $F^Q$ of the total transform $C^Q$. To get $F'$, remove the factor $z$ with multiplicity $m_p+m_{\hat{p}}+m_{\hat{\hat{p}}}$. These multiplicities are hard to find explicitly. It is easier to use the implicit approach to get the complete picture.

\section{Implicit Cremona transformations}
Explicitly applying Cremona transformations is useful for obtaining defining polynomials of strict transforms of curves. It hides, however, much of the ongoing action in the process of transforming a curve. By describing the Cremona transformation implicitly, much information about the strict transform of a curve, and in particular its singularities, can be deduced directly, i.e., without using the defining polynomial and without needing to worry about explicit orientation in $\mathbb{P}^2$. In this section we will first describe how we can find the degree of the strict transform of a curve. Then we will investigate each kind of Cremona transformations and describe them implicitly in order to clarify notation. Examples will be given in the next section.

\subsection{The degree of the strict transform}
Observe that we may estimate the degree $d'$ of the strict transform $C'$ of the curve $C$ of degree $d$ without explicitly calculating the strict transform. $F^Q$ has degree $2\cdot d$, and to get $F'$ we remove linear factors with a known multiplicity. Abusing notation, we let $p$, $q$ and $r$ denote the three base points of $\psi_j$, regardless of the nature of the point. Let $m_p$ denote the multiplicity of the point $p$ with respect to the curve $C_i$, where $i=0$ if $p$ is a proper base point, $i=1$ if $p$ is in the exceptional line of any proper base point, and $i=2$ if $p$ is in the exceptional line of any non-proper base point. The degree $d'$ of $C'$ is given by
$$d'=2 \cdot d - m_p-m_q-m_r.$$

\subsection{Three proper base points}
A quadratic Cremona transformation of type $\psi_3(p,q,r)$ can be regarded as nothing more than a simultaneous blowing-up of three points and a blowing-down of the strict transforms of the three lines connecting them. By the properties of monoidal transformations, it is possible to deduce intersection multiplicities and multiplicity sequences of points on the strict transform $C'$ of a curve $C$ under this type of Cremona transformation.

{\samepage{
\subsection{Elementary transformations}
Before we discuss the quadratic Cremona transformations with two, respectively one, proper base point, we introduce the concept of {\em elementary transformations} on {\em ruled surfaces}.}}

In this thesis we adopt the notation of Fenske in \cite{Fen99a} and define a {\em ruled surface} as a surface $X$ isomorphic to a Hirzebruch surface $$\Sigma_n=\mathbb{P}(\mathcal{O}_{\mathbb{P}^1}\oplus \mathcal{O}_{\mathbb{P}^1}(n)), \qquad \text{ for some } n>0.$$ $X$ has a horizontal section $E$ and vertical fibers $F$. The fibers $F$ have self-intersection $F^2=0$, and the horizontal section has self-intersection $E^2=-n$. Observe that for $n=1$ this is equivalent to the surface obtained via a blowing-up of a point $p$ in \Po. The horizontal section $E$, with $E^2=-1$, is the exceptional line from the blowing-up, and the fibers $F$ are the strict transforms of the lines through $p$.

An {\em elementary transformation} of a fiber $F$ at a point $q \in F$ on a ruled surface $X$ is a composition of the blowing-up of $X$ at the point $q$ and the blowing-down of the strict transform of the fiber $F$. To shorten notation, we say that we have an elementary transformation in $q$.

The blowing-up of $X$ at $q$ produces a surface which we denote $\bar{X}$. On this surface we have the strict transforms $\bar{E}$ and $\bar{F}$ of $E$ and $F$ respectively. Additionally we have the new exceptional line of the blowing-up of $q$, denoted $\bar{E}_q$. 

Blowing down $\bar{F}$, we obtain a new ruled surface $X'$, where the fiber $F$ is replaced by the strict transform of the exceptional line $\bar{E}_q$ from the blowing-up of $q$. Hopefully making notation simpler, we call this fiber $F'$. Furthermore, the strict transform of the horizontal section $\bar{E}$ is referred to as $E'$. Note that $E'$ has the property that if $q \in E$, then $E'^2=E^2-1$. If $q \notin E$, then $E'^2=E^2+1$.

By the properties of the blowing-up process on page \pageref{properelm}, we are able to calculate intersection multiplicities between the strict transform of a curve and the fibers of the ruled surface throughout an elementary transformation. We are additionally able to predict changes in the multiplicity sequence of a cusp caused by the elementary transformation.

\subsection{Two proper base points}
Let $p$ and $q$ be points in $\mathbb{P}^2$. Let $L^{pq}$ be the line between the two points, and let $L$ be another line through $q$. Then the quadratic Cremona transformation of type $\psi_2(p,q,L)$ can be decomposed into three main steps.\\ 

\noi{\bf Blowing up at $q$}\\
\noi The first step is blowing up the point $q \in \mathbb{P}^2$, hence producing a ruled surface $X_1$ with an exceptional line $E_1$. We have $E_1^2=-1$.\\

\noi{\bf Elementary transformations in $p_1$ and $\hat{q}$}\\
\noi The second step is performing two elementary transformations of two fibers at two points on $X_1$. The two fibers are the strict transform of the line $L^{pq}$, $L^{pq}_1$, and the strict transform of the line $L$, $L_1$. The two points are the strict transform of $p$, $p_1 \in L^{pq}_1$, and the point $\hat{q}$, which is uniquely determined by the fact that $\hat{q} \in E_1 \cap L_1$. 

After performing the elementary transformations, we obtain another ruled surface, $X_2$. The horizontal section $E_2$ of this surface is the strict transform of $E_1$ under the elementary transformations. On $X_2$ we have fibers $L_2$ and $L_2^{pq}$, the strict transforms of the exceptional lines $E_{\hat{q}}$ and $E_{p_1}$ under the elementary transformation.\\

\noi {\bf Blowing down $E_2$}\\
\noi By properties of elementary transformations, $E_2$ satisfies $E_2^2=-1$. Hence, the third step of the Cremona transformation is blowing down $E_2$, which leads us back to \Po. \\

\noi Since the Cremona transformation is a composition of blowing-ups and \linebreak blowing-downs, we are able to follow the intersection multiplicities and the multiplicity sequences of points of a curve and transformations of the curve in every step of the Cremona transformation. In particular, we can determine the invariants for the strict transform of the curve.

\subsection{One proper base point}
Let $p$ be a point and let $L$ be a line through $p$ in $\mathbb{P}^2$. A quadratic Cremona transformation of type $\psi_1(p,L,-)$ can then be decomposed into four main steps.\\


\noi {\bf Blowing up at $p$}\\
\noi The first step is blowing up the point $p \in \mathbb{P}^2$. This results in a ruled surface $X_1$ with horizontal section $E_1$, where $E_1^2=-1$. On $X_1$ we additionally have the fiber $L_1$, the strict transform of the line $L$, and the point $\hat{p}$, which is given by $\hat{p}=E_1 \cap L_1$.\\

\noi {\bf Elementary transformation in $\hat{p}$}\\
\noi The second step is performing an elementary transformation of the fiber $L_1$ at the point $\hat{p}$ on $X_1$. This results in a new ruled surface $X_2$, with horizontal section $E_2$ and the fiber $L_2$. $E_2$ denotes the transform of the horizontal section $E_1$, and $L_2$ denotes the transform of the exceptional line $E_{\hat{p}}$ of the elementary transformation. Note that $E_2^2=-2.$\\

{\samepage
\noi {\bf Elementary transformation in $\hat{\hat{p}}$}\\
\noi The third step is performing an elementary transformation of the fiber $L_2$ at a point $\hat{\hat{p}}$, where $\hat{\hat{p}} \notin E_2$. This results in a new ruled surface $X_3$, in which we have the transform of $E_2$, $E_3$, where $E_3^2=-1$. Analogous to the above, we denote by $L_3$ the transform of the exceptional line $E_{\hat{\hat{p}}}$.}\\

\noi {\bf Blowing down $E_3$}\\
\noi Since $E_3^2=-1$, blowing down the horizontal section $E_3 \subset X_3$ is the last step in the Cremona transformation, taking us back to \Po.\\

\noi As above, we are able to follow the intersection multiplicities and the multiplicity sequences of points of a curve and transformations of the curve in every step of the Cremona transformation.

\section{Constructing curves}
In this section we will make use of the quadratic Cremona transformations described above and construct the cuspidal cubic and quartic curves. 

The pictures of the implicit Cremona transformations of curves displayed in this thesis are merely illustrations of the action described in the text, and they do not give a complete overview of the curves and surfaces. 

Some of the points and lines are marked with notation in the illustrations. Points which are going to be blown up are marked with name and a black dot. Lines which are going to be blown down are marked with two black dots. The curves $C$ stand out in the images, hence, they are not marked. Furthermore, vertical fibers in the ruled surfaces will be marked by name only when they are exceptional lines of blowing-ups. The other vertical fibers are possible to identify from the context. Moreover, in order to avoid problems with the notation, we will write $T^p$ for the tangent line to $C$ at $p$.

\subsubsection{Cuspidal cubic -- $[(2)]$}
We construct the cuspidal cubic with the help of a quadratic Cremona transformation with two proper base points.

\begin{description}
\item
\parpic[l]{\includegraphics[width=0.3\textwidth]{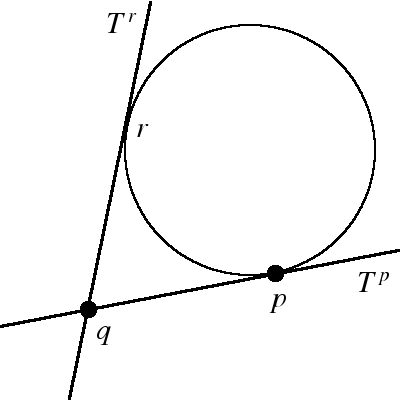}}\;
Let $C$ be an irreducible conic. Two arbitrary points $p$ and $r$ on $C$ have unique tangents $T^p$ and $T^r$, which only intersect $C$ at $p$ and $r$, respectively.
\begin{equation*}
\begin{split}
T^p \cdot C &= 2 \cdot p,\\
T^r \cdot C &= 2 \cdot r.
\end{split}
\end{equation*}
Let $q$ denote the intersection point $T^p \cap T^r$, and note that $q \notin C$.
\end{description}
\noi Applying the transformation $\psi_2(p,q,T^r)$ to $C$, we get the desired cubic.\\

\pagebreak
\noi {\bf Blowing up at $q$}
\begin{description}
\item
\parpic[l]{\includegraphics[width=0.3\textwidth]{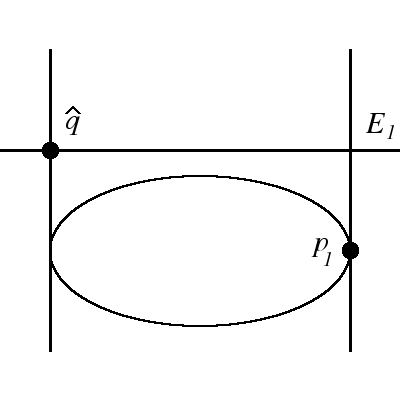}}\;
\noi Blowing up at $q$, we get the ruled surface $X_1$ with horizontal section $E_1$ and the transformed curve $C_1$.  No points or intersection multiplicities have been affected by this process. Note that since $q \notin C$, $E_1 \cap C_1= \emptyset$.

The points $\hat{q}=E_1 \cap T^r_1$ and $p_1=C_1 \cap T^p_1$ are marked in the figure. 
\end{description}

\noi{\bf Elementary transformations in $\hat{q}$ and $p_1$}
\begin{description}
\item 
\parpic[l]{\includegraphics[width=0.3\textwidth]{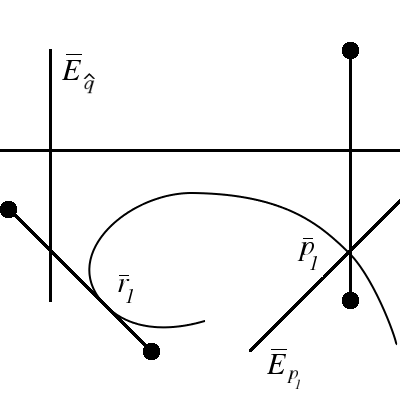}}\;
Blowing up at $\hat{q}$ and $p_1$ gives on $\bar{X}_1$ two exceptional lines $\bar{E}_{\hat{q}}$ and $\bar{E}_{p_1}$. We have the intersections 
\begin{equation*}
\begin{split}
\bar{T}^r_1 \cdot \bar{C}_1 &= 2 \cdot \bar{r}_1,\\
\bar{T}^r_1 \cap \bar{E}_1 &= \emptyset,\\
\bar{T}^p_1 \cdot \bar{C}_1 &= ((T^p_1 \cdot C_1)_{p_1} - m_{p_1}) \cdot \bar{p}_1\\
&=1 \cdot \bar{p}_1,\\
\bar{E}_{p_1} \cdot \bar{C}_1 &=1 \cdot \bar{p}_1.
\end{split}
\end{equation*}

\item 
\parpic[l]{\includegraphics[width=0.3\textwidth]{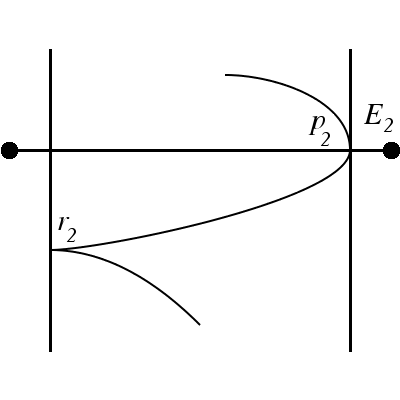}}\;
Blowing down $\bar{T}^r_1$ and $\bar{T}^p_1$ gives the surface $X_2$. On this surface we have $T^r_2$, the strict transform of $\bar{E}_{\hat{q}}$, and $T^p_2$, the strict transform of $\bar{E}_{p_1}$. 
\par Because of the above intersection concerning $\bar{r}_1$, $r_2$ is a cusp with multiplicity sequence $(2)$. Note that $r_2 \notin E_2$. 
\par The smooth point $\bar{p}_1$ is transformed into the smooth point $p_2$. Furthermore, $T^p_2$ is actually the tangent to $C_2$ at $p_2$, and $p_2 \in E_2$.
\begin{equation*}
\begin{split}
(T^p_2 \cdot C_2)_{p_2} &= ((\bar{T}^p_1 \cdot \bar{C}_1)_{\bar{p}_1}+(\bar{E}_{p_1} \cdot \bar{C}_1)_{\bar{p}_1})\\
&= 2,\\
(E_2 \cdot C_2)_{p_2} &= (\bar{T}^p_1 \cdot \bar{C}_1)_{\bar{p}_1}\\&= 1.
\end{split}
\end{equation*}

\end{description}

\pagebreak
\noi{\bf Blowing down $E_2$}
\begin{description}
\item
\parpic[l]{\includegraphics[width=0.3\textwidth]{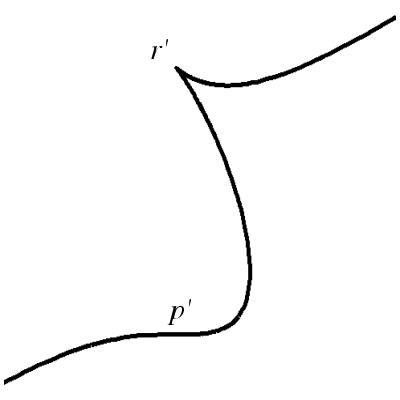}}\;
\noi Blowing down $E_2$ gives a curve $C'$ with one cusp $r'$ with multiplicity sequence $(2)$. Additionally, the above intersection multiplicities ensure that the point $p'$ is an inflection point.
\begin{equation*}
\begin{split}
T^{p'} \cdot C' &=((T^p_2 \cdot C_2)_{p_2} +(E_2 \cdot C_2)_{p_2}) \cdot p'\\
&= 3 \cdot p'.
\end{split}
\end{equation*}\medskip
\end{description}

\noi To see that $C'$ really is a cubic, observe that 
\begin{equation*}
\begin{split}
d'&=2 \cdot d - m_q-m_{\hat{q}}-m_{p_1}\\
&=2 \cdot 2 - 1\\
&=3.
\end{split}
\end{equation*}

\noi We may do this explicitly by using the Cremona transformation $\psi_2$ to transform the conic C, $$C=\mathcal{V}(y^2-xz).$$ We get the desired cubic curve $C'$ with an $A_2$-cusp, $$C'=\mathcal{V}(xy^2-z^3).$$

\subsubsection{$C_1$ -- Tricuspidal quartic -- $[(2),(2),(2)]$} 
Performing a quadratic Cremona transformation with three proper base points on an appropriately oriented irreducible conic produces the tricuspidal quartic curve. 

\begin{description}
\item
\parpic[l]{\includegraphics[width=0.3\textwidth]{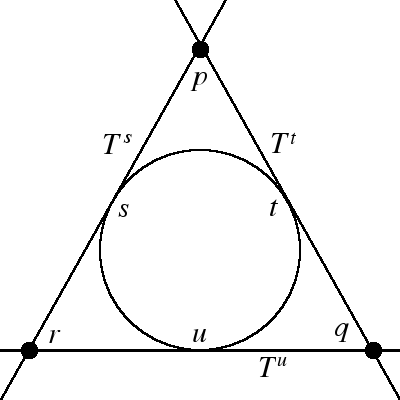}}\;
\noi Let $C$ be an irreducible conic. Choose three points $s$, $t$ and $u$ on $C$, and let $T^s$, $T^t$ and $T^u$ be the respective tangent lines to $C$ at the three points,
\begin{equation*}
\begin{split}
T^s \cdot C&=2 \cdot s,\\
T^t \cdot C&=2 \cdot t,\\
T^u \cdot C&=2 \cdot u.
\end{split}
\end{equation*}
By B\'{e}zout's theorem, the tangent lines intersect in three different points, $p$, $q$ and $r$, which can not be on $C$.
\end{description}

\pagebreak
\begin{description}
\item
\parpic[l]{\includegraphics[width=0.3\textwidth, height=0.3\textwidth]{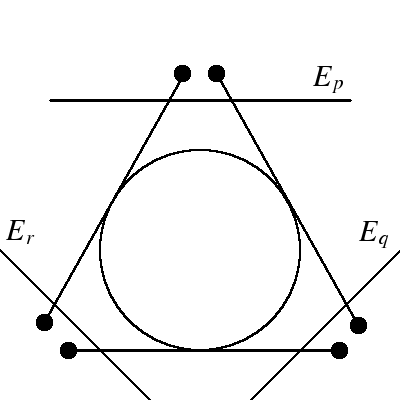}}\;
\noi Applying the Cremona transformation \linebreak $\psi_3(p,q,r)$, we get three exceptional lines $E_p$, $E_q$ and $E_r$ replacing $p$, $q$ and $r$.
\\
\\
\\
\\
\end{description}

\begin{description}
\item
\parpic[l]{\includegraphics[width=0.3\textwidth]{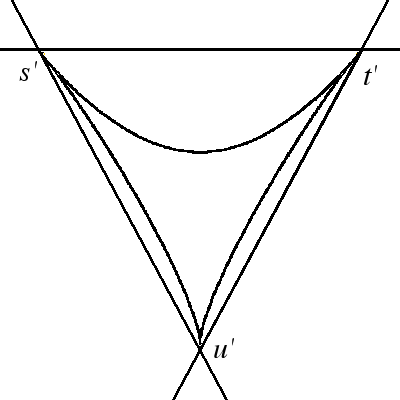}}\;
\noi Blowing down the tangent lines results in \linebreak points $s'$, $t'$ and $u'$ on $C'$, with multiplicity $$m_{s'}=m_{t'}=m_{u'}=2.$$ Since the points $s$, $t$ and $u$ were smooth points on $C$, $s'$, $t'$ and $u'$ are cusps with multiplicity sequence $(2)$.\medskip 
\end{description}

\vspace{4mm}
\noi The degree of $C'$ is $d'=2 \cdot d=4$ since the base points were not in $C$. Hence, we have constructed the tricuspidal quartic.\\

\noi Explicitly, we choose the appropriately oriented conic $C$, $$C=\mathcal{V}(x^2+y^2+z^2-2xy-2xz-2yz).$$ 
\noi Applying the Cremona transformation $\psi_3$ results in the desired tricuspidal quartic curve $C'$, $$C'=\mathcal{V}(y^2z^2+x^2z^2+x^2y^2-2xyz(x+y+z)).$$


\subsubsection{$C_2$ -- Bicuspidal quartic -- $[(2_2),(2)]$}\label{bcq}
We construct the bicuspidal quartic using a quadratic Cremona transformation with two proper base points.

\begin{description}
\item
\parpic[l]{\includegraphics[width=0.3\textwidth]{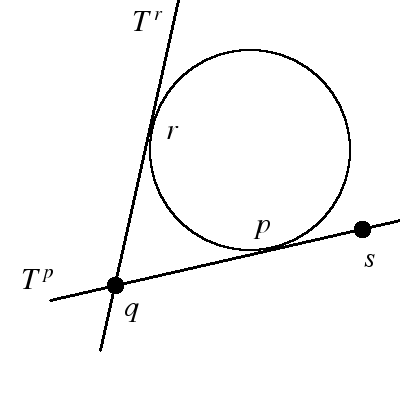}}\;
Let $C$ be an irreducible conic. Two arbitrary points $p$ and $r$ on $C$ have unique tangents $T^p$ and $T^r$, which only intersect $C$ at $p$ and $r$, respectively.
\begin{equation*}
\begin{split}
T^p \cdot C &= 2 \cdot p,\\
T^r \cdot C &= 2 \cdot r.
\end{split}
\end{equation*} 
Let $q$ denote the intersection point $T^p \cap T^r$, and let $s$ denote another point on $T^p$. Note that $q,s \notin C$.
\end{description}

\pagebreak
\noi Applying the transformation $\psi_2(s,q,T^r)$ to $C$, we get the desired quartic $C'$.\\

\noi {\bf Blowing up at $q$}
\begin{description}
\item
\parpic[l]{\includegraphics[width=0.3\textwidth]{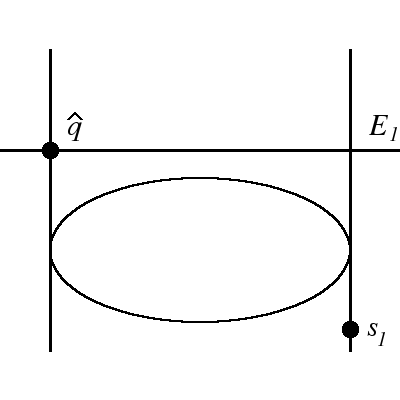}}\;
\noi Blowing up at $q$, we get the ruled surface $X_1$ with horizontal section $E_1$, fibers $T^p_1$ and $T^r_1$, and the transformed curve $C_1$.  No points or intersection multiplicities have been affected by this process.  Since $q \notin C$, we have $E_1 \cap C_1= \emptyset$.

The points $\hat{q}=E_1 \cap T^r_1$ and $s_1 \in T^p_1$ are marked in the figure.\\ 
\end{description}

\noi{\bf Elementary transformations in $\hat{q}$ and $s_1$}
\begin{description}
\item 
\parpic[l]{\includegraphics[width=0.3\textwidth]{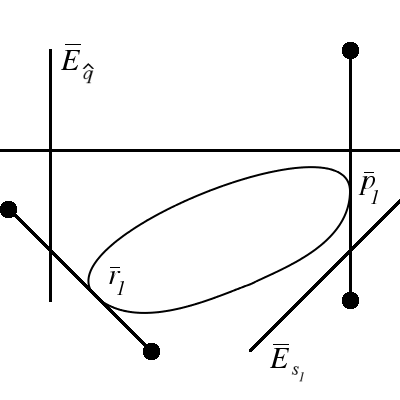}}\;
Blowing up at $\hat{q}$ and $s_1$ gives on $\bar{X}_1$ two exceptional lines $\bar{E}_{\hat{q}}$ and $\bar{E}_{s_1}$. We have the intersections 
\begin{equation*}
\begin{split}
\bar{T}^r_1 \cdot \bar{C}_1 &= 2 \cdot \bar{r}_1,\\
\bar{T}^r_1 \cap \bar{E}_1 &= \emptyset,\\
\bar{T}^p_1 \cdot \bar{C}_1 &=2 \cdot \bar{p}_1,\\
\bar{T}^p_1 \cap \bar{E}_1 &\neq \emptyset.\\
\end{split}
\end{equation*}
\end{description}

\begin{description}
\item 
\parpic[l]{\includegraphics[width=0.3\textwidth]{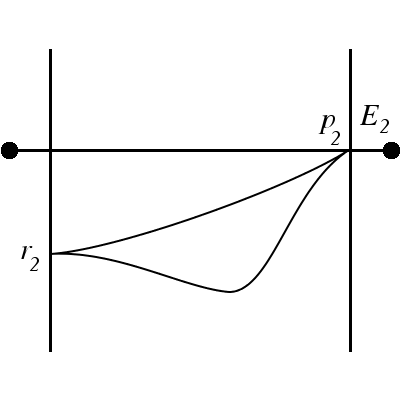}}\;
Blowing down $\bar{T}^r_1$ and $\bar{T}^p_1$ gives the surface $X_2$. On this surface we have fibers $T^r_2$, the strict transform of $\bar{E}_{\hat{q}}$, and $T^p_2$, the strict transform of $\bar{E}_{s_1}$. 
\par Because of the intersections above, $r_2$ and $p_2$ are cusps, both with multiplicity sequence $(2)$. Additionally, $r_2 \notin E_2$, but $p_2 \in E_2$. \par \noi We have the intersections
\begin{equation*}
\begin{split}
\qquad \qquad \qquad \qquad \qquad T^p_2 \cdot C_2 &= 2 \cdot p_2,\\
\qquad \qquad \qquad \qquad \qquad E_2 \cdot C_2 &= 2 \cdot p_2.
\end{split}
\end{equation*}
\end{description}

\pagebreak
\noi{\bf Blowing down $E_2$}
\begin{description}
\item
\parpic[l]{\includegraphics[width=0.3\textwidth]{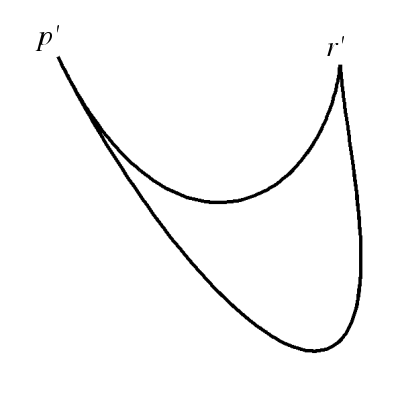}}\;
\noi Blowing down $E_2$ gives a curve $C'$ with two cusps $r'$ and $p'$. The cusp $r'$ is unaltered by the last blowing-down, and it has multiplicity sequence $(2)$. Because of the intersection $$E_2 \cdot C_2=2\cdot p_2,$$ the cusp $p'$ has multiplicity sequence $(2_2)$. 

\end{description}

\vspace{3mm}
\noi To see that $C'$ is a quartic, note that 
\begin{equation*}
\begin{split}
d'&=2 \cdot d - m_q-m_{\hat{q}}-m_{s_1}\\
&=2 \cdot 2 \\
&=4.
\end{split}
\end{equation*}

\noi We may get this curve explicitly by using the Cremona transformation $\psi_2$ to transform the conic $C$, $$C=\mathcal{V}(y^2-2xy+x^2-xz).$$ We get the desired quartic curve $C'$ with an $A_2$- and $A_4$-cusp, $$C'=\mathcal{V}(z^4-xz^3-2xyz^2+x^2y^2).$$

\subsubsection{$C_3$ -- Unicuspidal ramphoid quartic -- $[(2_3)]$}
We construct the unicuspidal ramphoid quartic using a quadratic Cremona transformation with one proper base point.

\begin{description}
\item
\parpic[l]{\includegraphics[width=0.3\textwidth]{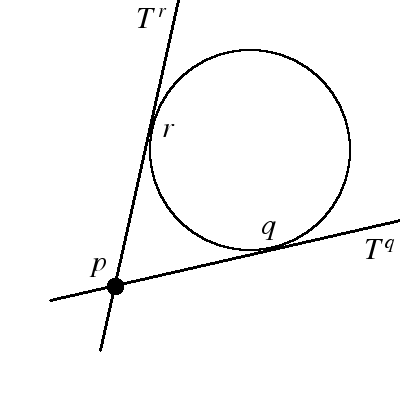}}\;
Let $C$ be an irreducible conic. Choose two arbitrary points $q$ and $r$ on $C$, and let $T^q$ and $T^r$ be the respective tangent lines to $C$ at these points, 
\begin{equation*}
\begin{split}
T^q \cdot C &= 2 \cdot q,\\
T^r \cdot C &= 2 \cdot r.
\end{split}
\end{equation*}
Let $p$ denote the intersection point $T^q \cap T^r$, and note that $p \notin C$.
\end{description}

\noi Applying the transformation $\psi_1(p,T^r,-)$ to $C$, we get the desired quartic.\\

\pagebreak
\noi {\bf Blowing up at $p$}
\begin{description}
\item
\parpic[l]{\includegraphics[width=0.3\textwidth]{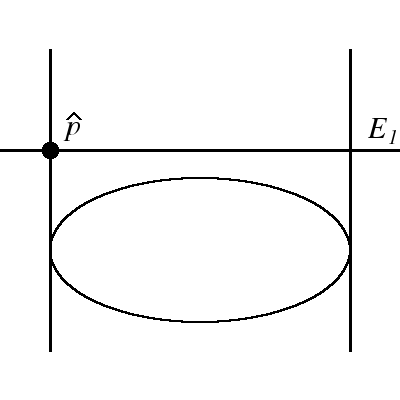}}\;
\noi Blowing up at $p$, we get the ruled surface $X_1$ with horizontal section $E_1$, fibers $T^r_1$ and $T^q_1$, and the transformed curve $C_1$. No points or intersection multiplicities have been affected by this process. Note that since $p \notin C$, we have $E_1 \cap C_1= \emptyset$.

The point $\hat{p}=E_1 \cap T^r_1$ is marked in the figure.
\end{description}

\noi{\bf Elementary transformation in $\hat{p}$}
\begin{description}
\item 
\parpic[l]{\includegraphics[width=0.3\textwidth]{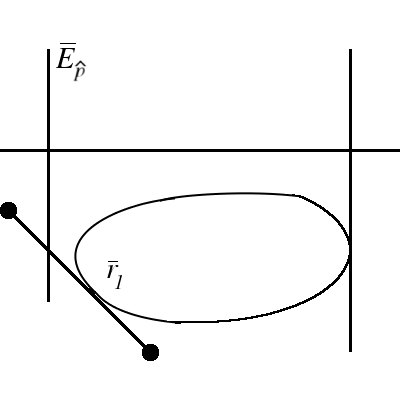}}\;
Blowing up at $\hat{p}$, we get on $\bar{X}_1$ an exceptional line $\bar{E}_{\hat{p}}$. We have the intersections 
\begin{equation*}
\begin{split}
\bar{T}^r_1 \cdot \bar{C}_1 &= 2 \cdot \bar{r}_1,\\
\bar{T}^r_1 \cap \bar{E}_1 &= \emptyset.\\
\end{split}
\end{equation*}
\\
\\
\end{description}

\begin{description}
\item 
\parpic[l]{\includegraphics[width=0.3\textwidth]{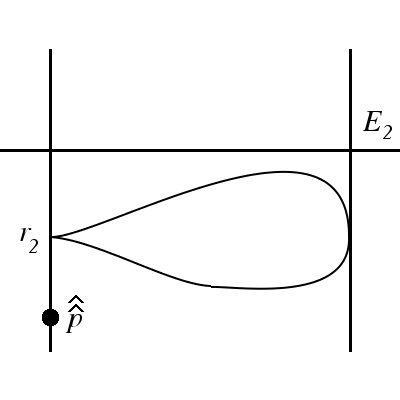}}\;
Blowing down $\bar{T}^r_1$ gives the surface $X_2$. On this surface we have horizontal section $E_2$ and the fiber $T^r_2$, the strict transform of $\bar{E}_{\hat{p}}$. On $T^r_2$ we have marked the point $\hat{\hat{p}}$.
\par Because of the intersection above, $r_2$ is a cusp with multiplicity sequence $(2)$. Note that $r_2 \notin E_2$ and that $E_2^2=-2$.
\\
\end{description}

\noi{\bf Elementary transformation in $\hat{\hat{p}}$}

\begin{description}
\item 
\parpic[l]{\includegraphics[width=0.3\textwidth]{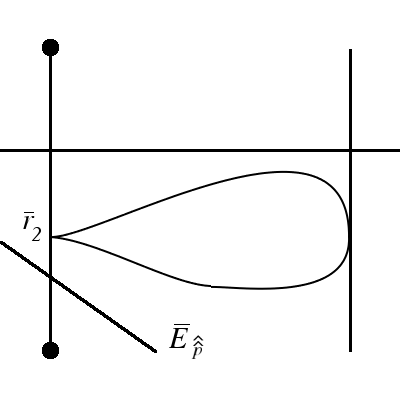}}\;
Blowing up at $\hat{\hat{p}}$ we get on $\bar{X}_2$ an exceptional line $\bar{E}_{\hat{\hat{p}}}$. We have the intersections 
\begin{equation*}
\begin{split}
\bar{T}^r_2 \cdot \bar{C}_2 &= 2 \cdot \bar{r}_2,\\
\bar{T}^r_2 \cap \bar{E}_2 &\neq \emptyset.
\end{split}
\end{equation*}
\\
\\
\end{description}

\pagebreak
\begin{description}
\item 
\parpic[l]{\includegraphics[width=0.3\textwidth]{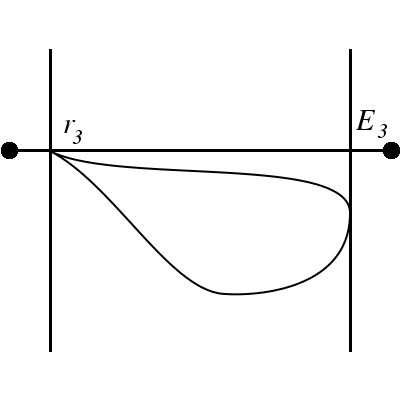}}\;
Blowing down $\bar{T}^r_2$ gives the surface $X_3$. On this surface we have the fiber $T^r_3$, the strict transform of $E_{\hat{\hat{p}}}$. The horizontal section $E_3$ has self-intersection $E_3^2=-1$.
\par Because of the intersection above concerning $\bar{r}_2$, $r_3$ is a cusp with multiplicity sequence $(2_2)$. Since $\bar{T}^r_2 \cap \bar{E}_2 \neq \emptyset$, $r_3 \in E_3$. Furthermore, we have the intersection $$ \qquad \qquad \qquad \qquad \qquad E_3 \cdot C_3  = 2 \cdot r_3.$$
\end{description}

\noi{\bf Blowing down $E_3$}
\begin{description}
\item
\parpic[l]{\includegraphics[width=0.3\textwidth]{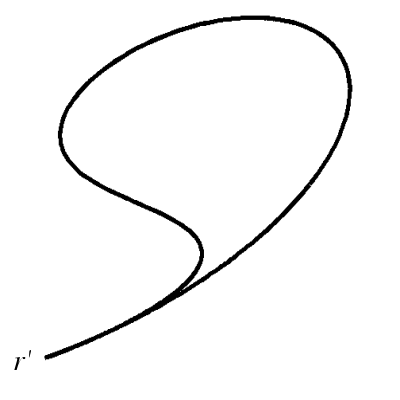}}\;
\noi Blowing down $E_3$ gives a curve $C'$ with one cusp $r'$. Because of the above intersection concerning $r_3$, the cusp $r'$ has multiplicity sequence $(2_3)$. 
\\
\\
\\
\end{description}

\vspace{4mm}
\noi To see that $C'$ is a quartic, note that 
\begin{equation*}
\begin{split}
d'&=2 \cdot d - m_p-m_{\hat{p}}-m_{\hat{\hat{p}}}\\
&=2 \cdot 2\\
&=4.
\end{split}
\end{equation*}

\noi We may construct the curve explicitly by using the Cremona transformation $\psi_1$ to transform the conic $C$, $$C=\mathcal{V}(yz+x^2).$$ We get the desired cuspidal quartic curve $C'$ with an $A_6$-cusp, $$C'=\mathcal{V}(y^4-2xy^2z+yz^3+x^2z^2).$$

\pagebreak
\subsubsection{$C_{4A}$-- Ovoid quartic A -- $[(3)]$} \label{oqa}
We construct the ovoid quartic with one inflection point of type 2 using a quadratic Cremona transformation with two proper base points.

\begin{description}
\item
\parpic[l]{\includegraphics[width=0.3\textwidth]{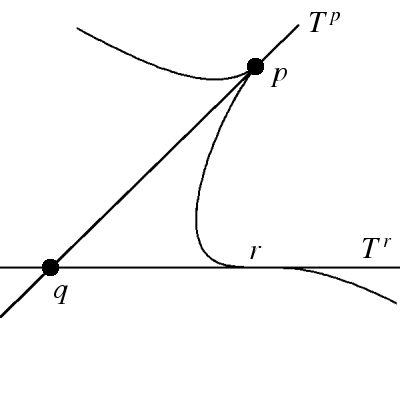}}\;
Let $C$ be the cuspidal cubic with a cusp $p$, $\ol{m}_p=(2)$, and an inflection point $r$ of type 1. The points $p$ and $r$ have tangents $T^p$ and $T^r$, which only intersect $C$ at $p$ and $r$ respectively.
\begin{equation*}
\begin{split}
T^p \cdot C &= 3 \cdot p,\\
T^r \cdot C &= 3 \cdot r. 
\end{split}
\end{equation*}
Let $q$ denote the intersection point $T^p \cap T^r$, and note that $q \notin C$.
\end{description}

\noi Applying the transformation $\psi_2(p,q,T^r)$ to $C$, we get the desired quartic.\\

\noi {\bf Blowing up at $q$}
\begin{description}
\item
\parpic[l]{\includegraphics[width=0.3\textwidth]{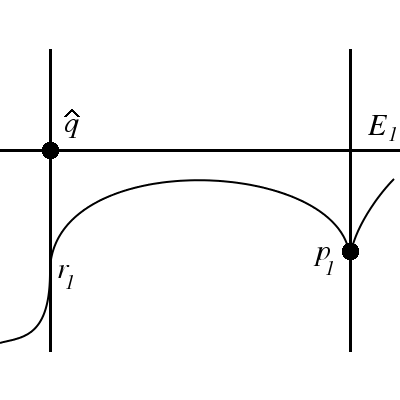}}\;
\noi Blowing up at $q$ we get on $X_1$ a horizontal section $E_1$, fibers $T^p_1$ and $T^r_1$, and the transformed curve $C_1$. No points or intersection multiplicities have been affected by this process. Since $q \notin C$, we have $E_1 \cap C_1= \emptyset$.

The points $\hat{q}=E_1 \cap T^r_1$ and $p_1 \in T^p_1$ are marked in the figure. 
\end{description}

\noi{\bf Elementary transformations in $\hat{q}$ and $p_1$}
\begin{description}
\item 
\parpic[l]{\includegraphics[width=0.3\textwidth]{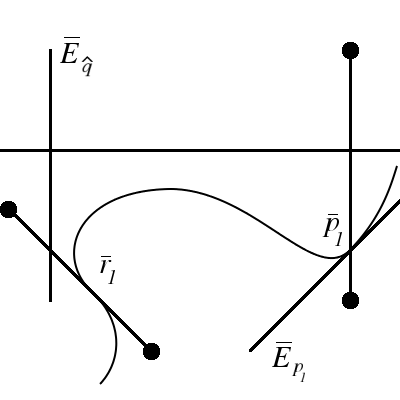}}\;
Blowing up at $\hat{q}$ and the cusp $p_1$ gives on $\bar{X}_1$ two exceptional lines $\bar{E}_{\hat{q}}$ and $\bar{E}_{p_1}$. We now have two smooth points, $\bar{r}_1$ and $\bar{p}_1$, and the intersections 
\begin{equation*}
\begin{split}
\bar{T}^r_1 \cdot \bar{C}_1 &= 3 \cdot \bar{r}_1,\\
\bar{T}^r_1 \cap \bar{E}_1 &= \emptyset,\\
\bar{T}^p_1 \cdot \bar{C}_1 &=((T^p_1 \cdot C_1)_{p_1}-m_{p_1}) \cdot \bar{p}_1\\
&=(3-2) \cdot \bar{p}_1\\
&=1 \cdot \bar{p}_1,\\
\bar{T}^p_1 \cap \bar{E}_1&\neq \emptyset,\\
\bar{E}_{p_1} \cdot \bar{C}_1 &= m_{p_1} \cdot \bar{p}_1\\
&=2 \cdot \bar{p}_1.
\end{split}
\end{equation*}
\end{description}

\pagebreak
\begin{description}
\item 
\parpic[l]{\includegraphics[width=0.3\textwidth]{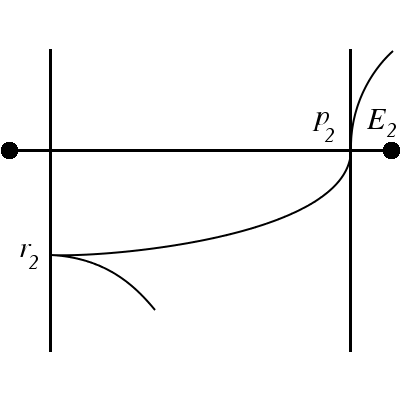}}\;
Blowing down $\bar{T}^r_1$ and $\bar{T}^p_1$ gives the surface $X_2$. On this surface we have $T^r_2$, the strict transform of $\bar{E}_{\hat{q}}$, and $T^p_2$, the strict transform of $\bar{E}_{p_1}$. 
\par Since $\bar{E}_1 \cap \bar{T}^p_1\neq \emptyset$, $p_2 \in E_2$. Moreover, since $\bar{T}^r_1 \cap \bar{E}_1 = \emptyset$, $r_2 \notin E_2$. While $p_2$ is a smooth point, $r_2$ is a cusp with multiplicity sequence $(3)$. The latter is a consequence of \begin{equation*}\qquad \qquad \qquad \qquad \qquad \bar{T}^r_1 \cdot \bar{C}_1 = 3 \cdot \bar{r}_1.\end{equation*}
\end{description}

\begin{description}
\item We also have the important intersections
\begin{equation*}
\begin{split}
T^p_2 \cdot C_2 &= ((\bar{T}^p_1 \cdot \bar{C}_1)_{\bar{p}_1}+(\bar{E}_{p_1} \cdot \bar{C}_1)_{\bar{p}_1})\cdot p_2\\
&=(1+2) \cdot p_2\\
&=3 \cdot p_2,\\
E_2 \cdot C_2 &= (\bar{T}^p_1 \cdot \bar{C}_1)_{\bar{p}_1} \cdot p_2 \\
&=1 \cdot p_2.
\end{split}
\end{equation*}
\end{description}

\noi{\bf Blowing down $E_2$}
\begin{description}
\item
\parpic[l]{\includegraphics[width=0.3\textwidth]{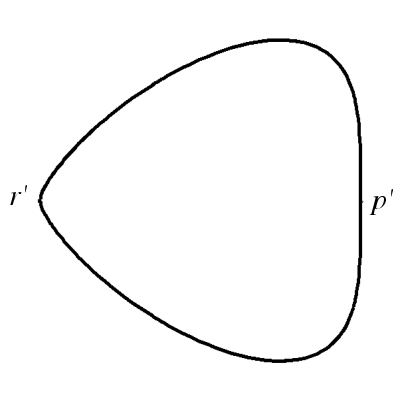}}\;
\noi Blowing down $E_2$ gives a curve $C'$ with one cusp $r'$. The cusp is unaltered by the last blowing-down, and it has multiplicity \linebreak sequence $(3)$. Because of the above intersections concerning $p_2$, $C$ has an inflection point $p'$ of type $2$. 
\\
\\
\end{description}

\noi To see that $C'$ is a quartic, note that 
\begin{equation*}
\begin{split}
d'&=2 \cdot d - m_q-m_{\hat{q}}-m_{p_1}\\
&=2 \cdot 3-2 \\
&=4.
\end{split}
\end{equation*}

\noi We may construct the curve explicitly by using the Cremona transformation $\psi_2$ to transform the cubic $C$, $$C=\mathcal{V}(xz^2-y^3).$$ We get the desired quartic curve $C'$ with one $E_6$-cusp and one inflection point of type 2, $$C'=\mathcal{V}(z^4-xy^3).$$

\pagebreak
\subsubsection{$C_{4B}$ -- Ovoid quartic B -- $[(3)]$}\label{oqb}
We construct the ovoid quartic with two inflection points using a quadratic Cremona transformation with two proper base points. By direct calculation and investigation of the curve and its orientation with respect to the Cremona transformation, we know that this construction is valid. There are, however, some unanswered questions concerning the number of inflection points constructed with Cremona transformations. This issue will be discussed in Section \ref{oninflection}.

\begin{description}
\item
\parpic[l]{\includegraphics[width=0.3\textwidth]{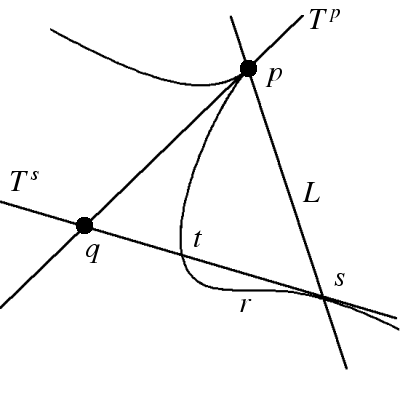}}\;
Let $C$ be the cuspidal cubic with a cusp $p$, $\ol{m}_p=(2)$, and an inflection point $r$. Choose a smooth point $s$, $s \neq r$. The tangents $T^p$ and $T^s$ intersect $C$,
\begin{equation*}
\begin{split}
T^p \cdot C &= 3 \cdot p,\\
T^s \cdot C &= 2 \cdot s +1 \cdot t,
\end{split}
\end{equation*}
where $t$ is a smooth point. 

Additionally, we denote by $L=L^{ps}$ the line through the cusp $p$ and the point $s$. Intersecting $L$ with $C$ gives $$L \cdot C = 2\cdot p+1 \cdot s.$$ 
Finally, let $q$ denote the intersection point $T^p \cap T^s$.
\end{description}

\noi Applying the transformation $\psi_2(q,p,L)$ to $C$, we get the desired quartic. Note that the inflection point $r$ will not be affected by this process.\\

\noi {\bf Blowing up at $p$}
\begin{description}
\item
\parpic[l]{\includegraphics[width=0.3\textwidth]{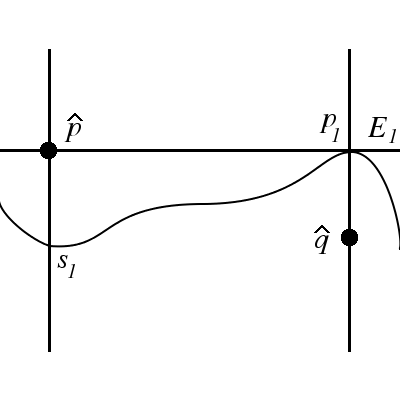}}\;
\noi Blowing up at the cusp $p$, we get $X_1$ with horizontal section $E_1$, fibers $L_1$ and $T^p_1$, and the transformed curve $C_1$. The points $\hat{p}=E_1 \cap L_1$ and $\hat{q}$ are marked in the figure. Since we blew up at $p$ with multiplicity sequence $(2)$, $p_1$ is smooth. Additionally,
\begin{equation*}
\begin{split}
E_1 \cdot C_1 &= 2 \cdot p_1,\\
T^p_1 \cdot C_1 &=((T^p \cdot C)_{p}-m_{p}) \cdot p_1\\
&=(3-2) \cdot p_1\\
&=1 \cdot p_1,\\
L_1 \cdot C_1&=1 \cdot s_1.\\
\end{split}
\end{equation*}
\end{description}

\pagebreak
\noi{\bf Elementary transformations in $\hat{p}$ and $\hat{q}$}
\begin{description}
\item 
\parpic[l]{\includegraphics[width=0.3\textwidth]{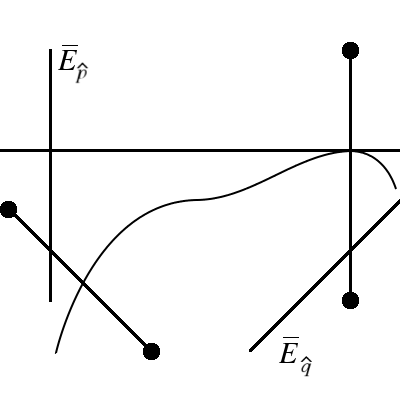}}\;
Blowing up at $\hat{p}$ and $\hat{q}$ gives on $\bar{X}_1$ two exceptional lines $\bar{E}_{\hat{p}}$ and $\bar{E}_{\hat{q}}$, which do not change any points or intersections.
\\
\\
\\
\\
\end{description}

\begin{description}
\item 
\parpic[l]{\includegraphics[width=0.3\textwidth]{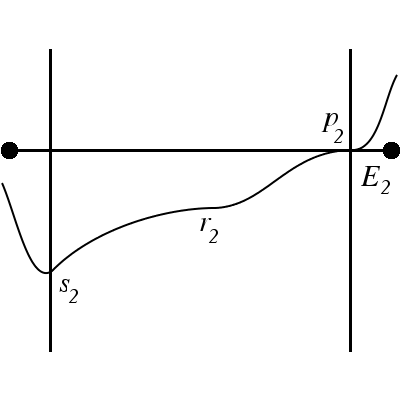}}\;
Blowing down $\bar{L}_1$ and $\bar{T}^p_1$ gives the surface $X_2$. On this surface we have the fibers $L_2$, the strict transform of $\bar{E}_{\hat{p}}$, and $T^p_2$, the strict transform of $\bar{E}_{\hat{q}}$. 
\par We get the intersection multiplicity 
\begin{equation*}
\begin{split}
E_2 \cdot C_2 &= ((\bar{E}_1 \cdot \bar{C}_1)_{\bar{p}_1}+(\bar{T}^p_1 \cdot \bar{C}_1)_{\bar{p}_1}) \cdot p_2\\
&=(2+1)\cdot p_2\\
&=3 \cdot p_2.
\end{split}
\end{equation*}
Additionally, note that we have affected the point $s_2$ with this elementary transformation. On $\bar{X}_1$, assuming $\bar{s}_1$ has a tangent $\bar{T}$ and then blowing down $\bar{L}_1$, we calculate the intersection multiplicity of the tangent $T^{s_2}$ to $C_2$ at $s_2$ and $C_2$,
\begin{equation*}
\begin{split}
(T^{s_2} \cdot C_2)_{s_2} &=((\bar{T} \cdot \bar{C}_1)_{\bar{s}_1}+(\bar{L}_{1} \cdot \bar{C}_1)_{\bar{s}_1})\\
&=(2+1)\\
&=3.
\end{split}
\end{equation*}
Hence, $s_2$ is an inflection point of $C_2$.
\end{description}

\noi{\bf Blowing down $E_2$}
\begin{description}
\item
\parpic[l]{\includegraphics[width=0.3\textwidth]{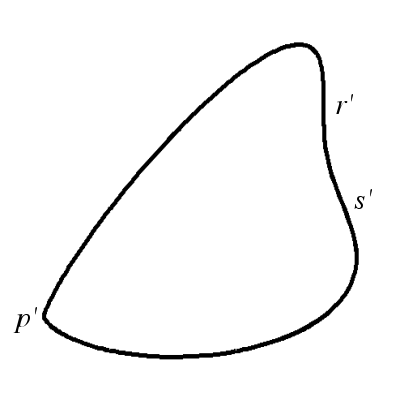}}\;
\noi Blowing down $E_2$ gives a curve $C'$ with one cusp $p'$ and two inflection points $r'$ and $s'$. Because of the intersection $E_2 \cdot C_2=3\cdot p_2$, the cusp $p'$ has multiplicity sequence $(3)$. 
\\
\\
\\
\end{description}

\vspace{2mm}
\noi To see that $C'$ is a quartic, observe that 
\begin{equation*}
\begin{split}
d'&=2 \cdot d - m_p-m_{\hat{p}}-m_{\hat{q}}\\
&=2 \cdot 3-2 \\
&=4.
\end{split}
\end{equation*}

\noi We may construct this curve explicitly by using the Cremona transformation $\psi_2$ to transform the cubic $C$, $$C=\mathcal{V}(yz^2-x^2z+x^3).$$ We get the desired quartic curve $C'$ with one $E_6$-cusp and two inflection points of type 1, $$C'=\mathcal{V}(x^3y-z^3x+z^4).$$

\section{A note on inflection points}\label{oninflection}
We observe that if two curves have the same cuspidal configuration, but a different number of inflection points, then they can not be projectively equivalent. We are therefore interested in the configuration of inflection points of a curve, in addition to the cuspidal configuration. In this section we give examples on how the configuration of inflection points can be determined by other properties of the curve. Furthermore, we give examples on how inflection points behave unpredictably under Cremona transformations.

\subsubsection{Configuration of inflection points}
The ovoid quartic curves are examples of curves which are not projectively equivalent, even though they have the same cuspidal configuration. The most prominent difference between the two curves is the configuration of inflection points. 

Fortunately, when $C$ is a curve of degree $d$ with a given cuspidal configuration, we can often determine which of the possible configurations of inflection points actually exist.



\begin{ex}
\noi Consider a cuspidal quartic curve $C$ with cusps $p$ and $q$ with multiplicity sequences $\ol{m}_p=(2_2)$ and $\ol{m}_q=(2)$. Since $C$ is a quartic, the intersection multiplicity of the curve and the tangents at the cusps must be $r_p=4$ and $r_q=3$. A quartic curve with this cuspidal configuration must, by Theorem \ref{Ghessian}, have precisely one inflection point of type 1. 
\begin{equation*}
\begin{split}
s&=3\cdot4\cdot(4-2)-6 \cdot (2+1)-(2+4-3)-(2+3-3)\\
&=1.
\end{split}
\end{equation*}

\end{ex}

\begin{ex}
Consider a quartic curve with one cusp $p$ with multiplicity sequence $(2_3)$. Since $C$ is a quartic, the intersection multiplicity of the curve and the tangent at the cusp is $r_p=4$. We want to determine the configuration of inflection points of this curve. \medskip \\
\noi The curve must have three inflection points, counted properly.
\begin{equation*}
\begin{split}
s&=3\cdot4\cdot(4-2)-6 \cdot 3- (2+4-3)\\
&=3.
\end{split}
\end{equation*}

\noi We show that the only existing curve with this cuspidal configuration has three inflection points of type 1. We will do this by excluding all other configurations of inflection points and checking that the mentioned configuration does not violate any restrictions. \medskip \\
\noi The dual curve $C^*$ has degree
\begin{equation*}
\begin{split}
d^*=&\;2 \cdot d-2-(m_p-1)\\
=&\;2 \cdot 4-2-1\\
=&\;5.
\end{split}
\end{equation*}
\noi The curve $C$ can at the cusp $p$ be given by the Puiseux parametrization $$(C,p)=(t^2:c_4t^4+c_7t^7+\ldots:1).$$ The dual curve $C^*$ can at the dual point $p^*$ be given by the Puiseux \linebreak parametrization $$(C^*,p^*)=(a^*t^2+\ldots:1:c_4^*t^4+c^*_7t^7\ldots).$$ Hence, $p^*$ is a cusp, $\ol{m}^*_p=(2_3)$, and $r_p^*=4$. \\ 

\noi Assume that $C$ only has one inflection point $q$ of type 3. Then the curve must locally around $q$ be parametrized by $$(C,q)=(t:c_5t^5+\ldots:1).$$ In this parametrization of the quartic we have that $r_q = 5$, but that is a violation of B\'{e}zout's theorem since $r_q=5>4=d$. Hence, $C$ can not have one inflection point of type 3.\\


\noi Assume that $C$ has two inflection points $q$ and $r$ of type 1 and 2, respectively. Then we have 
\begin{equation*}
\begin{split}
(C,q)&=(t:c_3t^3+\ldots:1),\\
(C,r)&=(t:c_4t^4+\ldots:1).\\
\end{split}
\end{equation*}
Moreover, 
\begin{equation*}
\begin{split}
(C^*,q^*)&=(a^*t^2+\ldots:1:c_3^*t^3+\ldots),\\
(C^*,r^*)&=(a^*t^3+\ldots:1:c_4^*t^4+\ldots),
\end{split}
\end{equation*}
implying that $q^*$ and $r^*$ are cusps with multiplicity sequences $(2)$ and $(3)$. Then the dual curve $C^*$ has three singularities, and $\sum_{j=p^*,q^*,r^*} \delta_j = 7$. This contradicts the genus formula (\ref{thm:genus}). Hence, $C$ can not have two inflection points. \\

\noi Assume that $C$ has three inflection points $q$, $r$ and $u$ of type 1. Then the dual quintic curve $C^*$ has four cusps $p^*$, $q^*$, $r^*$ and $u^*$ with multiplicity sequences $(2_3)$, $(2)$, $(2)$ and $(2)$, respectively. Indeed, $C^*$ does not contradict the genus formula,
$$\sum_{j=p^*,q^*,r^*,u^*} \delta_j = 6 = \frac{(5-1)(5-2)}{2}.$$
Additionally, we do not get any contradiction when we calculate the number of inflection points on $C^*$.
\begin{equation*}
\begin{split}
s^*&=3\cdot5\cdot(5-2)-6 \cdot 6- (2 +4-3)-3\cdot(2+3-3)\\
&=0.
\end{split}
\end{equation*}
We shall later see that $C^*$ is a remarkable curve.
\end{ex}

\subsubsection{Inflection points and Cremona transformations}
\noi Constructing inflection points with Cremona transformations is a bit more subtle than constructing cusps. The constructions of the cuspidal quartics give examples of this subtlety.\medskip

\begin{ex}
To construct the ovoid quartic with one inflection point on page \pageref{oqa}, we use a cuspidal cubic with one inflection point of type 1. In the transformation we transform the inflection point of the cubic into the cusp of the quartic. Furthermore, the cusp of the cubic is transformed into the inflection point of type 2 of the quartic. Moreover, the construction of the inflection point is directly visible in the action of the Cremona transformation.\medskip
\end{ex}

\begin{ex}
To construct the ovoid quartic with two inflection points on page \pageref{oqb}, we also use a cuspidal cubic with one inflection point of type 1. The cusp of the cubic is transformed into the cusp of the quartic. The inflection point on the cubic is unaffected by the Cremona transformation. The construction of the second inflection point, however, is not directly visible in the pictures of the Cremona transformation. We carefully argued for the construction of this inflection point based on elementary properties of monoidal transformations.\medskip
\end{ex}

\begin{ex}
To construct the bicuspidal quartic on page \pageref{bcq}, we transform a conic. We make the intriguing observation that there is no apparent reason why the Cremona transformation of the conic should give this curve an inflection point. We have, however, seen that the curve must have an inflection point earlier in this section.\medskip
\end{ex}

\noi The three examples clearly show the unpredictable behavior of the inflection points in the construction of curves by Cremona transformations. This is of course of little relevance when it comes to the construction of a particular cuspidal configuration of a curve of given degree. However, it is essential to note this unpredictable affair in the construction of all, up to projective equivalence, rational cuspidal curves of a given degree.

\section{The Coolidge--Nagata problem}
If a rational curve can be transformed into a line by successive Cremona transformations, then we call the curve {\em rectifiable}. The following problem was introduced by Coolidge in \cite[pp.396--399]{Coolidge}, but because of \cite{Nagata} it also has the name of Nagata.\\

\noi {\bf The Coolidge--Nagata problem.} {\em Which rational curves can birationally be transformed into a line? In particular, which rational cuspidal curves have this property?}\\

\noi Coolidge noted that most rational curves of low degree are rectifiable. There exist, however, rational curves which are {\em not} rectifiable. An example of such a curve is the sextic with ten singularities, each having two distinct tangents. 

There are actually no known rational {\em cuspidal} curves which are {\em not} rectifiable. In particular, all rational cuspidal curves we will encounter in this thesis are rectifiable by construction. This claim can easily be verified for the curves we have encountered this far. All these curves are somehow constructed from an irreducible conic. Hence, applying the inverse transformation to any of the cuspidal curves will produce the conic. Moreover, any irreducible conic can be transformed into a line by applying a quadratic Cremona transformation with three base points, just choose the base points on the conic.

The complement of a rational cuspidal curve in $\mathbb{P}^2$ is an open surface. For an open surface $\mathbb{P}^2 \setminus C$, where $C$ is a rational cuspidal curve, there exists assumptions with which it is possible to prove that $C$ is rectifiable. For a detailed discussion of these results, see \cite[pp.418--419]{Bobadillains}.

Matsuoka and Sakai observed that if a curve is transformable to a line by Cremona transformations, then the Matsuoka--Sakai inequality on page \pageref{MatsuokaSakai} is a consequence of one of Coolidge's results. Since they were able to independently prove that the inequality is valid for all rational cuspidal curves, they proposed the following conjecture \cite[p.234]{MatsuokaSakai}.

\begin{conj}
Every rational cuspidal curve can be transformed into a line by a Cremona transformation.
\end{conj}

\chapter{Rational cuspidal quintics}\label{rcq}
In this chapter we will describe the rational cuspidal quintic curves. Using the method provided by Flenner and Zaidenberg in \cite{FlZa97} and methods described in Chapter \ref{TB}, we will obtain a list of possible cuspidal configurations for a rational quintic. We will construct and describe examples of all these curves.

\section{The cuspidal configurations}
We want to find all possible cuspidal configurations for rational cuspidal quintics. The idea is the same as for lower degrees, only now getting the results requires more theoretical background. 

Flenner and Zaidenberg find in \cite[pp.104--105]{FlZa97} all possible cuspidal configurations for quintics with three or more cusps. Unfortunately, there are some incomplete explanations and one minor error in their work. In the following we will generalize the ideas of Flenner and Zaidenberg's proof. We will obtain a list of cuspidal configurations for all cuspidal quintics by using general restrictions and Flenner and Zaidenberg's ideas to exclude the nonexistent ones. Included in this work are Flenner and Zaidenberg's results, with complementary additions and appropriate corrections. 

\subsection{General restrictions}
\noi Assume that we have a rational cuspidal curve of degree $d=5$ with a finite number of cusps $p$, $q$, $r$, \ldots. Let the cusp $p$ have multiplicity sequence $\ol{m}_{p}=(m_p, m_{1},...,m_{n_p})$ and delta invariant $\delta_p$. Define the same invariants for the other cusps. By Theorem \ref{thm:genus} the cusps of $C$ satisfy \begin{equation*}\frac{(5-1)(5-2)}{2}=6=\sum_{p\in \mathrm{Sing}\,C}\sum_{i=0}^{n_p}\frac{m_i(m_{i}-1)}{2}.\label{5:genus}\end{equation*}\medskip


\noi We sort the cusps by their multiplicity such that $m_p \geq m_q \geq m_r \geq \ldots$. Then the largest multiplicity of the cusps on the curve is $\mu=m_{p}$. By B\'{e}zout's theorem we immediately get $$\mu \leq 4.$$ If there are two or more cusps on the curve, then we can improve this restriction. Let $L$ be a line which intersects the rational cuspidal curve $C$ in two of its cusps. Then the sum of the intersection multiplicities of $L$ and $C$ in these two points is less than or equal to the degree of the curve. For any pair of cusps, in particular the largest cusps $p$ and $q$, $$(C \cdot L)_{p}+(C \cdot L)_{q} \leq 5.$$ Since the minimal intersection multiplicity between a line and any cusp is equal to the minimal multiplicity of any cusp, i.e., $2$, $m_{p}$ can not be greater than 3. Hence, for curves with two or more cusps, we have the restriction $$\mu \leq 3.$$

\noi Further restrictions on the multiplicities can be found using the results in Chapter \ref{TB}, and we will refer to these in the following. 


\subsection{One cusp}
Assuming we only have one cusp $p$, the above restrictions leave us with four cuspidal configurations, of which two can be excluded. 

\subsubsection{Maximal multiplicity 4}
Let $\mu=4$. Then $\delta=6$. Hence, $(4)$ is the only possible multiplicity sequence. Curves with this multiplicity sequence exist.

\subsubsection{Maximal multiplicity 3}
Let $\mu=3$. By elementary properties of the multiplicity sequence, $m_1\leq3$. Furthermore, the only multiplicity sequences satisfying the genus formula are $(3_2)$ and $(3,2_3)$. Neither of these multiplicity sequences are possible on a rational cuspidal quintic curve. The sequence $(3_2)$ violates $m_0+m_1 \leq 5$ from (\ref{upperandlower}) on page \pageref{upperandlower}. The sequence $(3,2_3)$ does not satisfy the property of the multiplicity sequence given in Proposition \ref{multiseq} on page \pageref{multiseq}.

\subsubsection{Maximal multiplicity 2}
Let $\mu=2$. Then the only possible multiplicity sequence is $(2_6)$. Curves with this multiplicity sequence exist.

\subsubsection{Conclusion} 
We have two possible cuspidal configurations for unicuspidal quintics.
\begin{center}
	\setlength{\extrarowheight}{2pt}
\begin{tabular}{cc}
\hline
{\bf Curve} & {\bf Cuspidal configuration}\\
\hline
$C_1$ & $(4)$ \\
$C_2$ & $(2_6)$\\
\hline
\end{tabular}
\end{center}

\subsection{Two cusps}
Assuming we have two cusps $p$ and $q$, the above restrictions leave us with six cuspidal configurations, of which three can be excluded. 

\subsubsection{Maximal multiplicity 3}
\noi Let $\mu=m_p=3$. Then $m_q = 2$, and we have three cuspidal configurations satisfying the genus formula, $$[(3,2_2),(2)] \qquad [(3,2),(2_2)] \qquad [(3),(2_3)].$$ The configuration $[(3,2_2),(2)]$ violates properties of the multiplicity sequence given in Proposition \ref{multiseq}. The two remaining configurations actually exist.

\subsubsection{Maximal multiplicity 2}
\noi Let $\mu=m_p=2$. Then $m_q = 2$, and we have three cuspidal configurations satisfying the genus formula, $$[(2_5),(2)] \qquad [(2_4),(2_2)] \qquad [(2_3),(2_3)].$$ The first and last configurations can be excluded by a similar argument involving Cremona transformations.

\subsubsection{\bf Curve -- $[(2_5),(2)]$}
\begin{description}
\item
\parpic[l]{\includegraphics[width=0.3\textwidth]{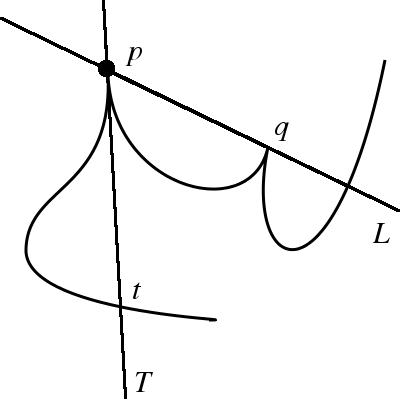}}\;
Assume that there exists a rational cuspidal quintic $C$ with two cusps $p$ and $q$ with multiplicity sequences $\ol{m}_p=(2_5)$ and $\ol{m}_q=(2)$. Let $T$ denote the tangent line of $p$ and let $L$ be the line between the two cusps.

Since $C$ is a quintic and $p$ has multiplicity sequence $(2_5)$, Lemma \ref{intmult} implies that $$\qquad \qquad \qquad \qquad \qquad(C \cdot T)_p = 4.$$ By B\'{e}zout's theorem, $T$ must intersect $C$ transversally in a smooth point $t$, $$C \cdot T = 4 \cdot p + 1 \cdot t.$$

\end{description}

\noindent Applying the Cremona transformation $\psi(p,T,-)$ to $C$ gives a quartic curve. We shall obtain a contradiction by observing that this particular quartic curve does not exist. This implies that the quintic curve $C$ can not exist either. 

Note that the Cremona transformation does not affect the $A_2$-cusp $q$, and we will therefore keep this point out of the discussion.\\

\noi {\bf Blowing up at $p$}
\begin{description}
\item \parpic[l]{\includegraphics[width=0.3\textwidth]{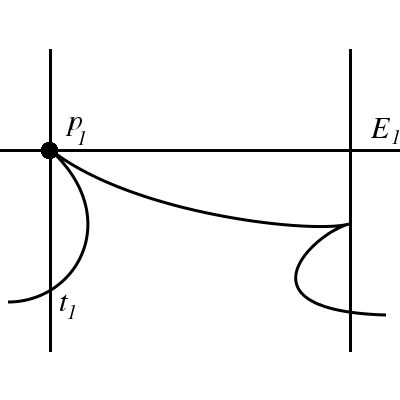}}\;
Blowing up at $p$, we get the ruled surface $X_1$, where the transform $p_1$ of $p$ is a cusp with multiplicity sequence $(2_4)$. We also have the intersections 
\begin{equation*}
\begin{split}
C_1 \cdot E_1 &= m_{p} \cdot p_1\\
&= 2 \cdot p_1,\\
C_1 \cdot T_1 &= ((C \cdot T)_p-m_p)\cdot p_1+ 1 \cdot t_1\\
&=2 \cdot p_1 + 1 \cdot t_1.
\end{split}
\end{equation*} 
\end{description}

\noi {\bf Elementary transformation in $p_1$}
\begin{description}
\item \parpic[l]{\includegraphics[width=0.3\textwidth]{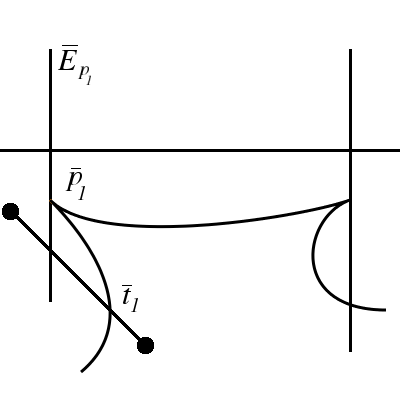}}\;
Blowing up at $p_1$ gives the cusp $\bar{p}_1$ with multiplicity sequence $(2_3)$ on the surface $\bar{X}_1$. The exceptional line $\bar{E}_{p_1}$ separates $E_1$ and $T_1$. We get the intersections 
\begin{equation*}
\begin{split}
\bar{E}_1 \cap \bar{C}_1 &= \emptyset,\\
\bar{E}_{p_1} \cdot \bar{C}_1 &= m_{p_1} \cdot \bar{p}_1 \\
&= 2 \cdot \bar{p}_1,\\
\bar{T}_1 \cdot \bar{C}_1 &= 1 \cdot \bar{t}_1,\\
\bar{T}_1 \cap \bar{E}_1 &=\emptyset.
\end{split}
\end{equation*}
\end{description}

\begin{description}
\item \parpic[l]{\includegraphics[width=0.3\textwidth]{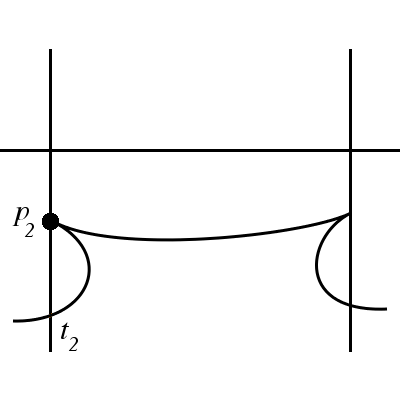}}\;
Blowing down $\bar{T}_1$ gives the surface $X_2$. The transform $T_2$ of $\bar{E}_{p_1}$ intersects $C_2$ in the $A_6$-cusp $p_2$ and the smooth point $t_2$. In particular, $t_2 \notin E_2$ because of the last intersection above. We have the intersection $$T_2 \cdot C_2 = 2 \cdot p_2 + 1\cdot t_2.$$ 
\\
\\
\end{description}

\pagebreak
\noi{\bf Elementary transformation in $p_2$}
\begin{description}
\item \parpic[l]{\includegraphics[width=0.3\textwidth]{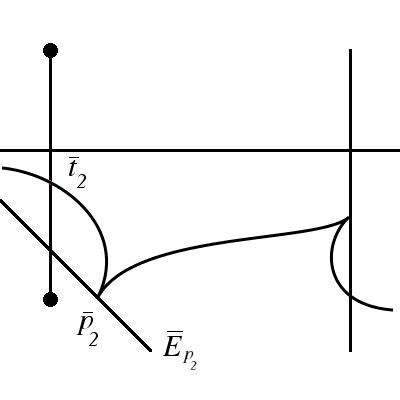}}\;
Blowing up at $p_2$ gives the cusp $\bar{p}_2$ with multiplicity sequence $(2_2)$ on the surface $\bar{X}_2$. Note that the exceptional line $\bar{E}_{p_2}$ does not separate $E_2$ and $T_2$, but we have the intersections 
\begin{equation*}
\begin{split}
\bar{E}_{p_2} \cdot \bar{C}_2 &= 2 \cdot \bar{p}_2,\\
\bar{T}_2 \cdot \bar{C}_2 &= 1 \cdot \bar{t}_2.
\end{split}
\end{equation*}
\\

\item \parpic[l]{\includegraphics[width=0.3\textwidth]{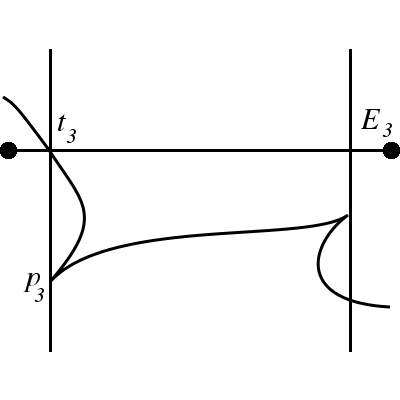}}\;
Blowing down $\bar{T}_2$ gives the surface $X_3$ with the $A_4$-cusp $p_3$, and the intersections
\begin{equation*}
\begin{split}
E_3 \cdot C_3 &= 1 \cdot t_3,\\
T_3 \cdot C_3 &= 1 \cdot t_3 + 2 \cdot p_3.
\end{split}
\end{equation*}
Note that both $E_3$ and $T_3$ intersect $C_3$ \linebreak transversally at $t_3$.
\\
\end{description}

\noi {\bf Blowing down $E_3$}\\
\begin{description}
\item \parpic[l]{\includegraphics[width=0.3\textwidth]{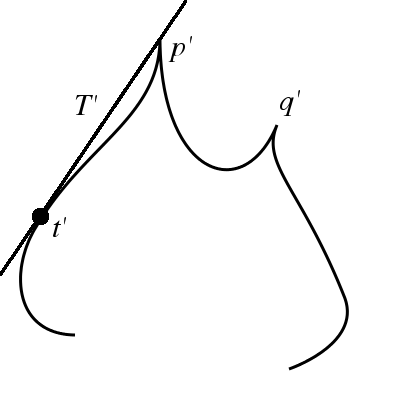}}\;
\noi Since neither of the cusps $p_3$ and $q_3$ lie on $E_3$, blowing down $E_3$ gives a curve $C'$ with one $A_4$-cusp $p'$  and one $A_2$-cusp $q'$. Since both $E_3$ and $T_3$ intersect $C_3$ transversally in the point $t_3$, this point is transformed into a smooth point $t' \in C'$, with the property that the tangent line at $t'$ will be the transform $T'$ of $T_3$. This is a result of the intersection multiplicity 
\begin{equation*}
\begin{split}
(T' \cdot C')_{t'}&=(T_3\cdot C_3)_{t_3} + (E_3 \cdot C_3)_{t_3}\\
&=2.
\end{split}
\end{equation*}
\noi Note that by the intersection $T_3 \cdot C_3$ above, $T'$ must intersect $C'$ at the $A_4$-cusp $p'$. We have $$T' \cdot C' = 2 \cdot t' + 2 \cdot p'.$$
\end{description}

\noi Observe that $C'$ is a quartic since 
\begin{equation*}
\begin{split}
d'&=2 \cdot d - m_p-m_{p_1}-m_{p_2}\\
&=2 \cdot 5 - 2-2-2 \\
&=4.
\end{split}
\end{equation*}

\pagebreak
\noindent Up to projective equivalence, there is only one quartic curve $D$ with cuspidal configuration $[(2_2),(2)]$. Let $D$ be the zero set of $F=z^4-xz^3-2xyz^2+x^2y^2$. $D$ has an $A_4$-cusp at $p=(0:1:0)$. The polar of $C$ at $p$ is given by $$P_pC=\V(-2xz^2+2x^2y).$$ Calculations in {\em Maple} reveal that the polar of $D$ at $p$ does not intersect $D$ at any smooth point. This implies that there is {\em no} point $t\in D$ such that the tangent $T_{t}D$ intersects $D$ at $p$. Hence, neither $C'$ nor the quintic curve with multiplicity sequence $[(2_5),(2)]$ exist.\\

\subsubsection{\bf Curve -- $[(2_3),(2_3)]$}

\begin{description}
\item
\parpic[l]{\includegraphics[width=0.3\textwidth]{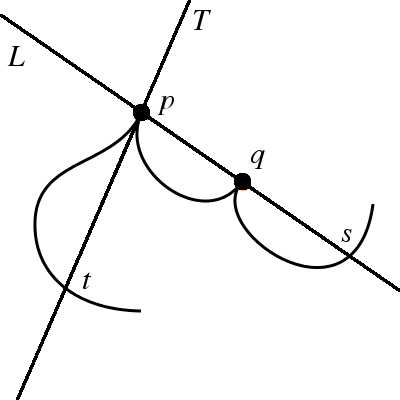}}\;
Assume that there exists a rational cuspidal quintic $C$ with two cusps $p$ and $q$ with multiplicity sequences $\ol{m}_p=(2_3)$ and $\ol{m}_q=(2_3)$. Let $T=T^p$ be the tangent line of $p$ and let $L=L^{pq}$ be the line between the two cusps. Let $s$ and $t$ be smooth points of $C$. Then we have
\begin{equation*}
\begin{split}
\qquad \qquad \qquad \qquad \qquad \qquad C \cdot T &= 4 \cdot p + 1 \cdot t,\\
C \cdot L &=2\cdot p+2 \cdot q+1 \cdot s.
\end{split}
\end{equation*}
\end{description}

\noi We proceed with the Cremona transformation $\psi_2(q,p,T)$. This will lead to a contradiction.\\

\noi {\bf Blowing up at $p$}

\begin{description}
\item
\parpic[l]{\includegraphics[width=0.3\textwidth]{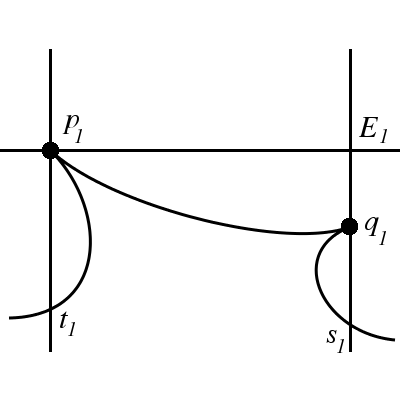}}\;
Blowing up at the cusp $p$ with multiplicity sequence $(2_3)$ gives a ruled surface $X_1$ with exceptional line $E_1$. The transform $p_1$ of $p$ is a cusp with multiplicity sequence $(2_2)$. We have
\begin{equation*}
\begin{split}
(C_1 \cdot T_1)_{p_1}&=(C \cdot T)_p-m_{p}\\
&=4-2\\
&=2,\\
E_1 \cdot C_1 &= m_{p} \cdot p_1\\
&=2 \cdot p_1.
\end{split}
\end{equation*}
\end{description}
\begin{description}
\item Additionally, with $q_1$ still an $A_6$-cusp, we have the intersections 
\begin{equation*}
\begin{split}
C_1 \cdot T_1&=2 \cdot p_1+ 1 \cdot t_1,\\ 
C_1 \cdot L_1&=2 \cdot q_1 + 1 \cdot s_1.
\end{split}
\end{equation*}
\end{description}

\pagebreak
\noi {\bf Elementary transformations in $p_1$ and $q_1$}
\begin{description}
\item
\parpic[l]{\includegraphics[width=0.3\textwidth]{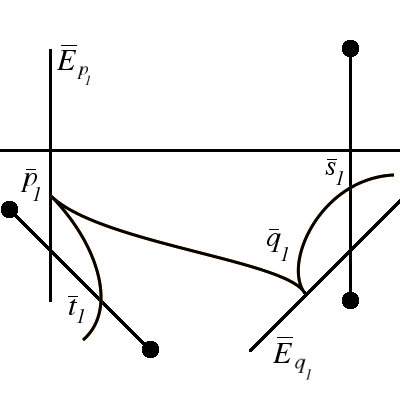}}\;
Blowing up at $p_1$ and $q_1$, we get on the surface $\bar{X}_1$ the exceptional lines $\bar{E}_{p_1}$ and $\bar{E}_{q_1}$. We get one $A_2$-cusp $\bar{p}_1$ and one $A_4$-cusp $\bar{q}_1$. Additionally, $\bar{E}_{p_1}$ separates $T_1$ and $E_1$, while $\bar{E}_{q_1}$ does not separate $L_1$ and $E_1$. We have the intersections 
\begin{equation*}
\begin{split}
\bar{E}_1 \cap \bar{C}_1&= \emptyset, \\
\bar{T}_1 \cap \bar{E}_1 &= \emptyset,\\
\bar{L}_1 \cap \bar{E}_1 &\neq \emptyset,\\
\bar{E}_{p_1} \cdot \bar{C}_1 &= 2 \cdot \bar{p}_1,\\
\bar{T}_1 \cdot \bar{C}_1&= 1 \cdot \bar{t}_1,\\
\bar{E}_{q_1} \cdot \bar{C}_1&=2 \cdot \bar{q}_1,\\
\bar{L}_1 \cdot \bar{C}_1 &=1 \cdot \bar{s}_1.
\end{split}
\end{equation*}
\end{description}

\begin{description}
\item
\parpic[l]{\includegraphics[width=0.3\textwidth]{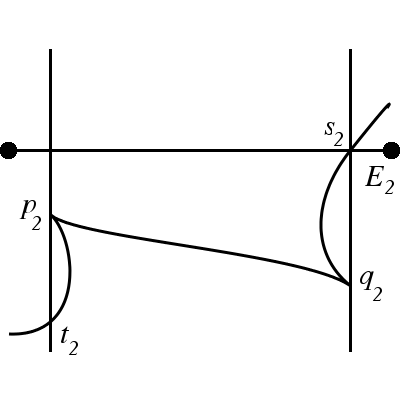}}\;
Blowing down the lines $\bar{T}_1$ and $\bar{L}_1$, we get the surface $X_2$ with horizontal section $E_2$, where $E_2^2=-1$. We let $T_2$ denote the transform of $\bar{E}_{p_1}$, and $L_2$ the transform of $\bar{E}_{q_1}$ on $X_2$. The curve $C_2$ has two cusps $p_2$ and $q_2$, of type $A_2$ and $A_4$ respectively. We have 
\begin{equation*}
\begin{split}
T_2 \cdot C_2 &= 2 \cdot p_2 + 1 \cdot t_2,\\
L_2 \cdot C_2 &= 2 \cdot q_2 + 1 \cdot s_2,\\
E_2 \cdot C_2 &= 1 \cdot s_2.
\end{split}
\end{equation*}
\end{description}

\noi{\bf Blowing down $E_2$}\\
\begin{description}
\item
\parpic[l]{\includegraphics[width=0.3\textwidth]{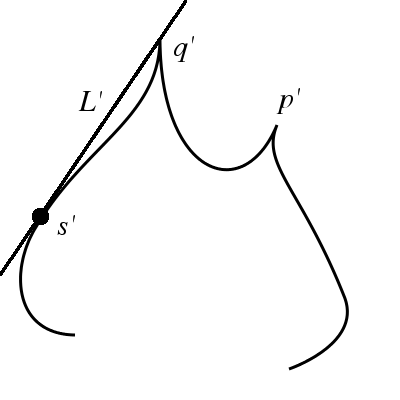}}\;
\noindent The curve $C'$ has two cusps, $p'$ of type $A_2$ and $q'$ of type $A_4$. Since both $E_2$ and $L_2$ intersect $C_2$ transversally in the point $s_2$, this point is transformed into a smooth point $s' \in C'$ with the property that the tangent line at $s'$ will be the transform $L'$ of $L_2$. This is a result of the intersection multiplicity 
\begin{equation*}
\begin{split}
\qquad \qquad \qquad \qquad \qquad (L' \cdot C')_{s'}&=(L_2\cdot C_2)_{s_2} + (E_2 \cdot C_2)_{s_2}\\
&=2.
\end{split}
\end{equation*}
\end{description}

\pagebreak
\noi Note that by the intersection $L_2 \cdot C_2$ above, $L'$ must also intersect $C'$ at the $A_4$-cusp $q'$. We have $$L' \cdot C' = 2 \cdot s' + 2 \cdot q'.$$

\vspace{2mm}
\noi We observe that $C'$ is a quartic curve, 
\begin{equation*}
\begin{split}
d'&=2 \cdot d - m_{p}-m_{p_1}-m_{q_1} \\
&= 2 \cdot 5-2-2-2\\
&=4.
\end{split}
\end{equation*}
 
\noindent As in the previous case, the quartic curve with this cuspidal configuration does not have a point like $s'$. Hence, the quintic curve $C$ does not exist.

\subsubsection{\bf Conclusion} We have three possible cuspidal configurations for bicuspidal quintics.
\begin{center}
	\setlength{\extrarowheight}{2pt}
\begin{tabular}{cc}
\hline
{\bf Curve} & {\bf Cuspidal configuration}\\
\hline
$C_3$ & $(3,2),(2_2)$ \\
$C_4$ & $(3),(2_3)$\\
$C_5$ & $(2_4),(2_2)$\\
\hline
\end{tabular}
\end{center}

\subsection{Three or more cusps}
Assuming we have three or more cusps $p$, $q$, $r$, $\ldots$, the general restrictions leave us with seven cuspidal configurations, of which four can be excluded.

{\samepage
\subsubsection{Maximal multiplicity 3}
\noindent Assuming that $\mu=3$ and letting $p$ be the cusp where $\mu=m_p=3$, Lemma \ref{rhof} on page \pageref{rhof} gives a maximum of two ramification points in the projection of $C$ from $p$. Thus, $C$ can have up to three cusps. By the assumption that $C$ has at least three cusps, $C$ must have exactly three cusps. By B\'{e}zout's theorem, the only possible configuration of cusps for such a curve is $$[(3),(2_a),(2_b)], \qquad a,b \in \mathbb{N}.$$

\noi Theorem \ref{thm:genus} gives the necessary restriction on $a$ and $b$. In this situation, the formula can be reduced to $6=3+a+b$. Therefore, assuming $a\geq b>0$, we must have $a=2$ and $b=1$. Hence, the only possible cuspidal configuration is $$[(3),(2_2),(2)].$$ This curve will be constructed by Cremona transformations at the end of this chapter.}

\subsubsection{Maximal multiplicity 2}
\noindent If $\mu=2$, then all the cusps must have multiplicity $m=2$. Let $p,q,r \ldots$ denote the cusps of $C$ with multiplicity sequences $(2_p),(2_q),(2_r) \ldots$ respectively, $p,q,r \ldots \in \mathbb{N}$. Although this notation is ambiguous, the nature of the objects we refer to by $p, q, r \ldots$ should be clear from the context. 
 
In this situation, Theorem \ref{thm:genus} leads to the sum $p+q+r+ \ldots=6$. By reordering the cusps, we may assume that $p \geq q \geq r \geq \ldots$. Projecting $C$ from the cusp $p$ gives, by Lemma \ref{rhof}, four ramification points. Hence, $C$ has at most five cusps, and we have the following cuspidal configurations.

\begin{description}
\item{\qquad \bf 1)}\;{\bf $C$ has 3 cusps, $p=(2_2),\; q=(2_2),\; r=(2_2).$}
\item{\qquad 2)}\;$C$ has 3 cusps, $p=(2_4),\; q=(2),\; r=(2).$
\item{\qquad 3)}\;$C$ has 3 cusps, $p=(2_3),\; q=(2_2),\; r=(2).$
\item{\qquad \bf 4)}\;{\bf $C$ has 4 cusps, $p=(2_3),\; q=(2),\; r=(2),\; s=(2).$}
\item{\qquad 5)}\;$C$ has 4 cusps, $p=(2_2),\; q=(2_2),\; r=(2),\; s=(2).$
\item{\qquad 6)}\;$C$ has 5 cusps, $p=(2_2),\; q=(2),\; r=(2),\; s=(2),\; t=(2).$
\end{description}

\noindent Only cases 1) and 4) exist. These curves will later be constructed and described. We will now exclude curves 2), 3), 5) and 6) from the list of possible cuspidal configurations.

\subsubsection{\bf Curve 6) -- $[(2_2),(2),(2),(2),(2)]$} 
The curve with 5 cusps is easily excluded from the above list because it contradicts Lemma \ref{rhof}. 
$$\sum_{j=q,r,s,t}(m_{j}-1)+(m_{p.1}-1) \leq 2(d-m_p-1)$$
is contradicted by
$$5=4 \cdot (2-1)+(2-1)\nleq2 \cdot (5-2-1)=4.$$

\subsubsection{\bf Curve 5) -- $[(2_2),(2_2),(2),(2)]$} 
An argument involving the dual curve contradicts the existence of this curve. Theorem \ref{plucker} implies that the dual curve of $C$ has degree 
\begin{equation*}
\begin{split}
d^*=&\;2d-2-\sum_{j=p,q,r,s}(m_j-1)\\
=&\;2\cdot 5-2-4\\
=&\;4. 
\end{split}
\end{equation*}\medskip

\begin{rem} Flenner and Zaidenberg claim that the dual curve has degree $d^*=3$. \end{rem}\medskip

\noindent The degree of the dual curve does not give a contradiction in itself. A contradiction additionally requires a local calculation of the corresponding configuration of singularities on the dual curve. 

We first analyze the simple cusps $r$ and $s$. By Lemma \ref{intmult}, $(C \cdot T_r)_r=3$. We then have a Puiseux parametrization for the germ $(C,r)$, $$(C,r)=(t^2: c_3t^3+ \ldots:1).$$ The dual germ $(C^*,r^*)$ consequently has Puiseux parametrization $$(C^*,r^*)=(a^*t+ \ldots:1:c^*_3t^3+ \ldots).$$ Hence, $r^*$ is an inflection point of type 1 on $C^*$. The same is true for $s^*$.

The fate of the $A_4$-cusps $p$ and $q$, with multiplicity sequences $(2_2)$, can be determined in the same way. Since $C$ is a quintic, the upper and lower bound on $(C \cdot T_p)_p$ are given by identity (\ref{upperandlower}) on page \pageref{upperandlower}, $$4 \leq (C \cdot T_p)_p \leq 5.$$ Thus, we have two possibilities for the intersection multiplicity and the Puiseux parametrization.

\begin{itemize}
\item[--] If $(C \cdot T_p)_p=5$, then the cusp $p$ has local Puiseux parametrization $$(C,p)=(t^2:c_5t^5+\ldots:1).$$ The dual point on the dual germ has Puiseux parametrization $$(C^*,p^*)=(a^*t^3+\ldots:1:c_5^*t^5+\ldots).$$ This implies that $p^*$ is a cusp with multiplicity sequence $(3,2)$.

\item[--] If $(C \cdot T_p)_p=4$, then the cusp $p$ has local Puiseux parametrization $$(C,p)=(t^2:c_4t^4+c_5t^5+\ldots:1).$$ The dual point on the dual germ has Puiseux parametrization $$(C^*,p^*)=(a^*t^2+\ldots:1:c_4^*t^4+c_5^*t^5+\ldots).$$ This implies that $p^*$ is a cusp of the same kind as $p$. They are both $A_4$-cusps with multiplicity sequence $(2_2)$.
\end{itemize}

\noindent Hence, the two cusps $p^*$ and $q^*$ on the dual curve $C^*$, corresponding to $p$ and $q$ on $C$, may be the following pairs.
\begin{description}
\item{\quad i)} Both $p^*$ and $q^*$ have multiplicity sequence $(2_2)$.
\item{\quad ii)} Both $p^*$ and $q^*$ have multiplicity sequence $(3,2)$.
\item{\quad iii)} One cusp, $p^*$, has multiplicity sequence $(2_2)$, and the other cusp, $q^*$, has multiplicity sequence $(3,2)$.\\
\end{description}

\noindent Case i) implies that the curve $C^*$ has genus $g \leq -1$, so it is not irreducible. Hence, $C$ is not a rational cuspidal quintic curve.\\

\noi To exclude case ii), observe that the line $L^*=L^*_{p^*q^*}$ between the two cusps $p^*$ and $q^*$ intersects the curve with multiplicity at least $$\sum_{j=p^*,q^*}(C^* \cdot L^*)_{j} \geq m^*_{p} + m^*_{q}=6.$$ But $6 > d^*=4$, which contradicts B\'{e}zout's theorem. Hence, $C^*$ and $C$ can not exist. A nearly identical argument rules out case iii). \\

\noi We may conclude that no rational cuspidal quintic curve with cuspidal configuration $[(2_2),(2_2),(2),(2)]$ can exist.

\subsubsection{\bf Curve 3) -- $[(2_3),(2_2),(2)]$}

\begin{description}
\item
\parpic[l]{\includegraphics[width=0.3\textwidth]{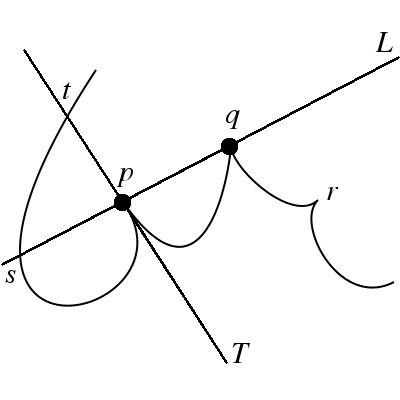}}\;

Assume that the curve $C$ with three cusps, $p$, $q$ and $r$, exists. The cusps have multiplicity sequences $(2_3)$, $(2_2)$ and $(2)$ respectively. Let $T=T_p$ denote the tangent line at the cusp $p$, and let $L=L^{pq}$ denote the line between the cusps $p$ and $q$. By Lemma \ref{intmult}, we have the intersections 
\begin{equation*}
\begin{split}
\qquad \qquad \qquad \qquad \qquad \qquad C \cdot T =&\; 4 \cdot p + 1 \cdot t,\\
C \cdot L=&\;2\cdot p+2 \cdot q+1 \cdot s.
\end{split}
\end{equation*}


\end{description}

\noi We proceed with the Cremona transformation $\psi_2(q,p,T)$. This will give a contradiction. The $A_2$-cusp $r$ will not be altered by this process, and it will therefore only be mentioned at the end of the discussion.\\

\pagebreak
\noi {\bf Blowing up at $p$}

\begin{description}
\item
\parpic[l]{\includegraphics[width=0.3\textwidth]{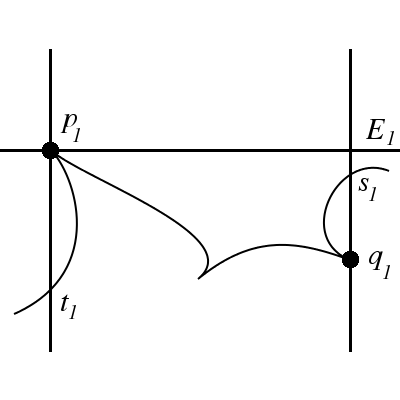}}\;
Blowing up at the cusp $p$ with multiplicity sequence $(2_3)$ gives a ruled surface $X_1$ with the exceptional line $E_1$. The transform $p_1$ of $p$ is a cusp with multiplicity sequence $(2_2)$. We have the intersection multiplicity and intersection
\begin{equation*}
\begin{split}
(C_1 \cdot T_1)_{p_1}&=(C \cdot T)_p-m_{p}\\
&=4-2\\
&=2,\\
E_1 \cdot C_1 &= m_{p} \cdot p_1\\
&=2 \cdot p_1.
\end{split}
\end{equation*}
Furthermore, we also have the intersection $$C_1 \cdot T_1=2 \cdot p_1+ 1 \cdot t_1.$$ 
For the $A_4$-cusp $q$ the situation is unaltered for the transform of the point, $q_1$, but the blowing up process gives the intersection $$C_1 \cdot L_1=2 \cdot q_1 + 1 \cdot s_1.$$
\end{description}

\noi {\bf Elementary transformations in $p_1$ and $q_1$}
\begin{description}
\item
\parpic[l]{\includegraphics[width=0.3\textwidth]{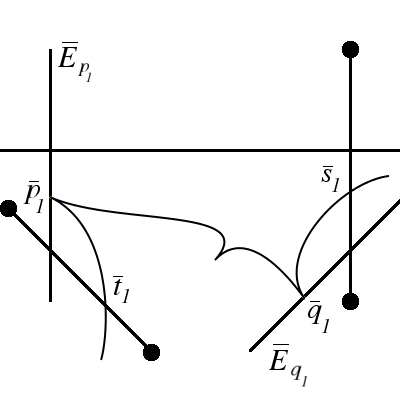}}\;
Blowing up at $p_1$ and $q_1$, we get on the surface $\bar{X}_1$ the exceptional lines $\bar{E}_{p_1}$ and $\bar{E}_{q_1}$. We get two $A_2$-cusps $\bar{p}_1$ and $\bar{q}_1$. Additionally, $\bar{E}_{p_1}$ separates $T_1$ and $E_1$, while $\bar{E}_{q_1}$ does not separate $L_1$ and $E_1$. Besides, we have the intersections 
\begin{equation*}
\begin{split}
\bar{E}_1 \cap \bar{C}_1&= \emptyset, \\
\bar{E}_{p_1} \cdot \bar{C}_1 &= 2 \cdot \bar{p}_1,\\
\bar{T}_1 \cdot \bar{C}_1&= 1 \cdot \bar{t}_1,\\
\bar{E}_{q_1} \cdot \bar{C}_1&=2 \cdot \bar{q}_1,\\
\bar{L}_1 \cdot \bar{C}_1 &=1 \cdot \bar{s}_1.
\end{split}
\end{equation*}
\end{description}

\begin{description}
\item
\parpic[l]{\includegraphics[width=0.3\textwidth]{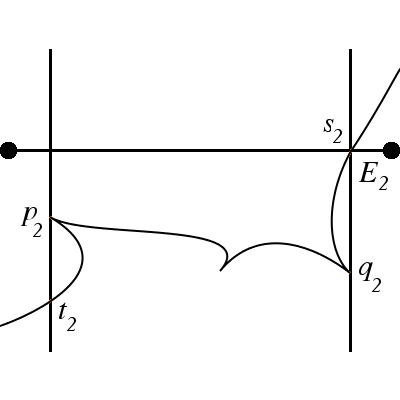}}\;
Blowing down the lines $\bar{T}_1$ and $\bar{L}_1$, we get a surface $X_2$ with horizontal section $E_2$. We let $T_2$ denote the transform of $\bar{E}_{p_1}$, and $L_2$ the transform of $\bar{E}_{q_1}$ on $X_2$. $C_2$ has three $A_2$-cusps, $p_2$, $q_2$ and $r_2$. We have
\begin{equation*}
\begin{split}
T_2 \cdot C_2 &= 2 \cdot p_2 + 1 \cdot t_2,\\
L_2 \cdot C_2 &= 2 \cdot q_2 + 1 \cdot s_2,\\
E_2 \cdot C_2 &= 1 \cdot s_2.
\end{split}
\end{equation*}
\end{description}

\noi{\bf Blowing down $E_2$}\\
\noindent Since neither of the $A_2$-cusps $p_2$, $q_2$ or $r_2$ lie on $E_2$, blowing down $E_2$ does not alter these. Hence, the curve $C'$ has three $A_2$-cusps $p'$, $q'$ and $r'$. However, since both $E_2$ and $L_2$ intersect $C_2$ transversally in the point $s_2$, this point is transformed into a point $s' \in C'$, with the property that the tangent line at $s'$ will be the transform $L'$ of $L_2$. This happens because  
\begin{equation*}
\begin{split}
(L' \cdot C')_{s'}&=(L_2\cdot C_2)_{s_2} + (E_2 \cdot C_2)_{s_2}\\
&=2.
\end{split}
\end{equation*}

\noi By the above intersection $L_2 \cdot C_2=2 \cdot q_2 + 1 \cdot s_2$, we have $$L' \cdot C' = 2 \cdot s' + 2 \cdot q'.$$

\noi Before we go on with the exclusion, we observe that $C'$ is a quartic curve, 
\begin{equation*}
\begin{split}
d'&=2 \cdot d - m_{p}-m_{p_1}-m_{q_1} \\
&= 2 \cdot 5-2-2-2\\
&=4.
\end{split}
\end{equation*}
 
\noindent Up to projective equivalence, there exists only one rational cuspidal quartic with three cusps, the tricuspidal quartic. By investigating the defining polynomial of this explicitly given curve directly, it can easily be proved that a smooth point like $s'$ does not exist. 

Let $F=x^2y^2+y^2z^2+x^2z^2-2xyz(x+y+z)$ be the defining polynomial of the tricuspidal quartic $D$ with cusps $p$, $q$ and $r$. By symmetry, these cusps must have similar properties. Hence, it is enough to investigate one of them. Let $q$ be the cusp in $(1:0:0)$. The polar of $D$ at $q$ is given by $$P_qD=\mathcal{V}(2xy^2+2xz^2-4xyz-2y^2z-2yz^2).$$ A calculation in {\em Maple} reveals that the only intersection points of the curves $P_qD$ and $D$ are the cusps. A point $s$ with tangent $T_s$ that intersects $D$ in any of the cusps, does therefore not exist on this curve. Hence, neither the curve $C'$ nor the quintic curve $C$ with cuspidal configuration $[(2_3),(2_2),(2)]$ exist.


 \pagebreak
\subsubsection{\bf Curve 2)  -- $[(2_4),(2),(2)]$}

\begin{description}
\item
\parpic[l]{\includegraphics[width=0.3\textwidth]{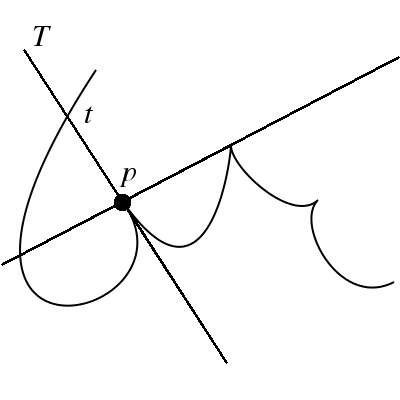}}\;
Assume that $C$ is a quintic with three cusps, $p$, $q$ and $r$, with multiplicity sequences $(2_4)$, $(2)$ and $(2)$ respectively. Let $T=T_p$ denote the tangent line at the $A_8$-cusp $p$. Since $C$ is a quintic and because of Lemma \ref{FlZa14}, we have $$C \cdot T = 4 \cdot p + 1 \cdot t,$$ 

for a smooth point $t \in C$.
\end{description}

\noindent The application of the Cremona transformation $\psi(p,T,-)$ will give a quartic curve similar to the one described above, and will therefore rule out this quintic curve. 

The Cremona transformation will not affect the two $A_2$-cusps $q$ and $r$, and they will therefore only be mentioned at the end of the discussion.\\

\begin{samepage}
\noi {\bf Blowing up at $p$}
\begin{description}
\item \parpic[l]{\includegraphics[width=0.3\textwidth]{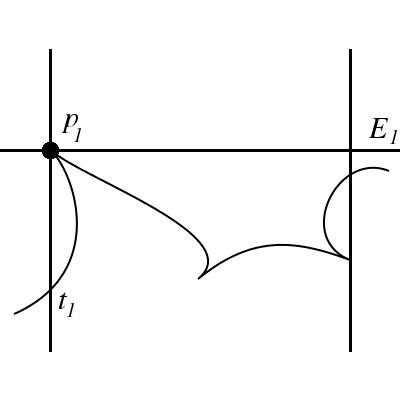}}\;
Blowing up at $p$, we get the ruled surface $X_1$, where the transform $p_1$ of $p$ is a cusp with multiplicity sequence $(2_3)$. By elementary properties of the blowing up process, we have the intersections 
\begin{equation*}
\begin{split}
C_1 \cdot E_1 &= m_{p} \cdot p_1\\
&= 2 \cdot p_1,\\
C_1 \cdot T_1 &= 2 \cdot p_1 + 1 \cdot t_1.
\end{split}
\end{equation*} 
\end{description}
\end{samepage}

\begin{samepage}
\noi {\bf Elementary transformation in $p_1$}
\begin{description}
\item \parpic[l]{\includegraphics[width=0.3\textwidth]{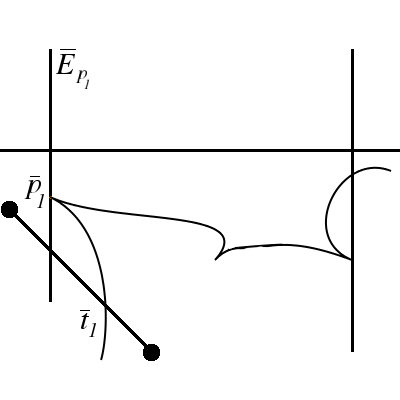}}\;
Blowing up at $p_1$ gives the cusp $\bar{p}_1$ with multiplicity sequence $(2_2)$ on the surface $\bar{X}_1$. The exceptional line $\bar{E}_{p_1}$ separates $E_1$ and $T_1$. Additionally, we get the intersections 
\begin{equation*}
\begin{split}
\bar{E}_{p_1} \cdot \bar{C}_1 &= m_{p_1} \cdot \bar{p}_1 \\
&= 2 \cdot \bar{p}_1,\\
\bar{T}_1 \cdot \bar{C}_1 &= 1 \cdot \bar{t}_1.
\end{split}
\end{equation*}
\end{description}
\end{samepage}

\begin{description}
\item \parpic[l]{\includegraphics[width=0.3\textwidth]{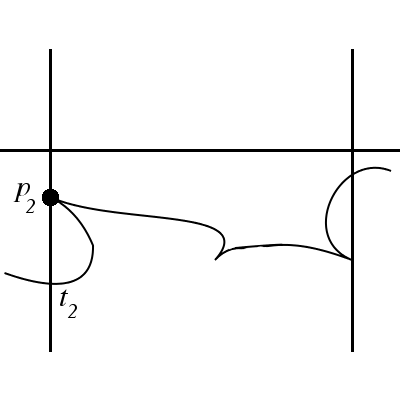}}\;
Blowing down $\bar{T}_1$ gives the surface $X_2$. The transform $T_2$ of $\bar{E}_{p_1}$ intersects $C_2$ in the $A_4$-cusp $p_2$ and a smooth point $t_2$.  In particular, since $E_1$ and $T_1$ were separated by $\bar{E}_{p_1}$, $t_2 \notin E_2$. We have the intersection $$T_2 \cdot C_2 = 2 \cdot p_2 + 1\cdot t_2.$$
\end{description}

\noi{\bf Elementary transformation in $p_2$}
\begin{description}
\item \parpic[l]{\includegraphics[width=0.3\textwidth]{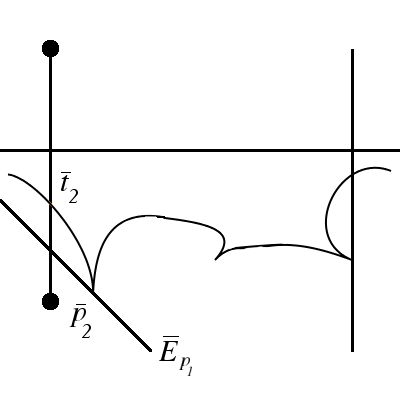}}\;
Blowing up at $p_2$ gives the $A_2$-cusp $\bar{p}_2$ with multiplicity sequence $(2)$ on the surface $\bar{X}_2$. Note that the exceptional line $\bar{E}_{p_2}$ does not separate $E_2$ and $T_2$. We have the intersections 
\begin{equation*}
\begin{split}
\bar{E}_{p_2} \cdot \bar{C}_2 &= 2 \cdot \bar{p}_2,\\
\bar{T}_2 \cdot \bar{C}_2 &= 1 \cdot \bar{t}_2.
\end{split}
\end{equation*}

\item \parpic[l]{\includegraphics[width=0.3\textwidth]{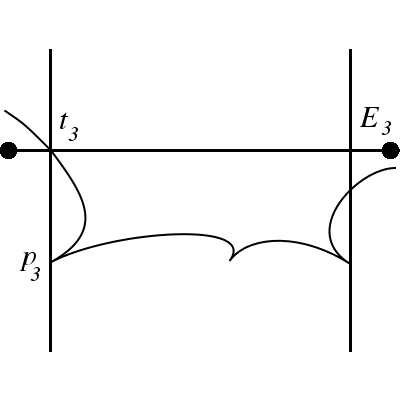}}\;
Blowing down $\bar{T}_2$ gives the surface $X_3$ with the $A_2$-cusp $p_3$ and the intersections 
\begin{equation*}
\begin{split}
T_3 \cdot C_3 &= 2 \cdot p_3+1 \cdot t_3,\\
E_3 \cdot C_3 &= 1 \cdot t_3. 
\end{split}
\end{equation*}
\\
\end{description}

\noi {\bf Blowing down $E_3$}\\
\noi Since neither of the cusps $p_3$, $q_3$ or $r_3$ lie on $E_3$, blowing down $E_3$ gives a curve $C'$ with three $A_2$-cusps $p'$, $q'$ and $r'$. However, because both $E_3$ and $T_3$ intersect $C_3$ transversally at $t_3$, the curve $C'$ also has a point $t'$ such that the tangent $T'_{t'}$ intersects $C'$ at the cusp $p'$.\\ 

\noi Observe that $C'$ is a quartic since 
\begin{equation*}
\begin{split}
d'&=2 \cdot d - m_p-m_{p_1}-m_{p_2}\\
&=2 \cdot 5 - 2-2-2 \\
&=4.
\end{split}
\end{equation*}

\noindent Apart from some notation, the exact same argument is used in the exclusion of this curve as for the exclusion of curve number 3). Hence, curve number 2) can also be excluded from the list of possible rational cuspidal quintics.

\subsubsection{Conclusion} We have three possible cuspidal configurations with three or more cusps.
\begin{center}
	\setlength{\extrarowheight}{2pt}
\begin{tabular}{cc}
\hline
{\bf Curve} & {\bf Cuspidal configuration}\\
\hline
$C_6$ & $(3),(2_2),(2)$ \\
$C_7$ & $(2_2),(2_2),(2_2)$\\
$C_8$ & $(2_3),(2),(2),(2)$\\
\hline
\end{tabular}
\end{center}
\vspace{15mm}
\section{Possible cuspidal configurations}
We now have eight possible cuspidal configurations for a rational cuspidal quintic curve. There exist curves with all these cuspidal configurations. Up to projective equivalence, there are a few more curves. We will not go into details concerning this classification here. 

Up to projective equivalence, all rational cuspidal quintic curves were described and found by Namba in \cite[Thm. 2.3.10., pp.179--182]{Namba}. An overview of the curves is presented in Table \ref{tab:degree5}.

\begin{table}[H]
  \renewcommand\thesubtable{}
  \setlength{\extrarowheight}{2pt}
\centering
	{\begin{tabular}{ccll}
	\hline
	{\bf \# Cusps}&	{\bf Curve} &{\bf Cuspidal configuration} & {\bf \# Curves}\\
	\hline 
	\multirow{2}{14mm}{1}&$C_1$&$\qquad(4)$&3 -- ABC\\
	&$C_2$&$\qquad (2_6)$&1\\
	\hline
	\multirow{3}{14mm}{2}&$C_3$&$\qquad (3,2),(2_2)$&2 -- AB\\
	&$C_4$&$\qquad (3), (2_3)$&1\\
	&$C_5$&$\qquad (2_4),(2_2)$&1\\
		\hline
		\multirow{2}{14mm}{3}&$C_6$&$\qquad (3),(2_2),(2)$&1\\
	&$C_7$&$\qquad (2_2),(2_2),(2_2)$&1\\
	\hline
	\multirow{1}{14mm}{4}	&$C_8$&$\qquad (2_3),(2),(2),(2)$&1\\
	\hline
	\end{tabular}}
	\caption {Rational cuspidal quintic curves.}
	\label{tab:degree5}
	\end{table}

\pagebreak
\section{Rational cuspidal quintics} 
We will now briefly describe how the rational cuspidal quintics can be constructed. Additionally, we will list some of the most important properties of each curve.

\subsubsection{Curve $C_{1A}$ -- $[(4)]$}
A cuspidal quintic with one cusp with multiplicity sequence $(4)$ and one inflection point of type 3 can be constructed using a Cremona transformation with two proper base points.\\ 

\parpic[l]{\includegraphics[width=0.3\textwidth]{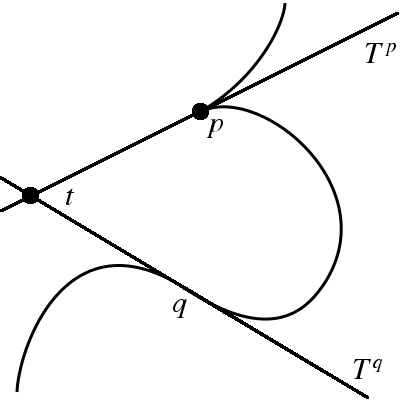}}
\picskip{10}
\noi Let $C$ be the cuspidal quartic with cusp $p$, $\ol{m}_p=3$, and one inflection point, $q$, of type 2. The tangent $T^p$ intersects $C$ at $p$, and the tangent $T^q$ intersects $C$ at $q$. Both lines intersect $C$ with intersection multiplicity $4$ in the respective points. Hence, the tangents do not intersect $C$ at any other points. Denote by $t$ the intersection point of $T^p$ and $T^q$.\\

\noi The Cremona transformation $\psi_2(p,t,T^q)$ transforms this curve into a unicuspidal quintic with cuspidal configuration $[(4)]$ and one inflection point of type 3.\\

{\samepage
\noi The curve $C_{1A}$ is given by the below parametrization, and it has the following properties.
\begin{center}
	\setlength{\extrarowheight}{2pt}
	{\begin{tabular}{ccc}
	\multicolumn{3}{c}{$(s^5:st^4:t^5)$}\\
			\multicolumn{3}{c}{} \\ 
	\hline
	\multicolumn{3}{c}{{\bf \# Cusps = 	1 }}\\
	{\bf Cusp $p_j$} & {\bf $(C \cdot T_{p_j})_{p_j}$}& {\bf $(C \cdot H_C)_{p_j}$ }\\
	\hline 
	$(4)$& 5 & 42\\
		\hline
		&&\\
		\hline
			\multicolumn{3}{c}{{\bf \# Inflection points =	1 }}\\
{\bf Inflection point $q_j$} & {\bf $(C \cdot T_{q_j})_{q_j}$}& {\bf $(C \cdot H_C)_{q_j}$ }\\	
	\hline
	$q_1$ & 5 & 3 \\
	\hline
	\end{tabular}}
	\end{center}}

\newpage
\subsubsection{Curve $C_{1B}$ -- $[(4)]$}
A cuspidal quintic with one cusp with multiplicity sequence $(4)$ and two inflection points of type 1 and 2 respectively, can be constructed using a Cremona transformation with two proper base points.\\ 

\parpic[l]{\includegraphics[width=0.3\textwidth]{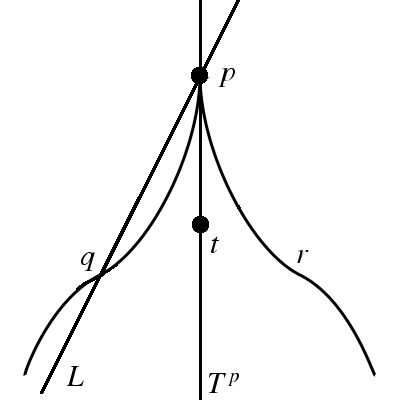}}
\picskip{10}
\noi Let $C$ be the cuspidal quartic with a cusp $p$, $\ol{m}_p=3$, and two inflection points, $q$ and $r$, of type 1. Let $T^p$ be the tangent to $C$ at $p$. Choose a point $t$ on $T^p$. Let $L=L^{pq}$ denote the line between $p$ and the inflection point $q$. By B\'{e}zout's theorem, $L$ and $C$ do not intersect in any other points. We have the intersections 
\begin{equation*}
\begin{split}
T^p \cdot C &= 4 \cdot p,\\
L \cdot C &=3 \cdot p + 1 \cdot q.
\end{split}
\end{equation*} 

\noi The Cremona transformation $\psi_2(t,p,L)$ transforms this curve into a unicuspidal quintic with cuspidal configuration $[(4)]$ and two inflection points of type 1 and 2.\\

{\samepage
\noi The curve $C_{1B}$ is given by the below parametrization, and it has the following properties.
\begin{center}
	\setlength{\extrarowheight}{2pt}
	{\begin{tabular}{ccc}
		\multicolumn{3}{c}{$(s^5-s^4t:st^4:t^5)$} \\ 
		\multicolumn{3}{c}{} \\ 
	\hline
	\multicolumn{3}{c}{{\bf \# Cusps = 	1 }}\\
	{\bf Cusp $p_j$} & {\bf $(C \cdot T_{p_j})_{p_j}$}& {\bf $(C \cdot H_C)_{p_j}$ }\\
	\hline 
	$(4)$& 5 & 42\\
	\hline
		&&\\
		\hline
			\multicolumn{3}{c}{{\bf \# Inflection points =	2 }}\\
{\bf Inflection point $q_j$} & {\bf $(C \cdot T_{q_j})_{q_j}$}& {\bf $(C \cdot H_C)_{q_j}$ }\\	
	\hline
	$q_1$ & 4 & 2 \\
	$q_2$ & 3 & 1 \\
	\hline
	\end{tabular}}
	\end{center}}

\newpage
\subsubsection{Curve $C_{1C}$ -- $[(4)]$}
A cuspidal quintic with one cusp with multiplicity sequence $(4)$ and three inflection points of type 1, can be constructed using a Cremona transformation with two base points.\\ 

\parpic[l]{\includegraphics[width=0.3\textwidth]{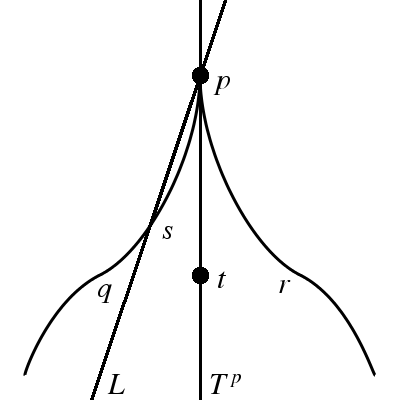}}
\picskip{10}
\noi Let $C$ be the cuspidal quartic with a cusp $p$, $\ol{m}_p=3$ and two inflection points $q$ and $r$ of type 1. Let $T^p$ denote the tangent line at $p$. Choose a point $t$ on $T^p$. Furthermore, choose a point $s \in C$ different from the above mentioned points. Let $L=L^{ps}$ denote the line between the cusp $p$ and the point $s$. We have the intersections 
\begin{equation*}
\begin{split}
T^p \cdot C &= 4 \cdot p,\\ 
L \cdot C &=3 \cdot p + 1 \cdot s. 
\end{split}
\end{equation*}

\noi The Cremona transformation $\psi_2(t,p,L)$ transforms this quartic curve into a unicuspidal quintic with cuspidal configuration $[(4)]$ and three inflection points of type 1.\\

\noi The curve $C_{1C}$ is given by the below parametrization, and it has the following properties.
\begin{center}
	\setlength{\extrarowheight}{2pt}
	{\begin{tabular}{ccc}
		\multicolumn{3}{c}{$(s^5+as^4t-(1+a)s^2t^2:st^4:t^5)$, \qquad $a\in \mathbb{C},\; a \neq -1.$} \\ 
		\multicolumn{3}{c}{} \\ 
	\hline
	\multicolumn{3}{c}{{\bf \# Cusps = 	1 }}\\
	{\bf Cusp $p_j$} & {\bf $(C \cdot T_{p_j})_{p_j}$}& {\bf $(C \cdot H_C)_{p_j}$ }\\
	\hline 
	$(4)$& 5 & 42\\
	\hline
		&&\\
		\hline
			\multicolumn{3}{c}{{\bf \# Inflection points =	1 }}\\
{\bf Inflection point $q_j$} & {\bf $(C \cdot T_{q_j})_{q_j}$}& {\bf $(C \cdot H_C)_{q_j}$ }\\	
	\hline
	$q_1$ & 3 & 1 \\
	$q_2$ & 3 & 1 \\
	$q_3$ & 3 & 1 \\
	\hline
	\end{tabular}}
	\end{center}

\newpage
\subsubsection{Curve $C_2$ -- $[(2_6)]$}
The rational cuspidal quintic with cuspidal configuration $[(2_6)]$ can be constructed by transforming the unicuspidal ramphoid quartic using a Cremona transformation with one base point.

Since the action of a Cremona transformation with one base point is hard to analyze in $\mathbb{P}^2$, we will construct this curve explicitly with a transformation that is known to work. Let $C$ be given by $$C=\V((yx-z^2)^2-x^3z).$$ Then the Cremona transformation $\psi_1$ transforms this curve into $$C_2=\V(xy^4-2x^2y^2z+x^3z^2-2y^3z^2+2xyz^3+z^5).$$\medskip

\parpic[l]{\includegraphics[width=0.3\textwidth]{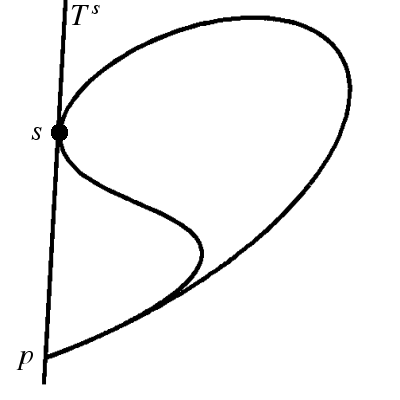}}
\picskip{10}
\noi The implicit construction can be explained adequately. Let $C$ be the quartic with one cusp $p$ with multiplicity sequence $(2_3)$. The polar of $C$ at the cusp $p$ intersects $C$ in a smooth point $s$. Denote by $T^s$ the tangent line at $s$, which intersects $C$ in two points, \begin{equation*}C \cdot T = 2 \cdot p+2 \cdot s.\end{equation*}
The appropriate Cremona transformation in this situation consists of blowing up $s$, performing elementary transformations in the two successive infinitely near points $s_1$ and $s_2$ on the strict transforms of $C$, and then blowing down the horizontal section.\\

\noi The curve $C_{2}$ is given by the below parametrization, and it has the following properties.

\begin{center}
	\setlength{\extrarowheight}{2pt}
	{\begin{tabular}{ccc}
		\multicolumn{3}{c}{$(s^4t:s^2t^3-s^5:t^5-2s^3t^2)$} \\ 
		\multicolumn{3}{c}{} \\ 
	\hline
	\multicolumn{3}{c}{{\bf \# Cusps = 	1 }}\\
	{\bf Cusp $p_j$} & {\bf $(C \cdot T_{p_j})_{p_j}$}& {\bf $(C \cdot H_C)_{p_j}$ }\\
	\hline 
	$(2_6)$& 4 & 39\\
	\hline
		&&\\
		\hline
			\multicolumn{3}{c}{{\bf \# Inflection points =	6 }}\\
{\bf Inflection point $q_j$} & {\bf $(C \cdot T_{q_j})_{q_j}$}& {\bf $(C \cdot H_C)_{q_j}$ }\\	
	\hline
	$q_j, \; j=1,\ldots,6$ & 3 & 1 \\
	\hline
	\end{tabular}}
	\end{center}

\newpage
\subsubsection{Curve $C_{3A}$ -- $[(3,2),(2_2)]$}
A cuspidal quintic with cuspidal configuration $[(3,2),(2_2)]$ and no inflection points can be constructed using a Cremona transformation with two base points.\\

\parpic[l]{\includegraphics[width=0.3\textwidth]{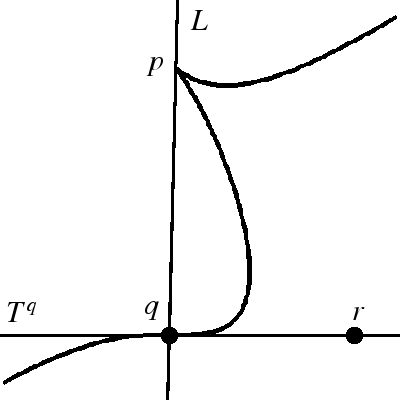}}
\picskip{10}
\noi Let $C$ be the cuspidal cubic with a cusp $p$, $\ol{m}_p=(2)$, and one inflection point, $q$, of type 1. Let $T^q$ be the tangent line at $q$. Let $L=L^{pq}$ be the line between $p$ and $q$, which by B\'{e}zout's theorem does not intersect $C$ in any other point. Choose an arbitrary point $r \in T^q$. We have the intersections \begin{equation*} \begin{split} T^q \cdot C &= 3 \cdot q,\\ L \cdot C &=2 \cdot p + 1 \cdot q. \end{split} \end{equation*}\medskip

\noi The Cremona transformation $\psi_2(r,q,L)$ transforms this cubic into the bicuspidal quintic with cuspidal configuration $[(3,2),(2_2)]$.\\

\noi The curve $C_{3A}$ is given by the below parametrization, and it has the following properties.
\begin{center}
	\setlength{\extrarowheight}{2pt}
	{\begin{tabular}{ccc}
	\multicolumn{3}{c}{$(s^5:s^3t^2:t^5)$}\\
			\multicolumn{3}{c}{} \\ 
	\hline
	\multicolumn{3}{c}{{\bf \# Cusps = 	2 }}\\
	{\bf Cusp $p_j$} & {\bf $(C \cdot T_{p_j})_{p_j}$}& {\bf $(C \cdot H_C)_{p_j}$ }\\
	\hline 
	$(3,2)$& 5 & 29\\
	$(2_2)$& 5 & 16\\
		\hline
		&&\\
		\hline
			\multicolumn{3}{c}{{\bf \# Inflection points =	0 }}\\
	\hline
	\end{tabular}}
	\end{center}

\newpage
\subsubsection{Curve $C_{3B}$ -- $[(3,2),(2_2)]$}
A cuspidal quintic with cuspidal configuration $[(3,2),(2_2)]$ and one inflection point of type 1 can be constructed using a Cremona transformation with two base points. \\

\parpic[l]{\includegraphics[width=0.3\textwidth]{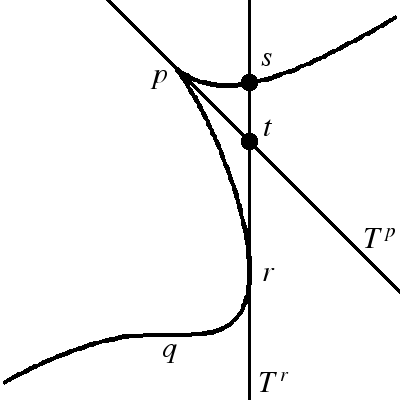}}
\picskip{10}
\noi Let $C$ be the cuspidal cubic with a cusp $p$, $\ol{m}_p=(2)$, and one inflection point, $q$, of type 1. Choose an arbitrary smooth point $r \in C,\; r\neq q$ and let $T^r$ denote the tangent line at this point. By B\'{e}zout's theorem, $$T^r \cdot C=2\cdot r+1\cdot s,$$ for another smooth point $s$. Make sure $s \neq q$. Furthermore, denote by $T^p$ the tangent line at the cusp $p$. The tangents $T^r$ and $T^p$ intersect at a point $t \notin C$.\\

\noi The Cremona transformation $\psi_2(s,t,T^p)$ transforms this cubic into the bicuspidal quintic with cuspidal configuration $[(3,2),(2_2)]$ and one inflection point of type 1.\\

\noi The curve $C_{3B}$ is given by the below parametrization, and it has the following properties.
\begin{center}
	\setlength{\extrarowheight}{2pt}
	{\begin{tabular}{ccc}
	\multicolumn{3}{c}{$(s^5:s^3t^2:st^4+t^5)$}\\
			\multicolumn{3}{c}{} \\ 
	\hline
	\multicolumn{3}{c}{{\bf \# Cusps = 	2 }}\\
	{\bf Cusp $p_j$} & {\bf $(C \cdot T_{p_j})_{p_j}$}& {\bf $(C \cdot H_C)_{p_j}$ }\\
	\hline 
	$(3,2)$& 5 & 29\\
	$(2_2)$& 5 & 15\\
		\hline
		&&\\
		\hline
			\multicolumn{3}{c}{{\bf \# Inflection points =	1 }}\\
			{\bf Inflection point $q_j$} & {\bf $(C \cdot T_{q_j})_{q_j}$}& {\bf $(C \cdot H_C)_{q_j}$ }\\
	\hline
		$q_1$& 3 & 1\\
			\hline
	\end{tabular}}
	\end{center}

\newpage
\subsubsection{Curve $C_4$ -- $[(3), (2_3)]$}
A cuspidal quintic with cuspidal configuration $[(3),(2_3)]$ and two inflection points of type 1 can be constructed by a Cremona transformation with two base points.\\ 

\parpic[l]{\includegraphics[width=0.3\textwidth]{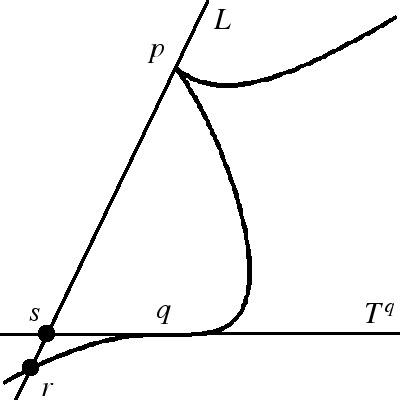}}
\picskip{9}
\noi Let $C$ be the cuspidal cubic with a cusp $p$, $\ol{m}_p=(2)$, and one inflection point, $q$, of type 1. Choose an arbitrary smooth point $r \in C,\; r\neq q$ and let $L=L^{pr}$ denote the line through $p$ and $r$. Furthermore, denote by $T^q$ the tangent line at the inflection point $q$. The lines $T^q$ and $L$ intersect at a point $s \notin C$. We have the intersections \begin{equation*} \begin{split} L \cdot C &=2 \cdot p + 1 \cdot r,\\ T^q \cdot C &=3 \cdot q. \end{split} \end{equation*} \medskip

\noi The Cremona transformation $\psi_2(r,s,T^q)$ transforms the cubic into the bicuspidal quintic with cuspidal configuration $[(3),(2_3)]$ and two inflection points of type 1.\\

\noi The curve $C_{4}$ is given by the below parametrization, and it has the following properties.
\begin{center}
	\setlength{\extrarowheight}{2pt}
	{\begin{tabular}{ccc}
	\multicolumn{3}{c}{$(s^4t-\frac{1}{2}s^5:s^3t^2:\frac{1}{2}st^4+t^5)$}\\
			\multicolumn{3}{c}{} \\ 
	\hline
	\multicolumn{3}{c}{{\bf \# Cusps = 	2 }}\\
	{\bf Cusp $p_j$} & {\bf $(C \cdot T_{p_j})_{p_j}$}& {\bf $(C \cdot H_C)_{p_j}$ }\\
	\hline 
	$(3)$& 4 & 22\\
	$(2_3)$& 4 & 21\\
		\hline
		&&\\
		\hline
			\multicolumn{3}{c}{{\bf \# Inflection points =	2 }}\\
			{\bf Inflection point $q_j$} & {\bf $(C \cdot T_{q_j})_{q_j}$}& {\bf $(C \cdot H_C)_{q_j}$ }\\
	\hline
		$q_1$& 3 & 1\\
		$q_2$& 3 & 1\\
			\hline
	\end{tabular}}
	\end{center}

\newpage
\subsubsection{Curve $C_5$ -- $[(2_4),(2_2)]$}
A cuspidal quintic with cuspidal configuration $[(2_4),(2_2)]$ and three inflection points of type 1 can be constructed by a Cremona transformation with two base points.\\ 

\parpic[l]{\includegraphics[width=0.3\textwidth]{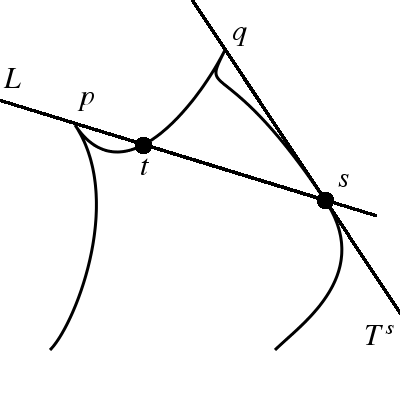}}
\picskip{9}
\noi Let $C$ be the bicuspidal quartic with two cusps, $p$ and $q$, with $\ol{m}_p=(2_2)$ and $\ol{m}_q=(2)$, and one inflection point, $r$, of type 1. The polar of $C$ at the cusp $q$ intersects $C$ in a smooth point $s$. Denote by $T^s$ the tangent line at $s$, which has the property $$T^s \cdot C =2 \cdot q+ 2\cdot s.$$ Let $L=L^{ps}$ be the line through the cusp $p$ and the smooth point $s$. Since $T^p \cdot C=4\cdot p$, we have that $L \neq T^p$. By B\'{e}zout's theorem, $L$ must intersect $C$ in yet another smooth point $t$. We have the intersection $$L \cdot C=2 \cdot p+ 1 \cdot s + 1 \cdot t.$$ \medskip

\noi The Cremona transformation $\psi_2(t,s,T^s)$ transforms the quartic into the bicuspidal quintic with cuspidal configuration $[(2_4),(2_2)]$ and three inflection points of type 1.\\

\noi The curve $C_{5}$ is given by the below parametrization, and it has the following properties.
\begin{center}
	\setlength{\extrarowheight}{2pt}
	{\begin{tabular}{ccc}
	\multicolumn{3}{c}{$(s^4t-s^5:s^2t^3-\frac{5}{32}s^5:-\frac{47}{128}s^5+\frac{11}{16}s^3t^2+st^4+t^5)$}\\
			\multicolumn{3}{c}{} \\ 
	\hline
	\multicolumn{3}{c}{{\bf \# Cusps = 	2 }}\\
	{\bf Cusp $p_j$} & {\bf $(C \cdot T_{p_j})_{p_j}$}& {\bf $(C \cdot H_C)_{p_j}$ }\\
	\hline 
	$(2_4)$& 4 & 27\\
	$(2_2)$& 4 & 15\\
		\hline
		&&\\
		\hline
			\multicolumn{3}{c}{{\bf \# Inflection points =	3 }}\\
			{\bf Inflection point $q_j$} & {\bf $(C \cdot T_{q_j})_{q_j}$}& {\bf $(C \cdot H_C)_{q_j}$ }\\
	\hline
		$q_j, \; j=1,2,3 $& 3 & 1\\
			\hline
	\end{tabular}}
	\end{center}

\newpage
\subsubsection{Curve $C_6$ -- $[(3),(2_2),(2)]$} 
A cuspidal quintic with cuspidal configuration $[(3),(2_2),(2)]$ and no inflection points can be constructed using a Cremona transformation with two proper base points. \\ 

\begin{description}
\item
\parpic[l]{\includegraphics[width=0.3\textwidth]{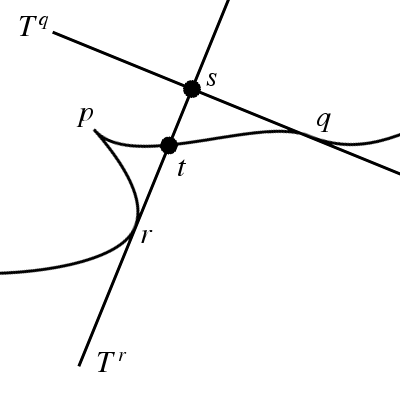}}\;
Let $C$ be the cuspidal cubic with a cusp $p$, where $\ol{m}_p=(2)$, and one inflection point, $q$, of type 1. Denote by $T^q$ the tangent line at $q$. Choose an arbitrary smooth point $r \in C,\; r \neq q$. Then the tangent line at $r$, $T^r$, intersects $T^q$ in a point $s \notin C$. Furthermore, $T^r$ intersects $C$ in another smooth point $t$. We have the intersection $$\qquad \qquad \qquad \qquad \qquad T^r\cdot C=2 \cdot r + 1 \cdot t.$$
\end{description}

\noi The Cremona transformation $\psi_2(t,s,T^q)$ transforms the cubic into the tricuspidal quintic with cuspidal configuration $[(3),(2_2),(2)]$.\\

\noi {\bf Blowing up at $s$}
\begin{description}
\item
\parpic[l]{\includegraphics[width=0.3\textwidth]{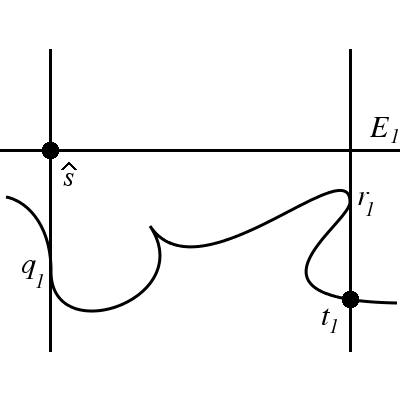}}\;
\noi Blowing up at $s$, we get the ruled surface $X_1$ with horizontal section $E_1$ and the transformed curve $C_1$. We have the intersections
\begin{equation*}
\begin{split}
E_1 \cap C_1&=\emptyset,\\
T^q_1 \cdot C_1&= 3\cdot q_1,\\
T^r_1 \cdot C_1&= 2\cdot r_1 + 1 \cdot t_1.\\
\end{split}
\end{equation*}   
\\
\end{description}

\noi{\bf Elementary transformations in $t_1$ and $\hat{s}$}
\begin{description}
\item 
\parpic[l]{\includegraphics[width=0.3\textwidth]{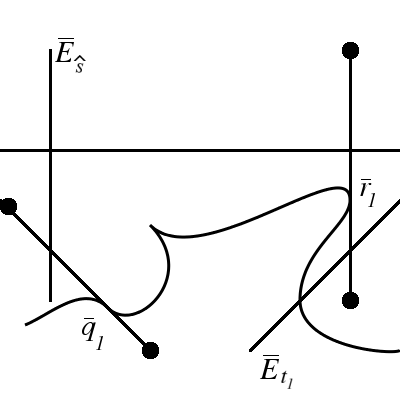}}\;
Blowing up at $t_1$ and $\hat{s}$ gives two exceptional lines $\bar{E}_{t_1}$ and $\bar{E}_{\hat{s}}$. We have the intersections 
\begin{equation*}
\begin{split}
\bar{E}_1 \cap \bar{C}_1&=\emptyset,\\
\bar{E}_{\hat{s}} \cap \bar{C}_1&=\emptyset,\\
\bar{T}^q_1 \cdot \bar{C}_1&=3 \cdot \bar{q}_1,\\
\bar{E}_{t_1} \cdot \bar{C}_1&=1 \cdot \bar{t}_1,\\
\bar{T}^r_1 \cdot \bar{C}_1&=2 \cdot \bar{r}_1.\\
\end{split}
\end{equation*}
\\
\end{description}

\begin{description}
\item 
\parpic[l]{\includegraphics[width=0.3\textwidth]{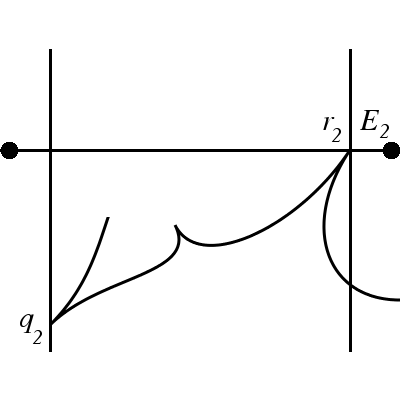}}\;
Blowing down $\bar{T}^q_1$ and $\bar{T}^r_1$ gives the surface $X_2$. On this surface we have $T^q_2$, the strict transform of $\bar{E}_{\hat{s}}$, and $T^r_2$, the strict transform of $\bar{E}_{t_1}$. 
\par Because of the intersection multiplicities \linebreak above, $q_2$ is a cusp with $\ol{m}_{q_2}=(3)$, and $r_2$ is a cusp with $\ol{m}_{r_2}=(2)$. Note that $r_2 \in E_2$. We have the important intersection
\begin{equation*}
\qquad \qquad \qquad \qquad \qquad \qquad E_2 \cdot C_2 = 2 \cdot r_2.
\end{equation*}
\end{description}

\noi{\bf Blowing down $E_2$}\\
\noi Blowing down $E_2$ gives a curve $C'$ with three cusps $p',\,q'$ and $r'$, where $\ol{m}_{p'}=(2)$, $\ol{m}_{q'}=(3)$ and $\ol{m}_{r'}=(2_2)$. \\

\noi To see that $C'$ is a quintic, note that
\begin{equation*}
\begin{split}
d'&=2 \cdot d - m_{s}-m_{\hat{s}}-m_{t_1}\\
&=2 \cdot 4-1-1-1 \\
&=5.
\end{split}
\end{equation*}

\noi The curve $C_{6}$ is given by the below parametrization, and it has the following properties.
\begin{center}
	\setlength{\extrarowheight}{2pt}
	{\begin{tabular}{ccc}
	\multicolumn{3}{c}{$(s^4t-\frac{1}{2}s^5:s^3t^2:-\frac{3}{2}st^4+t^5)$}\\
			\multicolumn{3}{c}{} \\ 
	\hline
	\multicolumn{3}{c}{{\bf \# Cusps = 	3 }}\\
	{\bf Cusp $p_j$} & {\bf $(C \cdot T_{p_j})_{p_j}$}& {\bf $(C \cdot H_C)_{p_j}$ }\\
	\hline 
	$(3)$& 4 & 22\\
	$(2_2)$& 4 & 15\\
	$(2)$& 3 & 8\\
		\hline
		&&\\
		\hline
			\multicolumn{3}{c}{{\bf \# Inflection points =	0 }}\\
	\hline
	\end{tabular}}
	\end{center}

\newpage
\subsubsection{Curve $C_7$ -- $[(2_2),(2_2),(2_2)]$} 
A cuspidal quintic with cuspidal configuration $[(2_2),(2_2),(2_2)]$ and no inflection points can be constructed using a Cremona transformation with two base points.\\

\begin{description}
\item
\parpic[l]{\includegraphics[width=0.3\textwidth]{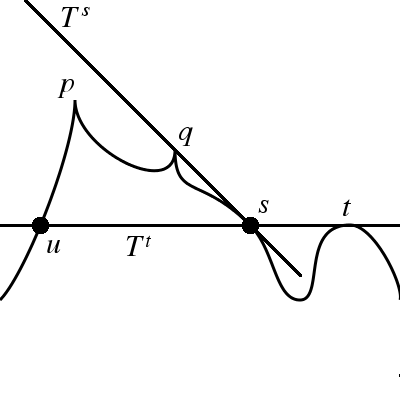}}\;
Let $C$ be the bicuspidal quartic with two cusps $p$ and $q$, with $\ol{m}_p=(2_2)$ and $\ol{m}_q=(2)$, and one inflection point, $r$, of type 1. The polar of $C$ at the point $q$ intersects $C$ in a smooth point $s$. Denote by $T^s$ the tangent line at $s$, which has the property $T^s \cdot C =2 \cdot q+ 2\cdot s.$ The polar of $C$ at the point $s$ intersects $C$ in a smooth point $t, \; t\neq r$. Denote by $T^t$ the tangent line at $t$, which by B\'{e}zout's theorem intersects $C$ in a smooth point $u$. We have the intersection $$T^t \cdot C = 2 \cdot t+ 1 \cdot s+ 1 \cdot u.$$
\end{description}

\noi The Cremona transformation $\psi_2(u,s,T^s)$ transforms the bicuspidal quartic into the tricuspidal quintic with cuspidal configuration $[(2_2),(2_2),(2_2)]$.\\

\noi {\bf Blowing up at $s$}
\begin{description}
\item
\parpic[l]{\includegraphics[width=0.3\textwidth]{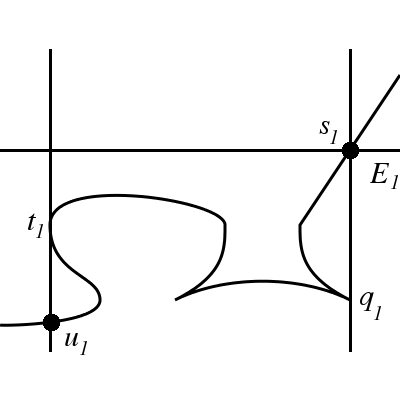}}\;
\noi Blowing up at $s$, we get the ruled surface $X_1$ with horizontal section $E_1$ and the transformed curve $C_1$. We have the intersections
\begin{equation*}
\begin{split}
E_1 \cdot C_1&=1 \cdot s_1,\\
T^t_1 \cdot C_1&= 1\cdot u_1+2\cdot t_1,\\
T^s_1 \cdot C_1&= 1\cdot s_1 + 2 \cdot q_1.
\end{split}
\end{equation*}   
\end{description}

\noi{\bf Elementary transformations in $u_1$ and $s_1$}
\begin{description}
\item 
\parpic[l]{\includegraphics[width=0.3\textwidth]{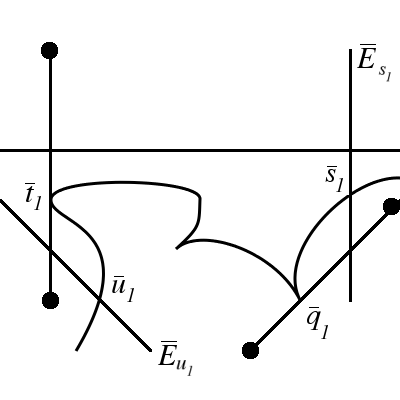}}\;
Blowing up at $u_1$ and $s_1$ gives on $\bar{X}_1$ two exceptional lines $\bar{E}_{u_1}$ and $\bar{E}_{s_1}$. We then have the intersections 
\begin{equation*}
\begin{split}
\bar{E}_1 \cap \bar{C}_1&=\emptyset,\\
\bar{E}_{u_1} \cdot \bar{C}_1&=1 \cdot \bar{u}_1,\\
\bar{T}^t_1 \cdot \bar{C}_1&=2 \cdot \bar{t}_1,\\
\bar{E}_{s_1} \cdot \bar{C}_1&=1 \cdot \bar{s}_1,\\
\bar{T}^s_1 \cdot \bar{C}_1&=2 \cdot \bar{q}_1.\\
\end{split}
\end{equation*}
\\
\end{description}

\begin{description}
\item 
\parpic[l]{\includegraphics[width=0.3\textwidth]{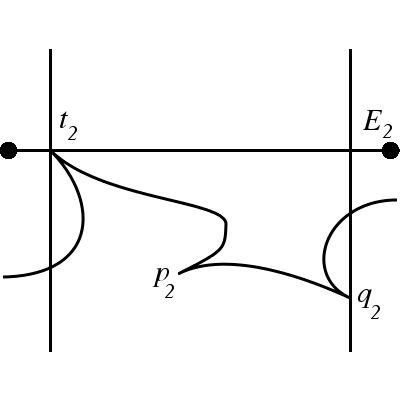}}\;
Blowing down $\bar{T}^t_1$ and $\bar{T}^s_1$ gives the surface $X_2$. On this surface we have $T^t_2$, the strict transform of $\bar{E}_{u_1}$, and $T^s_2$, the strict transform of $\bar{E}_{s_1}$. 
\par Because of the intersection multiplicities \linebreak above, $q_2$ is a cusp with $\ol{m}_{q_2}=(2_2)$, and $t_2$ is a cusp with $\ol{m}_{t_2}=(2)$. Note that $t_2 \in E_2$. \par \noi We have the important intersection
\begin{equation*}
\qquad \qquad \qquad \qquad \qquad E_2 \cdot C_2 = 2 \cdot t_2.
\end{equation*}
\end{description}

\noi{\bf Blowing down $E_2$}\\
\noi Blowing down $E_2$ gives a curve $C'$ with three cusps $p'$, $q'$ and $t'$. All these cusps have multiplicity sequence $(2_2)$.  \\

\noi To see that $C'$ is a quintic, note that
\begin{equation*}
\begin{split}
d'&=2 \cdot d - m_{s}-m_{s_1}-m_{u_1}\\
&=2 \cdot 4-1-1-1 \\
&=5.
\end{split}
\end{equation*}

\noi The curve $C_{7}$ is given by the below parametrization, and it has the following properties.
\begin{center}
	\setlength{\extrarowheight}{2pt}
	{\begin{tabular}{ccc}
	\multicolumn{3}{c}{$(s^4t-s^5:s^2t^3-\frac{5}{32}s^5:-\frac{125}{128}s^5-\frac{25}{16}s^3t^2-5st^4+t^5)$}\\
			\multicolumn{3}{c}{} \\ 
	\hline
	\multicolumn{3}{c}{{\bf \# Cusps = 	3 }}\\
	{\bf Cusp $p_j$} & {\bf $(C \cdot T_{p_j})_{p_j}$}& {\bf $(C \cdot H_C)_{p_j}$ }\\
	\hline 
	$(2_2)$& 4 & 15\\
	$(2_2)$& 4 & 15\\
	$(2_2)$& 4 & 15\\
		\hline
		&&\\
		\hline
			\multicolumn{3}{c}{{\bf \# Inflection points =	0 }}\\
	\hline
	\end{tabular}}
	\end{center}

\newpage
\subsubsection{Curve $C_8$ - $[(2_3),(2),(2),(2)]$} 
The rational cuspidal quintic with four cusps is the dual curve of the unicuspidal ramphoid quartic curve $C$. $C$ has a cusp $p$ of type $A_6$ and three inflection points of type 1. The explicit calculation can be done in the following way.\\

\noi We find the defining polynomial of the unicuspidal ramphoid quartic by using the parametrization given in chapter \ref{rccq} and eliminating $s$ and $t$ with {\it Singular}.

{\scr
\begin{verbatim}
ring R=0, (x,y,z,a,b,c,d,e,s,t), dp;
ideal A6=x-(a-d),y-(c),z-(e);
ideal ST=a-s4,b-s3t,c-s2t2,d-st3,e-t4;
ideal A6ST=A6,ST;
short=0;
eliminate(std(A6ST),abcdest);
_[1]=y^4-2*x*y^2*z+x^2*z^2-y*z^3
\end{verbatim}
}

\noi We find the defining polynomial of the dual curve $C_8$ in {\it Singular}.

{\scr
\begin{verbatim}
ring r=0,(x,y,z,s,t,u),dp;
poly f=y^4-2*x*y^2*z+x^2*z^2-y*z^3;
ideal I=f,s-diff(f,x),t-diff(f,y),u-diff(f,z);
short=0;
eliminate(std(I),xyz);
_[1]=27*s^5+4*s^2*t^3-144*s^3*t*u-16*t^4*u+128*s*t^2*u^2-256*s^2*u^3
\end{verbatim}
}

\noi The above output is a polynomial in $s,t,u$ defining the rational cuspidal quintic with four cusps. This can be verified by the following code in {\it Maple}.

{\scr
\begin{verbatim}
with(algcurves):
f := 27*s^5+4*s^2*t^3-144*s^3*t*u-16*t^4*u+128*s*t^2*u^2-256*s^2*u^3;
u := 1;
singularities(f, s, t);
\end{verbatim}
\begin{equation*}
\begin{split}
&[[x,y,z],m_p,\delta_p,\# \text{Branches}],\\
&[[-\frac{16}{3},8,1],2,1,1],\\
&[[{\it RootOf} \left( 9\,{{\it \_Z}}^{2}-
48\,{\it \_Z}+256 \right) ,-8+\frac{3}{2}\,{\it RootOf} \left( 9\,{{\it \_Z}}^
{2}-48\,{\it \_Z}+256 \right) ,1],2,1,1],\\
&[[0,0,1],2,3,1].
\end{split}
\end{equation*}
}

\noindent Although the above construction is by far the most elegant one, it is also possible to get $C_8$ by using a Cremona transformation to transform the unicuspidal ramphoid quartic. The easiest way of showing this is to apply the standard Cremona transformation with three proper base points to $C_8$ rotated such that the three simple cusps are placed in $(1:0:0)$, $(0:1:0)$ and $(0:0:1)$. This gives a rational quartic with one single $A_6$-cusp and three inflection points. The inverse transformation gives the quintic. We show the latter transformation implicitly.

\begin{description}
\item
\parpic[l]{\includegraphics[width=0.3\textwidth]{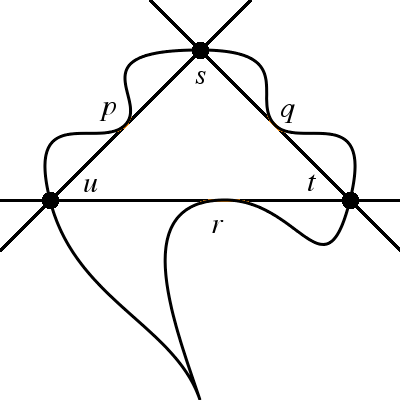}}\;
Let $C$ be the unicuspidal ramphoid quartic. Then rotate this curve such that we have the specific arrangement shown on the left. We have three lines $T^p$, $T^q$ and $T^r$ which are tangent lines to $C$ at three points $p$, $q$ and $r$. The three lines additionally intersect in three points $s$, $t$ and $u \in C$. 
\\
\end{description}

\begin{description}
\item
\parpic[l]{\includegraphics[width=0.3\textwidth, height=0.3\textwidth]{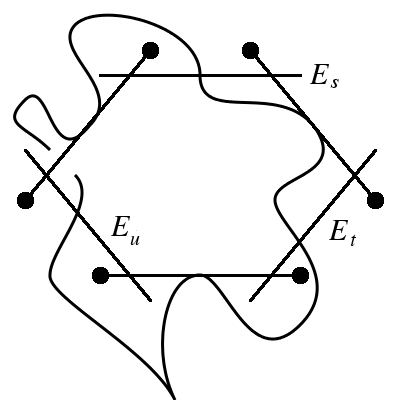}}\;
\noi Using the Cremona transformation $\psi_3(s,t,u)$, we get three exceptional lines $E_s$, $E_t$ and $E_u$, replacing $s$, $t$ and $u$. Furthermore, the transforms of $T^p$, $T^q$ and $T^r$ intersect the transform of $C$ with multiplicity $2$ in the transforms of the points $p$, $q$ and $r$.
\\
\\
\end{description}
\begin{description}
\item
\parpic[l]{\includegraphics[width=0.3\textwidth]{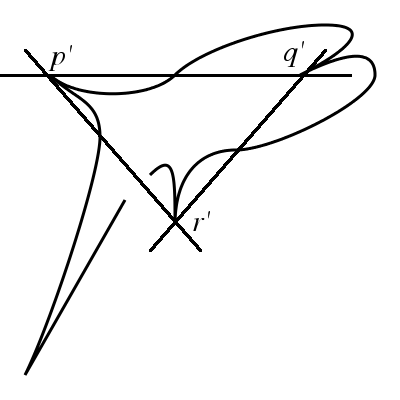}}\;
\noi Blowing down the transforms of the tangent lines results in points $p'$, $q'$ and $r'$ on $C'$, which by elementary properties of the blowing-down process have multiplicity $$m_{p'}=m_{q'}=m_{r'}=2.$$ Since the points originally were smooth points on $C$, $p'$, $q'$ and $r'$ are cusps with multiplicity sequence $(2)$. Notice that the cusp of $C$ is unaffected by the Cremona transformation.
\\
\end{description}

\noi The degree of $C'$ is 
\begin{equation*}
\begin{split}
d'=&2 \cdot d-m_s-m_t-m_q\\
=&2 \cdot 4-1-1-1\\
=&5.
\end{split}
\end{equation*}

\noi Thus, we have constructed $C_8$, a rational cuspidal quintic with cuspidal configuration $[(2_3),(2),(2),(2)]$. \\

{\samepage
\noi The curve $C_{8}$ is given by the below parametrization, and it has the following properties.

\begin{center}
	\setlength{\extrarowheight}{2pt}
	{\begin{tabular}{ccc}
	\multicolumn{3}{c}{$(s^4t:s^2t^3-s^5:t^5+2s^3t^2)$}\\
			\multicolumn{3}{c}{} \\ 
	\hline
	\multicolumn{3}{c}{{\bf \# Cusps = 	4 }}\\
	{\bf Cusp $p_j$} & {\bf $(C \cdot T_{p_j})_{p_j}$}& {\bf $(C \cdot H_C)_{p_j}$ }\\
	\hline 
		$(2_3)$& 4 & 21\\
	$(2)$& 3 & 8\\
	$(2)$& 3 & 8\\
	$(2)$& 3 & 8\\
		\hline
		&&\\
		\hline
			\multicolumn{3}{c}{{\bf \# Inflection points =	0 }}\\
	\hline
	\end{tabular}}
	\end{center}
}

\chapter{More cuspidal curves}\label{mcc}
The rational cuspidal cubics, quartics and quintics presented in this thesis, have been known for a while because of Namba's classification of curves in \cite{Namba}. For some time there were only a few known examples of rational cuspidal curves of degree $d\geq6$. In recent years this has changed. In this chapter we will present more cuspidal curves. 

\section{Binomial cuspidal curves}
Amongst the most simple cuspidal curves are curves which are given as the zero set of a binomial, a homogeneous polynomial with two terms, 
$$F_{(mn)}=z^{n-m}y^m-x^n, \qquad n>m \geq 1,\; \mathrm{gcd}(m,n)=1.$$

\noi These curves are cuspidal and are called {\em binomial curves}. They are unicuspidal if $m=1$ or symmetrically $n-m=1$, and they are bicuspidal if $m\geq2$. In the first case, the curves have one cusp and one inflection point of type $n-2$, which can be verified by direct calculation. In the latter case, the curves have one cusp in $p=(0:0:1)$ with multiplicity $m_p=m$, one cusp in $q=(0:1:0)$ with multiplicity $m_q=n-m$, and no inflection points. The two cusps have different multiplicity sequences, but are quite similar in the sense that they can be investigated by identical methods and arguments. Note that since $m$ and $n$ are relatively prime, so are $n-m$ and $n$. The observation that a general bicuspidal binomial curve does not have inflection points can be verified by a direct calculation of the Hessian curve and the intersection multiplicity of the Hessian curve and the curve at the cusps. We will come back to this in Section \ref{ICHC}.\\

\noi A binomial curve can be parametrized by $$(s^{n-m}t^m:t^n:s^n).$$ The projection center $V$, by which it can be projected from the rational normal curve in $\mathbb{P}^n$, is given by
$$V=\V(x_0,x_m,x_n),$$
$$A_V=\left[ \begin {array}{ccccccc} 1&0&\ldots&0&\ldots&0\\\noalign{\medskip}0&0&\ldots&1&\ldots&0
\\\noalign{\medskip}0&0&\ldots&0&\ldots&1\end {array} \right].$$\\

\noi We want to investigate a cusp on a binomial curve. Because of the symmetry of the unicuspidal binomial curves we may always assume that the binomial curve $C$ has a cusp in $p=(0:0:1)$. $C$ can then be represented around $p$ by the affine curve given by $f$, $$f=y^m-x^n.$$ The cusp has multiplicity $m$, and the multiplicity sequence can be calculated by the algorithm given in Theorem \ref{thm:mp} on page \pageref{thm:mp}. The cusp has tangent $\V(y)$, which intersects the curve with multiplicity $(T_p \cdot C)_p=n.$ Hence, the curve can locally around the cusp be parametrized by $$(C,p)=(t^m:t^n:1).$$

\subsubsection{Fibonacci curves}\label{def:fibcurves} 
We define the {\em Fibonacci numbers} $\f_k$ recursively. For $k \in \mathbb{N}$ we define \begin{equation*} \begin{split} \f_0 &= 0,\\ \f_1 &= 1, \\\f_{k+1}&=\f_{k}+\f_{k-1}. \end{split} \end{equation*}

\noi With $\f_k$ as above we define a subseries of binomial cuspidal curves called the {\em Fibonacci curves} \cite{Rolv}. The $k$th Fibonacci curve is defined as $$C_k=\mathcal{V}(y^{\f_{k}}-x^{\f_{k-1}}z^{\f_{k-2}}), \qquad \text{for all }k \geq 2.$$  This curve has degree $\f_k$. For $k=2$ we have a line. For $k=3$ we have an irreducible conic. For $k=4$ the curve is a cuspidal cubic with a simple cusp in $(0:0:1)$. For all $k \geq 5$ the curve $C_k$ has two cusps. The two cusps are located in $(1:0:0)$ and $(0:0:1)$. They have multiplicity sequences $(\f_{k-2},\verb#{#\f_{i}\verb#}#_{i=1}^{k-2})$ and $(\verb#{#\f_{i}\verb#}#_{i=1}^{k-1})$ respectively. The curves are bicuspidal since any two successive Fibonacci numbers are relatively prime integers. If they did have a common factor $\neq 1$, so would their sum and their difference, which is the next and the previous number in the series. This indicates that all Fibonacci numbers have a common factor $\neq 1$, but that is certainly not true.

These curves will be revisited in Section \ref{fibprod}.  
  
\subsubsection{Semi-binomial curves}
\noi Carefully adding a few terms to the defining polynomial of the binomial curves produces series of curves strongly related to the binomial curves. Adding terms can be done such that it simply corresponds to a linear change of coordinates. These curves provide nothing new. However, it is possible to add terms in such a way that the resulting curves have precisely the same cuspidal configuration as the original binomial curve, but such that they are not projectively equivalent. This can happen because the introduction of new terms in the defining polynomial sometimes leads to inflection points. Because they are strongly related to the binomial curves, we choose to call these curves {\em semi-binomial curves}. 

We have presented several examples of this phenomenon in this thesis. For example, the ovoid cuspidal quartic with one inflection point is binomial, $$C_{4A}=\V(z^4-xy^3).$$ The ovoid cuspidal quartic with two inflection point is semi-binomial, $$C_{4B}=\V(x^3y-z^3x+z^4).$$

\section{Orevkov curves}
In \cite{Ore}, Orevkov gave a proof of the existence of certain unicuspidal rational curves, which from here on will be referred to as {\em Orevkov curves}. The proof explains a way of constructing series of such curves explicitly by applying a product of Cremona transformations to simple, nonsingular algebraic curves of low degree. In the following we will describe the series of Orevkov curves and their construction.

Let $\f_k$ denote the $k$th Fibonacci number. An algebraic curve $C_k$ of degree $d_k=\f_{k+2}$ with a single cusp of multiplicity $m_k=\f_{k}$ is called an {\em Orevkov curve}. Moreover, an algebraic curve $C^{\star}_k$ of degree $d^{\star}_k=2\f_{k+2}$ with a single cusp of multiplicity $m^{\star}_k=2\f_k$ is also called an {\em Orevkov curve}. 

Let $N$ be a rational cubic curve with a singularity with two distinct tangents. The singularity is commonly referred to as a node, and the curve itself is called a nodal cubic. Furthermore, let $C_{-3}$ and $C_{-1}$ be the tangents to the branches of $N$ at the singular point, and let $C_{0}$ be the inflectional tangent. Each of the three lines intersects $N$ in only one point. Additionally, let $C^{\star}_{0}$ be a conic which intersects $N$ in one smooth point. Then a Cremona transformation $\psi$ which is biregular on $\mathbb{P}^2 \setminus N$ will transform the mentioned nonsingular curves into unicuspidal Orevkov curves.

\pagebreak
\begin{figure}[H]\centering
\renewcommand{\thesubfigure}{}
\subfigure[$N$, $C_{-3}$ and $C_{-1}$]{\includegraphics[width=0.27\textwidth]{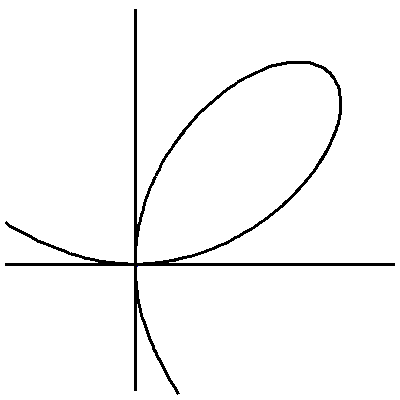}} \qquad
\subfigure[$N$ and $C_0$]{\includegraphics[width=0.27\textwidth]{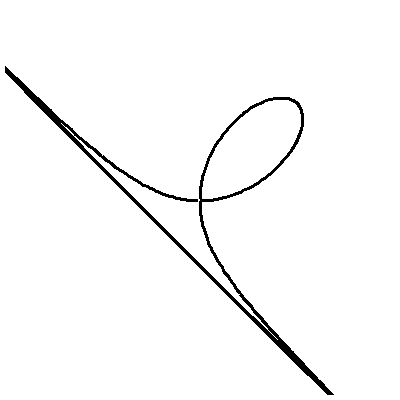}} \qquad
\subfigure[$N$ and $C_{0}^{\star}$]{\includegraphics[width=0.27\textwidth]{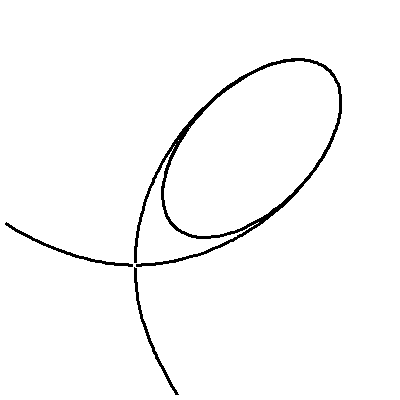}} 
\caption{Initial curves in the series of Orevkov curves, $z=1$.}
\label{FIG:O}
\end{figure}
\noi Note that in Figure \ref{FIG:O}, $N$ and $C_0$ intersect in the point $(1:1:0)$, which is not in the chosen affine covering.\\

\noi We define recursively,
\begin{equation*}
\begin{split}
&C_k=\psi(C_{k-4}), \qquad k \geq 3, \qquad k \neq 2 \;{(\mathrm{mod}\,4)}, \\ 
&C^{\star}_k=\psi(C^{\star}_{k-4}), \qquad k > 0, \qquad k \equiv 0 \;{(\mathrm{mod}\,4)}.
\end{split}
\end{equation*}

\begin{thm} \textup{\cite[p.658]{Ore}} The following series of Orevkov curves exist.
 \begin{description}
      \item{\qquad Orevkov I} -- $C_k$ for any $k > 1$, $k \equiv 1 \;(\mathrm{mod}\,4)$,
      \item{\qquad Orevkov II} -- $C_k$ for any $k > 0$, $k \equiv 3 \;(\mathrm{mod}\,4)$,
      \item{\qquad Orevkov III} -- $C_k$ for any $k > 0$, $k \equiv 0 \;(\mathrm{mod}\,4)$,
      \item{\qquad Orevkov IV} -- $C^{\star}_k$ for any $k > 0$, $k \equiv 0 \;(\mathrm{mod}\,4)$.\\
\end{description}
\end{thm}

\noi Choose $$N=\V(xyz-x^3-y^3).$$ Then the Orevkov curves can explicitly be constructed by the Cremona transformation $$\psi=\sigma_5 \circ \sigma_4 \circ \sigma_3 \circ \sigma_2 \circ \sigma_1.$$

{\small
\begin{center}
  \setlength{\extrarowheight}{2pt}
\begin{tabular}{lllll}
\hline
${\bf \sigma_1}$ & ${\bf \sigma_2}$ &${\bf \sigma_3}$ & ${\bf \sigma_4}$ & ${\bf \sigma_5}$ \\ 	
\hline
$x_1=x^2$ &$ x_2=x_1z_1 $& $x_3=y_2z_2 $&$ x_4=x_3z_3+y_3^2 $&$ x_5=x_4y_4$\\ 	
$y_1=xy $&$ y_2=y_1z_1-x_1^2 $&$ y_3=x_2z_2 $&$ y_4=y_3z_3 $&$ y_5=y_4^2$\\ 	
$z_1=xz-y^2 $&$ z_2=z_1^2 $&$ z_3=x_2y_2 $&$ z_4=z_3^2 $&$ z_5=y_4z_4+x_4^2$\\
\hline
\end{tabular} 
\end{center}	
}

\noindent Note that the five above Cremona transformations are not their own inverses; they are compositions of linear changes of coordinates and the quadratic Cremona transformations given in Chapter \ref{ccc}. Therefore, in applying the above transformations to defining polynomials, we must perform the necessary substitutions with the inverse transformations. The fact that the transformations are not on standard form complicates finding the strict transform. Although this causes problems in the direct calculations given below, we will not discuss this problem here.

{\small
\begin{center}
  \setlength{\extrarowheight}{2pt}
\begin{tabular}{lllll} 
\hline
${\bf \sigma_1^{-1}}$ &$ {\bf \sigma_2^{-1}} $&${\bf \sigma_3^{-1}}$ &$ {\bf \sigma_4^{-1}} $&$ {\bf \sigma_5^{-1}} $\\ 	
\hline
$x=x_1^2$ &$ x_1=x_2z_2 $&$ x_2=y_3z_3 $&$ x_3=x_4z_4-y_4^2 $&$ x_4=x_5y_5$\\ 	
$y=x_1y_1 $&$ y_1=y_2z_2+x_2^2 $&$ y_2=x_3z_3 $&$ y_3=y_4z_4 $&$ y_4=y_5^2$\\ 	
$z=x_1z_1+y_1^2 $&$ z_1=z_2^2 $&$ z_2=x_3y_3 $&$ z_3=z_4^2 $&$ z_4=y_5z_5-x_5^2$\\
\hline
\end{tabular} 	
\end{center}
}\medskip

\noindent {\bf Orevkov I.}\;
Choose the tangent of one branch of the nodal cubic at the node. It can be given by $C_{-3}=\V(F_{-3})$, where $F_{-3}=x$. Applying $\psi$ to this curve results in the curve $C_1=\V(F_1)$, where $F_1=yz-x^2$. This curve is nonsingular, hence not an Orevkov curve. Another application of $\psi$ gives the Orevkov curve $C_5$, a rational cuspidal curve of degree $d_5=\f_7=13$, with a single cusp of multiplicity $m_5=\f_5=5$. Successive applications of $\psi$ produce the series of Orevkov I curves.\medskip

\noindent {\bf Orevkov II.}\;
To produce the second series of Orevkov curves, choose the tangent to the other branch of the nodal cubic at the node. It can be given by $C_{-1}=\V(F_{-1})$, where $F_{-1}=y$. Applying $\psi$ to this curve gives the Orevkov curve $C_3$, a rational cuspidal curve of degree $d_3=\f_5=5$, with a single cusp of multiplicity $m_3=\f_3=2$. Successive applications of $\psi$ produce the series of Orevkov II curves.\medskip

\noindent {\bf Orevkov III.}\;
Produce the third series of Orevkov curves by choosing an inflectional tangent of the nodal cubic $N$. We let the curve $C_0=\V(F_0)$, where $F_{0}=3x+3y+z$, be this inflectional tangent. Applying $\psi$ to this curve gives the Orevkov curve $C_4$, a rational cuspidal curve of degree $d_4=\f_6=8$, with a single cusp of multiplicity $m_4=\f_4=3$. Successive applications of $\psi$ produce the series of Orevkov III curves.\medskip

\noi The cusp of an Orevkov curve $C_k$ of any of these three series has multiplicity sequence $$\ol{m}_k=(\f_k,S_k,S_{k-4},\ldots,S_{\nu}),$$ where $\nu=3,4,5$ is determined by $k=4j+\nu, \;j \in \mathbb{N}_0$, and $S_i$ denotes the subsequence $$S_i=((\f_i)_5,\f_i-\f_{i-4}).$$ 

\noindent {\bf Orevkov IV.}\;
To get the fourth series of Orevkov curves, choose a conic $C^{\star}_{0}=\V(F^{\star}_{0}),$ where $F^{\star}_{0}=21x^2-22xy+21y^2-6xz-6yz+z^2$. This particular conic intersects the nodal cubic in exactly one smooth point. Applying $\psi$ to this curve gives the Orevkov curve $C^{\star}_4$, a rational cuspidal curve of degree $d^{\star}_4=2\f_6=16$, with a single cusp of multiplicity $m^{\star}_4=2\f_4=6$. Successive applications of $\psi$ produce the series of Orevkov IV curves.\medskip

\noi The cusp of an Orevkov curve $C^{\star}_k$ has almost the same multiplicity sequence as the cusp of the curve $C_k$. The only difference is that every multiplicity is multiplied by 2.\\

\begin{rem}
Note that the birationality of the Cremona transformations allows the application of the inverse transformation $\psi^{-1}$ to all Orevkov curves. This immediately implies that all Orevkov curves of type I, II and III can be transformed into lines by a sequence of Cremona transformations. This is also the case for Orevkov curves of type IV. These curves can be transformed into an irreducible conic, and any such curve can in turn be transformed into a line. Hence, all Orevkov curves are rectifiable.
\end{rem}

\section{Other uni- and bicuspidal curves}
The search for more cuspidal curves led Fenske in \cite{Fen99b} to the discovery of essentially eight different series of rational uni- and bicuspidal curves. The curves were found using suitable Cremona transformations to transform the following binomial and semi-binomial cuspidal curves of degree $d$,$$\V(xy^{d-1}-z^d) \qquad \text{and} \qquad \V(xy^{d-1}-z^d-yz^{d-1}).$$ An overview is given in Table \ref{tab:fenskeunibi}. Note that since the curves $C_i$ are strict transforms of the above curves, the degree of $C_i$ is given as a function of $d$. 

\begin{table}[htb]
  \renewcommand\thesubtable{}
  \setlength{\extrarowheight}{2pt}
\centering
	{\begin{tabular}{ccl}
	\hline
	{\bf Curve} &{\bf $\deg C_{a,d}$}&{\bf Cuspidal configuration}\\ 
	\hline 
	$C_1$&$da+d$&$(da,d_{a+b},d-1),(d_{a-b})$\\
	$C_{1a}$&$da+d$&$(da,d_{2a},d-1)$\\
	$C_2$&$da+d$&$(da,d_{a+b}),(d_{a-b},d-1)$\\
	$C_{2a}$&$da+d$&$(da,d_{2a}),(d-1)$\\
	$C_3$&$da+d+1$&$(da+1,d_a),(d_{a+1})$\\
	$C_4$&$da+d+1$&$(da,d_{a+1}),((d+1)_a)$\\
	$C_5$&$da+d+1$&$(da,d_a),((d+1)_a,d)$\\
	$C_6$&$da+d+2$&$(da+1,d_a),((d+a)_{a+1})$\\
	$C_7$&$da+2d-1$&$(da+d-1,d_a,d-1),(d_{a+1},d-1)$\\
	$C_8$&$a+2$&$(a),(2_a)$\\
	\hline
	\multicolumn{3}{l}{\footnotesize{$d \geq 2$ and $0\leq b<a$ integers.}}
	\end{tabular}}\\
		\caption {Fenske's uni- and bicuspidal curves \cite[Thm. 1.1., p.310]{Fen99b}.}
	\label{tab:fenskeunibi}
	\end{table}

\section{Tricuspidal curves}
Curves with three or more cusps are quite rare. We have three series of tricuspidal curves, which were found by Flenner, Zaidenberg and Fenske in \cite{FlZa95} \cite{FlZa97} \cite{Fen99a}. The series of curves can be produced from simple cuspidal curves of low degree using successive quadratic Cremona transformations. 

\subsection{Curves with $\mu=d-2$}
In \cite{FlZa95}, Flenner and Zaidenberg constructed a series of tricuspidal curves where the maximal multiplicity of any cusp of the curves was $\mu=d-2$. They found that for any degree $d \geq 4$, $a \geq b \geq 1$, with $a+b=d-2$, there is a unique, up to projective equivalence, rational cuspidal curve $C=C_{d,a} \subset \Po$ of degree $d$ with three cusps. The curve $C_{d,a}$ has cuspidal configuration 
$$[(d-2),(2_a),(2_b)].$$
The restrictions on $a$ and $b$ imply that for every $d$ there are exactly $\left \lfloor \frac{d-2}{2} \right \rfloor$ tricuspidal curves of this kind. Furthermore, Flenner and Zaidenberg proved that any rational cuspidal curve with one cusp with multiplicity $d-2$ has at most three cusps. Additionally, they proved that if there are three cusps on the curve, then it has the given cuspidal configuration.

\subsubsection{The parametrization}
In \cite{FlZa95}, the tricuspidal curves with $\mu=d-2$ were painstakingly constructed by parametrization methods. They started out with a general parametrization. Then they imposed restrictions on the parametrization in order to get the desired type of cusps. They found that the curve is parametrized by $$((s-t)^{d-2}s^2:(s-t)^{d-2}t^2:s^2t^2\cdot q(s,t))$$ where the polynomial $q(1,t)$, and thus the homogeneous polynomial $q(s,t)$, can be found by feeding {\em Maple} the below code.

{\scr
\begin{verbatim}
c[k] := (product(2*a-2*i+1, i = 1 .. k))/2^k;
f(t) := expand(1+ sum((c[k])/(k!)(t^(2)-1)^(k),k=1..(d-3)));
q[d,a](t) := sort((factor(f(t)+t^((2 a-1))))/((1+t)^((d-2))));
\end{verbatim}
}

\noi The parametrization enables us to find and investigate the projection center for every such curve. Further investigating these projection centers could be interesting, but goes beyond the scope of this thesis. 

\pagebreak
\subsubsection{Construction by Cremona transformations}
In \cite{FlZa97}, Flenner and Zaidenberg presented an elegant way of constructing the tricuspidal curves using Cremona transformations.

Let $C$ be a curve with cusps $q$, $p$ and $r$, with multiplicity sequences $\ol{m}_q=(d-2)$, $\ol{m}_p=(2_a)$ and $\ol{m}_r=(2_b)$. Then let $L^{pq}$ denote the line through $p$ and $q$, $L^{qr}$ the line through $q$ and $r$, and $T^q$ the tangent line at $q$. 

By B\'{e}zout's theorem and Lemma \ref{intmult}, the tangent line $T^q$ intersects $C$ with multiplicity $d-1$ at $q$. Hence, $T^q$ intersects $C$ transversally in a smooth point $s$. Applying the Cremona transformation $\psi_2(s,q,L^{pq})$ to $C$ results in a curve $C'$ of degree $d+1$ with three cusps with multiplicity sequences $(d-1)$, $(2_{a+1})$ and $(2_b)$.

We may also perform the Cremona transformation $\psi_2(p,q,L^{qr})$ to $C$, this time getting a tricuspidal curve of degree $d$ and multiplicity sequences $(d-2)$, $(2_{a+1})$ and $(2_{b-1})$.

With these transformations, Flenner and Zaidenberg proved that all tricuspidal curves of degree $d$ with $\mu=d-2$ can be constructed from the tricuspidal quartic curve. This implies that all these curves are rectifiable.


\subsection{Curves with $\mu=d-3$}
The work on the tricuspidal curves with maximal multiplicity $\mu=d-2$ was soon followed by a parallel result for tricuspidal curves with a cusp with $\mu=d-3$. In \cite{FlZa97}, Flenner and Zaidenberg proved that for a tricuspidal curve of degree $d=2a+3$ with a cusp with maximal multiplicity $\mu=d-3$, the only possible cuspidal configuration is $$[(2a,2_a),(3_a),(2)].$$ Furthermore, they proved that for each $a \geq 1$ there exists such a rational cuspidal curve. The existence was proved performing a sequence of Cremona transformations to the rational cuspidal cubic and constructing the curves inductively. The construction proves that the curves are rectifiable. It was also proved that the curves are unique up to projective equivalence.\\

\begin{rem}
Observe that when $a=1$, then $d=5$. Then the curve has the same cuspidal configuration as curves in this series, but the maximal multiplicity of the cusps is $\mu=3=5-2$. Hence, this curve does by definition not belong to this series. Rather, it is an example of a curve from the series of curves with $\mu=d-2$.
\end{rem}


\pagebreak
\subsection{Curves with $\mu=d-4$}
The last series of tricuspidal curves was presented by Fenske in \cite{Fen99a}. In this article, Fenske found that for every $d=3a+4$, $a\geq 1,$ there exists a tricuspidal curve with cuspidal configuration $$[(3a,3_a),(4a,2_2),(2)].$$ The existence was proved by induction and performing a sequence of Cremona transformation to the bicuspidal quartic with one inflection point. Fenske also found that these curves are unique up to projective equivalence and that they are rectifiable. 

Fenske was {\em not} able to generally prove that if a curve of degree $d$ has three cusps and $\mu=d-4$, then the curve has this cuspidal configuration. However, he was able to prove this result under the condition that $\chi(\Theta_V \langle D \rangle) \leq 0$, where $\chi(\Theta_V \langle D \rangle)$ is the {\em Euler characteristic} and $D$ is the divisor with simple normal crossings of the minimal embedded resolution of the singularities of $C$.\\

\noi By \cite[Cor. 2.5., p.446]{FlZa95}, we can compute the Euler characteristic $\chi(\Theta_V \langle D \rangle)$ of a rational cuspidal curve using the degree and the multiplicity sequences of the cusps of the curve.
$$\chi(\Theta_V \langle D \rangle)=-3(d-3)+\sum_{p \in \mathrm{Sing}\,C} \bar{M}_p,$$ where $\bar{M}_p=\eta_p+\omega_p-1$, $$\eta_p=\sum_{i=0}^{n_p}(m_{p.i}-1)\qquad \text{and} \qquad \omega_p=\sum_{i=1}^{n_p}\left(\left \lceil \frac{m_{p.i-1}}{m_{p.i}} \right \rceil-1\right).$$ \medskip


\noi There are, so far, no known rational cuspidal curves with $\chi(\Theta_V \langle D \rangle) > 0$. Moreover, $\chi(\Theta_V \langle D \rangle) \leq 0$ is a consequence of a conjecture given by Fern�ndez de Bobadilla et al. in \cite[Conj. 4.1., pp.420--424]{Bobadillains}. 
\begin{conj}
Let $C$ be a rational cuspidal curve with $\bar{M}_p$ and $d$ as above. Furthermore, let $\dim{\mathrm{Stab}}_{PGL(3)}(C)$ denote the dimension of the group of transformations in $PGL_3(\mathbb{C})$ which do not move $C$. Then $$\sum_{p \in \mathrm{Sing}\,C} \bar{M}_p \leq 3d-9+\dim{\mathrm{Stab}}_{PGL(3)}(C).$$  
\end{conj}

\begin{rem}
Observe that when $a=1$, then $d=7$. Then the curve has the same cuspidal configuration as the curves in this series, but the maximal multiplicity of the cusps is $\mu=4=7-3$. Hence, this curve does by definition not belong to this series. Rather, it is an example of a curve from the series of curves with $\mu=d-3$.
\end{rem}


\pagebreak
\subsection{Overview}
We have three series of rational tricuspidal curves. For curves with degree $d \geq 6$, assuming $\chi(\Theta_V \langle D \rangle) \leq 0$, these are actually the only tricuspidal curves with these cuspidal configurations, up to projective equivalence. 
{\small
\begin{table}[H]
\begin{center}
\setlength{\extrarowheight}{3pt}
{\small
		\begin{tabular}[c]{ccccccc}
			\hline
			{\bf Series}&$d$ & $\bar{m}_p$ & $\bar{m}_q$ & $\bar{m}_r$ & {\bf Valid for} &  \\
			\hline
			I&$d$ & $(d-2)$ & $(2_a)$ & $(2_{d-2-a})$ & $d \geq 4$ & $d-3 \geq a \geq 1$\\
			II&$2a+3$ & $(d-3, 2_a)$ & $(3_a)$ & $(2)$ & $d \geq 5$ & $a \geq 1$\\
			III&$3a+4$ &$(d-4,3_a)$ & $(4a,2_2)$ & $(2)$ & $d \geq 7$ & $a \geq 1$\\
			\hline
		\end{tabular}}
		\caption{Tricuspidal curves.}
	\label{tab:TricuspidalCurves}
	\end{center}
	\end{table}}
\noi These results allow us to count the number $N$ of tricuspidal curves for each degree $d \geq 4$, up to projective equivalence. For degrees $d=4$ and $d=5$, we have seen that $N=1$ and $N=2$. For degree $d=6$ we see in Table \ref{tab:TricuspidalCurves} that we know $N=2$ tricuspidal curves, and for $d=7$ we know $N=3$ tricuspidal curves. For any degree $d \geq 8, \; d \equiv k \;(\mathrm{mod}\,6)$ the number $N$ of known tricuspidal curves is given in Table \ref{tab:NOTricuspidalCurves}. To simplify notation we write $N_0=\left \lfloor \frac{d-2}{2} \right \rfloor,$ so that $N_0$ is the number of curves in series I for each $d$.
\begin{table}[htb]
\begin{center}
\setlength{\extrarowheight}{3pt}
		\begin{tabular}[c]{ccc}
			\hline
			$k$ & {\bf \;\# Tricuspidal curves $N$\;} & {\bf Series represented}\\
			\hline
			$0,2$ & $N_0$ & I\\
			$1$ & $N_0+2$ &I, II, III\\
			$3,5$ & $N_0+1$ & I, II\\
			$4$ & $N_0+1$ &I, III\\
			\hline
		\end{tabular}
		\caption{The number of known tricuspidal curves for $d\geq 8,\; d \equiv k \;(\mathrm{mod}\,6)$.}
	\label{tab:NOTricuspidalCurves}
	\end{center}
	\end{table}

\section{Rational cuspidal sextics}
The search for rational cuspidal curves of a given degree is made slightly easier with all the above results. With the results it was possible for Fenske in \cite[Cor. 1.5., p.312]{Fen99b} to present a list of all existing cuspidal configurations of rational cuspidal sextics. He also gave explicit parametrizations of all rational cuspidal sextic curves with one and two cusps, up to projective equivalence \cite[pp.327--328]{Fen99b}. For the rational cuspidal sextic curves with three cusps, the explicit parametrizations of the curves are given by Flenner and Zaidenberg in \cite[Thm. 3.5., p.448]{FlZa95}. 

\pagebreak
\begin{table}[H]
  \renewcommand\thesubtable{}
  \setlength{\extrarowheight}{2pt}
\centering
	{\begin{tabular}{ccll}
	\hline
	{\bf \# Cusps}&	{\bf Curve} &{\bf Cuspidal configuration} & {\bf \# Curves}\\
	\hline 
	\multirow{3}{14mm}{1}&$C_1$&$\qquad(5)$&4 -- ABCD\\
	&$C_2$&$\qquad (4,2_4)$&2 -- AB\\
	&$C_3$&$\qquad (3_3,2)$&3 -- ABC\\
	\hline
	\multirow{6}{14mm}{2}&$C_4$&$\qquad (3_3),(2)$&2 -- AB\\
	&$C_5$&$\qquad (3_2,2),(3)$&1\\
	&$C_6$&$\qquad (3_2),(3,2)$&1\\
	&$C_7$&$\qquad (4,2_3),(2)$&1\\
	&$C_8$&$\qquad (4,2_2),(2_2)$&1\\
	&$C_9$&$\qquad (4),(2_4)$&1\\
		\hline
		\multirow{2}{14mm}{3}&$C_{10}$&$\qquad (4),(2_3),(2)$&1\\
	&$C_{11}$&$\qquad (4),(2_2),(2_2)$&1\\
	\hline
	\end{tabular}}
	\caption {Rational cuspidal sextic curves \cite[Cor. 1.5., p.312]{Fen99b}.}
	\label{tab:degree6}
	\end{table}

\chapter{On the number of cusps}\label{onc}
All the examples of rational cuspidal curves presented in this thesis are interesting in themselves. Some are, however, more intriguing than the others in the search for an upper bound on the number of cusps on a rational cuspidal curve. The particularly interesting curves are the curves with three or more cusps. Indeed, we so far only know one rational cuspidal curve with four cusps, the quintic with cuspidal configuration $[(2_3),(2),(2),(2)]$. Furthermore, we only know one tricuspidal curve not contained in table \ref{tab:TricuspidalCurves}, the quintic with cuspidal configuration $[(2_2),(2_2),(2_2)]$.

\section{A conjecture}
The above observation was made by Piontkowski, and in 2007 he proposed the following conjecture after investigating almost all cuspidal curves of degree $\leq 20$ \cite[Conj. 1.4., p.252]{Piontkowski}.
\begin{conj}[On the number of cusps of rational cuspidal curves]
A rational cuspidal plane curve of degree $d \geq 6$ has at most three cusps. The curves of degree $d \geq 6$ with precisely three cusps occur in the three series in Table \ref{tab:TricuspidalCurves} on page \pageref{tab:TricuspidalCurves}.
\end{conj}

\section{An upper bound}
The general research on rational cuspidal curves has not only circled around constructing curves and series of curves. The most recent progress was made in 2005, when Tono lowered the upper bound for the number of cusps on a cuspidal plane curve of genus $g$. Thus, he also found a new upper bound for the number of cusps on a rational cuspidal plane curve \cite[Thm. 1.1., p.216]{Tono05}.

\begin{thm}[An upper bound]
A cuspidal plane curve of genus $g$ has no more than $\frac{21g+17}{2}$ cusps.
\end{thm}

\begin{cor}[An upper bound for rational curves]
A rational cuspidal plane curve has no more than 8 cusps.
\end{cor}

\section{Particularly interesting curves}
Noting Tono's theorem and keeping Piontkowski's conjecture in mind, we are naturally led to the investigation of the two particular quintic curves in Table \ref{tab:interestingrcc}.
\begin{table}[ht!]
\begin{center}
\setlength{\extrarowheight}{2pt}
\begin{tabular}{c c c c}
\hline
{\bf Curve} & {\bf \# Cusps} & {\bf Cuspidal configuration} & {\bf Degree $d$} \\
\hline
$C_{\triangle}$ &$3$& $[(2_2),(2_2),(2_2)]$ & $d=5$\\
$C_{\square}$ &$4$& $[(2_3),(2),(2),(2)]$& $d=5$\\
\hline
\end{tabular}
\end{center}
\caption{Particularly interesting cuspidal curves.}
\label{tab:interestingrcc}
\end{table}

\subsection{All about $C_{\triangle}$}
The quintic curve with cuspidal configuration $[(2_2),(2_2),(2_2)]$ is special since it is the only tricuspidal curve which is not found in any of the three series of tricuspidal curves in Table \ref{tab:TricuspidalCurves}. Apart from that, however, it is hard to find ways in which this curve stands out. In Chapter \ref{rcq} we saw that this curve can be constructed by a Cremona transformation of the bicuspidal quartic, and there was nothing worth noting in this construction. Next we will see that projection and dualization do not reveal any secrets of this curve either.

\subsubsection{Projection}
The parametrization given in \cite[Thm. 2.3.10., pp. 179--182]{Namba} gives us the projection center $V$.
$$V=\V(-x_0+x_1,-\textstyle{\frac{5}{32}}x_0+x_3,-\textstyle{\frac{125}{128}}x_0-\textstyle{\frac{25}{16}}x_2-5x_4+x_5),$$
$$A_V= \left[ \begin {array}{cccccc} -1&1&0&0&0&0\\\noalign{\medskip}-{
\textstyle{\frac {5}{32}}}&0&0&1&0&0\\\noalign{\medskip}-{\textstyle{\frac {125}{128}}}&0&-{
\textstyle{\frac {25}{16}}}&0&-5&1\end {array} \right] , \qquad K_V=\begin{bmatrix}0&0&1&0&0&\textstyle{\frac{25}{16}}\\\noalign{\medskip}\textstyle{\frac{32}{5}}&\textstyle{\frac{32}{5}}&0&1&0&\textstyle{\frac{25}{4}} \\ \nom 0&0&0&0&1&5\end{bmatrix}.$$

\noi With this information on the projection center we observe that the curve $C_{\triangle}$ has three cusps $p_i$ with $(T_{p_i} \cdot C)_{p_i}=4$.
\begin{itemize}
\item[--] $V \cap C_5=\emptyset.$
\item[--] $V \cap T_5=\{p_1,p_2,p_3\}$ for $$(s:t)=\{(0:1),(\textstyle{\frac{-10+6\sqrt{5}}{5}}:1),(\textstyle{\frac{-10-6\sqrt{5}}{5}}:1)\}.$$
\item[--] $V \cap O^2(s,t)=\{p_1,p_2,p_3\}$ for $$(s:t)=\{(0:1),(\textstyle{\frac{-10+6\sqrt{5}}{5}}:1),(\textstyle{\frac{-10-6\sqrt{5}}{5}}:1)\}.$$
\item[--] $V \nsubseteq O^2(s,t)$ for any $(s:t)$.
\item[--] $V \subset O^3(s,t)=\{p_1,p_2,p_3\}$ for $$(s:t)=\{(0:1),(\textstyle{\frac{-10+6\sqrt{5}}{5}}:1),(\textstyle{\frac{-10-6\sqrt{5}}{5}}:1)\}.$$
\item[--] $V \nsubseteq O^4(s,t)$ for any $(s:t)$.
\end{itemize}

\subsubsection{The dual curve}
Using the parametrization, we find, with the help of {\em Singular}, that the curve $C_{\triangle}$ can be given by the following defining polynomial.

{\scr
\begin{verbatim}
ring R=0, (x,y,z,a,b,c,d,e,f,s,t), dp;
ideal C=x-(b-a),y-(d-5/32*a),z-(-125/128*a-25/16*c-5e+f);
ideal ST=a-s5,b-s4t,c-s3t2,d-s2t3,e-st4,f-t5;
ideal CS=C,ST;
short=0;
eliminate(std(CS),abcdefst);
_[1]=709375*x^5+4800000*x^4*y+54560000*x^3*y^2-199424000*x^2*y^3+265420800*x*y^4
-8126464*y^5-3360000*x^4*z-17664000*x^3*y*z+18022400*x^2*y^2*z+49807360*x*y^3*z
-1048576*y^4*z+4915200*x^3*z^2-3932160*x^2*y*z^2+2097152*x*y^2*z^2-1048576*x^2*z^3
\end{verbatim}}

\noi Calculating the dual curve with {\em Singular} gives $C_{\triangle}^*$.

{\scr
\begin{verbatim}
ring R=0,(x,y,z,s,t,u),dp;
poly f=709375*x^5+4800000*x^4*y+54560000*x^3*y^2-199424000*x^2*y^3+265420800*x*y^4
-8126464*y^5-3360000*x^4*z-17664000*x^3*y*z+18022400*x^2*y^2*z+49807360*x*y^3*z
-1048576*y^4*z+4915200*x^3*z^2-3932160*x^2*y*z^2+2097152*x*y^2*z^2-1048576*x^2*z^3;
ideal I=f,s-diff(f,x),t-diff(f,y),u-diff(f,z);
short=0;
eliminate(std(I),xyz);
_[1]=65536*s*t^4+12288*t^5-524288*s^2*t^2*u-163840*s*t^3*u+83200*t^4*u+1048576*s^3*u^2
+19333120*s^2*t*u^2+972800*s*t^2*u^2-256000*t^3*u^2-101888000*s^2*u^3+2560000*s*t*u^3
-2900000*t^2*u^3-152400000*s*u^4-17250000*t*u^4-66796875*u^5
\end{verbatim}}

\noi With {\em Maple} we find that $C_{\triangle}^*$ has exactly the same cuspidal configuration as $C_{\triangle}$.

{\scr
\begin{verbatim}
with(algcurves);
f := 65536*s*t^4+12288*t^5-524288*s^2*t^2*u-163840*s*t^3*u+83200*t^4*u+1048576*s^3*u^2
+19333120*s^2*t*u^2+972800*s*t^2*u^2-256000*t^3*u^2-101888000*s^2*u^3+2560000*s*t*u^3
-2900000*t^2*u^3-152400000*s*u^4-17250000*t*u^4-66796875*u^5
u := 1
singularities(f, s, t);
\end{verbatim}}
{\tiny
\begin{equation*}
\begin{split}
&[[1,0,0],2,2,1],\\
&[[{\it RootOf} \left( 65536\,{{\it \_Z}}^{2}+
416000\,{\it \_Z}+90625 \right) ,-{\frac {32}{15}}\,{\it RootOf}
 \left( 65536\,{{\it \_Z}}^{2}+416000\,{\it \_Z}+90625 \right) -{
\frac {125}{24}},1],2,2,1].
\end{split}
\end{equation*}
}

\subsection{All about $C_{\square}$}
The curve $C_{\square}$ is truly unique since it is the only known rational cuspidal curve with more than three cusps. We have already seen that this curve is the dual curve of the unicuspidal ramphoid quartic, and that it can be constructed from this curve by a Cremona transformation with three proper base points. Unfortunately, there are no surprising properties to be found when analyzing the curve from the perspective of projection.

\subsubsection{Projection} 
The parametrization given in \cite[Thm. 2.3.10., pp. 179--182]{Namba} gives us the projection center $V$.
$$V=\V(x_1,-x_0+x_3,2x_2+x_5),$$
$$A_V=  \left[ \begin {array}{cccccc} 0&1&0&0&0&0\\\noalign{\medskip}-1&0&0&1
&0&0\\\noalign{\medskip}0&0&2&0&0&1\end {array} \right] , \qquad K_V=\begin{bmatrix}0&0&0&0&1&0\\\noalign{\medskip}1&0&0&1&0&0 \\ \nom 0&0&1&0&0&-2\end{bmatrix}.$$

\noi With this information on the projection center we observe that the curve $C_{\square}$ has four cusps $p_1$, $p_2$, $p_3$ and $p_4$. For the three cusps $p_2$, $p_3$ and $p_4$, with multiplicity sequences $(2)$, $(T_{p_i} \cdot C)_{p_i}=3$. For the ramphoid cusp $p_1$, $(T_{p_1} \cdot C)_{p_1}=4$.
\begin{itemize}
\item[--] $V \cap C_5=\emptyset.$
\item[--] $V \cap T_5=\{p_1,p_2,p_3,p_4\}$ for $$(s:t)=\{(0:1),(-4^{-\frac{1}{3}}:1),(4^{-\frac{1}{3}}e^{\frac{\pi i}{3}}:1),(4^{-\frac{1}{3}}e^{\frac{5\pi i}{3}}:1)\}.$$
\item[--] $V \cap O^2(s,t)=\{p_1,p_2,p_3,p_4\}$ for $$(s:t)=\{(0:1),(-4^{-\frac{1}{3}}:1),(4^{-\frac{1}{3}}e^{\frac{\pi i}{3}}:1),(4^{-\frac{1}{3}}e^{\frac{5\pi i}{3}}:1)\}.$$
\item[--] $V \nsubseteq O^2(s,t)$ for any $(s:t)$.
\item[--] $V \subset O^3(s,t)=\{p_1\}$ for $(s:t)=\{(0:1)\}.$
\item[--] $V \nsubseteq O^4(s,t)$ for any $(s:t)$.
\end{itemize}

\pagebreak
\subsubsection{The dual curve}
Using the parametrization, we find, with the help of {\em Singular}, that $C_{\square}$ can be given by the following defining polynomial.

{\scr
\begin{verbatim}
ring R=0, (x,y,z,a,b,c,d,e,f,s,t), dp;
ideal C=x-(b),y-(d-a),z-(f+2c);
ideal ST=a-s5,b-s4t,c-s3t2,d-s2t3,e-st4,f-t5;
ideal CS=C,ST;
short=0;
eliminate(std(CS),abcdefst);
_[1]=27*x^5-2*x^2*y^3+18*x^3*y*z-y^4*z+2*x*y^2*z^2-x^2*z^3
\end{verbatim}}

\noi Calculating the dual curve with {\em Singular} gives $C_{\square}^*$.

{\scr
\begin{verbatim}
ring R=0,(x,y,z,s,t,u),dp;
poly f=27*x^5-2*x^2*y^3+18*x^3*y*z-y^4*z+2*x*y^2*z^2-x^2*z^3;
ideal I=f,s-diff(f,x),t-diff(f,y),u-diff(f,z);
short=0;
eliminate(std(I),xyz);
t^4-8*s*t^2*u+16*s^2*u^2+128*t*u^3
\end{verbatim}}

\noi With {\em Maple} we find that $C_{\square}^*$ is the unicuspidal ramphoid quartic with three inflection points.

{\scr
\begin{verbatim}
with(algcurves);
f := t^4-8*s*t^2*u+16*s^2*u^2+128*t*u^3
singularities(f, s, t);
\end{verbatim}
$$[[1,0,0],2,3,1] $$}

\section{Projections and possibilities}
In the search for an upper bound on the number of cusps of a rational cuspidal curve in $\mathbb{P}^2$, we are led to investigate the problem from $\mathbb{P}^n$ with the language of projections. If we could find a way to estimate the maximal number of intersections of a $n-3$-dimensional projection center $V$ of a cuspidal projection and the tangent developable $T_n$ of the rational normal curve, then we would have an upper bound on the number of cusps of a rational curve of degree $d=n$.

We give the projection center $V$ as the intersection of the zero set of three linearly independent linear polynomials, three hyperplanes. The hyperplanes are often represented by the coefficient matrix $A_V$.
$$V=\V(\sum_{k=0}^na_{0k}x_k, \sum_{k=0}^na_{1k}x_k, \sum_{k=0}^na_{2k}x_k),$$
$$A_V=\begin{bmatrix} a_{00}&a_{01}&a_{02}& \ldots &a_{0n}\\\noalign{\medskip}  a_{10}&a_{11}&a_{12}& \ldots &a_{1n}\\\noalign{\medskip} a_{20}&a_{21}&a_{22}& \ldots &a_{2n}\end{bmatrix}.$$

\noi Since we for any $n$ can find the tangent developable $T_n$ on polynomial form by elimination, it should be possible to intersect it with the three hyperplanes of $V$. In theory this sounds promising. Unfortunately, the number of variables and constants soon gets out of hand. And additionally, a general result is hard to extract, since we have to deal with each degree separately.

Although it does not solve the above problems, we are able to slightly improve this result. We know that the tangent developable $T_n$ of the rational normal curve $C_n$ has similar properties for all values $(s:t)$ by the homogeneity of $C_n$. Therefore, we may freely fix one of the intersection points of the projection center and $T_n$. Let this point be $(0:1:0:\ldots:0)$. Then the projection center $V_1$ of this projection can be represented by the coefficient matrix $A_{V_1}$,

$$A_{V_1}=\begin{bmatrix} a_{00}&0&a_{02}& \ldots &a_{0n}\\\noalign{\medskip}  a_{10}&0&a_{12}& \ldots &a_{1n}\\\noalign{\medskip} a_{20}&0&a_{22}& \ldots &a_{2n}\end{bmatrix}.$$

\noi Keeping Piontkowski's conjecture in mind, we propose the following conjecture concerning the number of intersection points of a projection center and the tangent developable.

\begin{conj}
Let $V$ be a projection center of a cuspidal projection from $\mathbb{P}^n$ to $\mathbb{P}^2$ such that $V$ intersects neither the rational normal curve $C_n$ nor the secant variety $S_n$ outside the tangent developable $T_n$. Then $V$ intersects the tangent developable $T_n$ in maximally three points for all $n\geq4,\; n\neq 5$. For $n=5$, the maximal number of intersection points is four.
\end{conj}

\noi Notice that we, assuming Piontkowski's conjecture and recalling the information given in Table \ref{tab:NOTricuspidalCurves} on page \pageref{tab:NOTricuspidalCurves}, additionally know how many different kinds of intersections are possible for each $n$.

\chapter{Miscellaneous related results}\label{mrr}

\section{Cusps with real coordinates}
An interesting question concerning rational cuspidal curves is whether or not all the cusps of a cuspidal curve can have real coordinates.

For all cuspidal curves with three cusps or less it is elementary that we, by a linear change of coordinates, can assign real coordinates to all the cusps.

For $C_{\square}$, the only known curve with more than three cusps, an answer to the above question is much harder to find. Although $C_{\square}$ in Chapter \ref{onc} was presented with two cusps with real and two cusps with complex coordinates, we can not a priori exclude the possibility that there might exist a linear change of coordinates that will give us a curve where all four cusps have real coordinates. However, we can prove a partial result.\\

\noi We call a curve $C=\V(F)$ {\em real} if $F[x,y,z] \in \mathbb{R}[x,y,z]$. The real image of $C$, denoted by $C(\mathbb{R})$, is defined as $C(\mathbb{R})=\V(F) \cap \mathbb{P}^2_{\mathbb{R}}.$

We will prove that if $C_{\square}$ is a real curve, then all the cusps can not have real coordinates. Let $C_{\square}$ be a real curve.\\

\noi For real algebraic curves the Klein--Schuh theorem holds \cite[Thm. 3.2.2., p.23]{Heidi}.


\begin{thm}[Klein--Schuh]
Let $C$ be a real algebraic curve of degree $d$ with real singularities $p_j \in C(\mathbb{R}),\; j=1,\ldots,s_p$. Let $m_j$ denote the multiplicity of $p_j$ and let $b_j$ be the number of real branches of $C$ at $p_j$. Let $C^*$ be the dual curve of $C$ with degree $d^*$ and real singularities $q_i,\; i=1,\ldots,s_q$. Let $m_i$ denote the multiplicity of $q_i$, and let $b_i$ denote the number of real branches of $C^*$ at $q_i$. For every real algebraic curve $C \subset \mathbb{P}^2$ with dual curve $C^* \subset \mathbb{P}^{2*}$ we have
$$d- \sum_{p_j \in \mathrm{Sing}\,C(\mathbb{R})}(m_j-b_j) = d^*-\sum_{q_i \in \mathrm{Sing}\,C^*(\mathbb{R})}(m_i-b_i).$$
\end{thm}

\noi We know that $C_{\square}$ is a quintic with cuspidal configuration $[(2_3),(2),(2),(2)]$. We also know that $C_{\square}^*$ is a quartic with cuspidal configuration $[(2_3)]$. Assume that all cusps on the curve $C_{\square}$ has real coordinates. Since $C_{\square}$ is a real curve, the cusp on $C_{\square}^*$ has real coordinates. Then $C_{\square}$ contradicts the Klein--Schuh theorem, 
\begin{equation*}
\begin{split}
5 - (2-1)-3\cdot(2-1) &\neq 4-(2-1),\\
2 &\neq 4.
\end{split}
\end{equation*}
Hence, all cusps on $C_{\square}$ can not be real.



\section{Intersecting a curve and its Hessian curve}\label{ICHC}

\noi Calculating the intersection points and intersection multiplicities of a curve $C$ and its Hessian curve $H_C$ is a classical problem. Nevertheless, if we wish to calculate the intersection multiplicity of a curve and its Hessian curve at a point, we are often overwhelmed by the complexity of the expressions we need to deal with.\\ 

\noi We can simplify the calculations by using an affine form of the Hessian curve when we search for intersection multiplicities at the point $(0:0:1)$. Let $C$ be given by a defining polynomial $F(x,y,z)$. Setting $z=1$, we have $F(x,y,1)=f(x,y)$. We can then  derive the following matrix from the expression of the Hessian matrix on page \pageref{Hessian}  using Euler's identity \cite[p.66]{Fisc:01}.
$$\mathcal{H}= \begin{bmatrix} df & f_x & f_y \\ \\ (d-1)f_x & f_{xx} & f_{yx} \\ \\ (d-1)f_y & f_{xy} & f_{yy} \end{bmatrix}.$$
This radically simplifies the defining polynomial for the affine part of $H_C$.
\begin{equation*}
\begin{split}
H_f=&\;df(f_{xx}f_{yy}-f_{xy}^2)-(d-1)f_x^2f_{yy}\\
&+2(d-1)f_xf_yf_{xy}-(d-1)f_y^2f_{xx}.
\end{split}
\end{equation*}

\noi Following Fulton \cite[pp.74--75]{Fulton} we can now calculate $(C \cdot H_C)_{(0:0:1)}$ directly. We introduce the notation $(C \cdot H_C)_{(0:0:1)}=\mathrm{I}(f,H_f)$. For $p=(0:0:1)$ we can, using properties given by Fulton, simplify $\mathrm{I}(f,H_f)$ to
$$\mathrm{I}(f,2f_xf_yf_{xy}-f_x^2f_{yy}-f_y^2f_{xx}).$$
This is rarely enough simplification to find the intersection multiplicity.\\ 

\noi To further simplify the calculation of $(C \cdot H_C)_p$, observe that Corollary \ref{CHC} on page \pageref{CHC} hints to the following conjecture proposed by Ragni Piene, concerning the intersection multiplicity of $C$ and $H_C$ at a cusp or an inflection point $p$.

\begin{conj}[Intersection multiplicity]\label{Lhessian}
The intersection multiplicity $(C\cdot H_C)_p$ of a rational cuspidal curve $C$ and its Hessian curve $H_C$ in a cusp or an inflection point $p$ is given by
\begin{equation*}
\begin{split}
(C \cdot H_C)_p &= 6 \delta_p+ 2(m_p-1)+m_p^*-1\\
&=6 \delta_p+ m_p+r_p-3,
\end{split}
\end{equation*}
\noi where $m_p$ and $m_p^*$ denotes the multiplicity of $p$ and the dual point $p^*$, $\delta_p$ denotes the delta invariant of $p$, and $r_p$ denotes the intersection multiplicity $(T_p \cdot C)_p$.
\end{conj}
\noi Although all our examples imply that the conjecture holds, we have not been able to give a local proof. For specific curves and cusps, however, the result can be verified.

\begin{thm}[Intersection multiplicity for binomial curves]\label{LShessian} Let $p$ be a cusp on a binomial cuspidal curve $C$ given by a defining polynomial $F$, $$F=z^{n-m}y^m-x^n, \qquad \gcd(m,n)=1.$$ Then Conjecture \ref{Lhessian} holds. 
\end{thm}

\begin{pf}
The proof consists of calculating the expressions on each side of the equation in the conjecture, and subsequently observe that they coincide. The intersection multiplicity on the left hand side of the equation will be calculated by the method presented by Fischer in \cite[p.156]{Fisc:01}. The right hand side will be calculated directly.\\

\noi Let $C$ be a binomial cuspidal curve. The proof of Theorem \ref{LShessian} is identical for any cusp on this curve. We let $p$ be a cusp on $C$ with coordinates $(0:0:1)$. Then the affine part of $C$ around $p$ can be given by $f=y^m-x^n$, and we calculate the expressions in the conjecture.\\ 

\noi {\bf Left -- $\bf{(C \cdot H_C)_p}$}\\
\noi The polynomial $f$ has partial derivatives and double derivatives
\begin{equation*}
\begin{split}
f_x&=-nx^{n-1},\\
f_y&=my^{m-1},\\
f_{xx}&=-n(n-1)x^{n-2},\\
f_{xy}&=0,\\
f_{yy}&=m(m-1)y^{m-2}.
\end{split}
\end{equation*}
This leads to the following polynomial defining the Hessian curve,
\begin{equation*}
\begin{split}
H_f=&\;nff_{xx}f_{yy}-nff_{xy}^2-(n-1)f_x^2f_{yy}\\
&+2(n-1)f_xf_yf_{xy}-(n-1)f_y^2f_{xx},
\end{split}
\end{equation*}
\begin{equation*}
\begin{split}
H_f(x,y)=&-m(m-1)n^2(n-1)x^{n-2}y^{2m-2}\\
&+m(m-1)n^2(n-1)x^{2n-2}y^{m-2}\\
&-m(m-1)n^2(n-1)x^{2n-2}y^{m-2}\\
&+m^2n(n-1)^2x^{n-2}y^{2m-2}.\\
\end{split}
\end{equation*}

\noi The curve has a Puiseux parametrization around $p$ given by $$(C,p)=(t^m:t^n:1).$$ Substituting, we get
\begin{equation*}
\begin{split}
H_f(t)=&-m(m-1)n^2(n-1)t^{m(n-2)+n(2m-2)}\\
&+m(m-1)n^2(n-1)t^{m(2n-2)+n(m-2)}\\
&-m(m-1)n^2(n-1)t^{m(2n-2)+n(m-2)}\\
&+m^2n(n-1)^2t^{m(n-2)+n(2m-2)}.\\
\end{split}
\end{equation*}

\noi All terms in this polynomial in $t$ has degree $3mn-2m-2n$. Hence, this is the value of the intersection multiplicity, 
$$(C \cdot H_C)_p=3mn-2m-2n.$$

\noi {\bf Right -- {$\bf 6 \delta_p+ 2(m_p-1)+m_p^*-1$}}\\
\noi By Fischer \cite[pp.207,214]{Fisc:01}, we have that, for this particular point and curve, $$\delta_p=\frac{(n-1)(m-1)}{2}.$$ We know that $m_p=m$. We may also find $m_p^*$ by using the Puiseux parametrization of $p$ and finding the Puiseux parametrization of the dual point $p^*$, $$(C^*,p^*)=(t^{n-m}:1:t^n).$$ We find that $m_p^*=n-m$. 

The calculation is henceforth straightforward,
\begin{equation*}
\begin{split}
6 \delta_p+ 2(m_p-1)+m_p^*-1&=3(n-1)(m-1)+2(m-1)+(n-m)-1\\
&=3mn-3m-3n+3+2m-2+n-m-1\\
&=3mn-2m-2n,
\end{split}
\end{equation*} 
which is exactly what we wanted.
\begin{flushright}
$\square$
\end{flushright}
\end{pf}

\section{Reducible toric polar Cremona transformations}\label{fibprod}
Cremona transformations are phenomena which have been carefully studied and thoroughly described ~\cite{Hudson} ~\cite{Coolidge} ~\cite{Albe} ~\cite{Dolgachev}. There are, however, still unanswered questions concerning birational maps. A discussion with Professor Kristian Ranestad concerning Fibonacci curves resulted in an example of a Cremona transformation of a particular kind.\\ 

\noi Let $F_x$, $F_y$ and $F_z$ denote the partial derivatives of a homogeneous polynomial $F(x,y,z) \in \mathbb{C}[x,y,z]$. Then a Cremona transformation $\phi$ which can be written on the form $$\phi: (x:y:z) \longmapsto (xF_x:yF_y:zF_z)$$ is called a {\em{toric polar Cremona transformation}}. The polynomial $F$ is called a {\em toric polar Cremona polynomial}. Note that $F$ can be irreducible or reducible. The zero set $C=\mathcal{V}(F)$ of the polynomial $F$ is called a {\em toric polar Cremona curve}.\\

\noi In this section we will give an example of a series of reducible toric polar Cremona curves, where each irreducible component actually is a cuspidal curve. Additionally, we will give a constructive proof of this claim.\\

\noindent The $k$th Fibonacci curve was defined as $C_k=\mathcal{V}(y^{\f_{k}}-x^{\f_{k-1}}z^{\f_{k-2}})$ for all $k \geq 2$ on page \pageref{def:fibcurves}. Additionally, define $C_1=\mathcal{V}(y-z)$. $C_k$ is then a curve of degree $\f_k$ for all $k \geq 1$.

\begin{prp}[The Fibonacci example] The curve $C_k \cup C_{k-1}$ with defining polynomial 
\begin{equation*}
\begin{split}
F(k+1)&=(y^{\f_k}-x^{\f_{k-1}}z^{\f_{k-2}})(y^{\f_{k-1}}-x^{\f_{k-2}}z^{\f_{k-3}})\\
&=y^{\f_{k+1}}+x^{\f_k}z^{\f_{k-1}}-x^{\f_{k-2}}y^{\f_k}z^{\f_{k-3}}-x^{\f_{k-1}}y^{\f_{k-1}}z^{\f_{k-2}}.
\end{split}
\end{equation*}
is a reducible toric polar Cremona curve for every $k \geq 3$.
\end{prp}

\noindent To simplify notation we will henceforth write $F=F(k+1)$. With $F$ as above, the proposition has the consequence that the map $$\phi: (x:y:z) \longmapsto (xF_x:yF_y:zF_z)$$ is a Cremona transformation for $k \geq 3$,
\begin{equation*}
\begin{split}
xF_x&=\f_kx^{\f_k}z^{\f_{k-1}}-\f_{k-2}x^{\f_{k-2}}y^{\f_k}z^{\f_{k-3}}-\f_{k-1}x^{\f_{k-1}}y^{\f_{k-1}}z^{\f_{k-2}},\\
yF_y&=\f_{k+1}y^{\f_{k+1}}-\f_kx^{\f_{k-2}}y^{\f_k}z^{\f_{k-3}}-\f_{k-1}x^{\f_{k-1}}y^{\f_{k-1}}z^{\f_{k-2}},\\
zF_z&=\f_{k-1}x^{\f_k}z^{\f_{k-1}}-\f_{k-3}x^{\f_{k-2}}y^{\f_k}z^{\f_{k-3}}-\f_{k-2}x^{\f_{k-1}}y^{\f_{k-1}}z^{\f_{k-2}}.\\
\end{split}
\end{equation*}

\begin{pf} The proposition will be proved by construction. First we construct three polynomials $G_{100}$, $G_{010}$ and $G_{001}$ which define a Cremona transformation. They will define a Cremona transformation since they will be constructed by a sequence of linear changes of coordinates and standard Cremona transformations. Furthermore, we will show that there exists a linear change of coordinates which sends $xF_x$, $yF_y$ and $zF_z$ to the three polynomials. Then $xF_x$, $yF_y$ and $zF_z$ also define a Cremona transformation, and $F$ is a toric polar Cremona polynomial.\\

\noindent Let $L_{abc}=ax+by+cz$ be linear polynomials in $\mathbb{C}[x,y,z]$. Polynomials $G_{a'b'c'}$ of degree $\f_{k+1}$ may be constructed from these polynomials using successive linear changes of coordinates and standard Cremona transformations.

To the polynomials $ax+by+cz$, first apply a standard Cremona transformation. After removing linear factors, apply a linear transformation $\tau$ to the resulting polynomials. The transformation $\tau$ is given by the matrix $\mathcal{T}$, $$\mathcal{T} = \begin{bmatrix}-1&1&0\\-1&1&-1\\0&1&-1 \end{bmatrix}.$$ 
Repeat this process another $k-4$ times. A total of $k-3$ repetitions are required in order to eventually achieve the appropriate degree.

To the resulting polynomials apply another standard Cremona transformation. Follow this transformation by removing linear factors and then apply a linear transformation $\tau_1$ to the polynomials. Let $\tau_1$ be given by the matrix $\mathcal{T}_1$, $$\mathcal{T}_1 = \begin{bmatrix}1&-1&1\\1&-1&0\\0&1&0 \end{bmatrix}.$$

Last, apply another standard Cremona transformation and remove linear factors in the resulting polynomials. Then apply a linear transformation $\tau_0$. Let $\tau_0$ be given by the matrix $\mathcal{T}_0$, $$\mathcal{T}_0 = \begin{bmatrix}1&0&0\\-1&1&1\\0&0&-1 \end{bmatrix}.$$\\

\noi After applying the sequence of Cremona transformations, we get polynomials $G_{abc}(k+1)$ of degree $\f_{k+1}$. $G_{abc}$ can be given explicitly, and it has a nearly identical form for all $k \geq 3$. The only difference occurring is that the coefficients $c$ and $a$ switch places.\smallskip

{\samepage
\noindent For $k+1$ odd we have $$G_{abc}=-cy^{\f_{k+1}}+ax^{\f_k}z^{\f_{k-1}}+(b-a)x^{\f_{k-2}}y^{\f_k}z^{\f_{k-3}}+(c-b)x^{\f_{k-1}}y^{\f_{k-1}}z^{\f_{k-2}}.$$

\noindent For $k+1$ even we have
$$G_{abc}=-ay^{\f_{k+1}}+cx^{\f_k}z^{\f_{k-1}}+(b-c)x^{\f_{k-2}}y^{\f_k}z^{\f_{k-3}}+(a-b)x^{\f_{k-1}}y^{\f_{k-1}}z^{\f_{k-2}}.$$}

\noindent In either case we get polynomials of degree $\f_{k+1}$ on the form $G_{a'b'c'}$,
{\small
\begin{equation*}
\begin{split} G_{a'b'c'}=&\;a'y^{\f_{k+1}}+c'x^{\f_k}z^{\f_{k-1}}+(b'-c')x^{\f_{k-2}}y^{\f_k}z^{\f_{k-3}}-(a'+b')x^{\f_{k-1}}y^{\f_{k-1}}z^{\f_{k-2}}\\
=&\;a'(y^{\f_{k+1}}-x^{\f_{k-1}}y^{\f_{k-1}}z^{\f_{k-2}})\\
&+b'(x^{\f_{k-2}}y^{\f_k}z^{\f_{k-3}}-x^{\f_{k-1}}y^{\f_{k-1}}z^{\f_{k-2}})\\
&+c'(x^{\f_k}z^{\f_{k-1}}-x^{\f_{k-2}}y^{\f_k}z^{\f_{k-3}}).
\end{split}
\end{equation*}
}

\noi By construction, the three polynomials $G_{100}$, $G_{010}$ and $G_{001}$ define a Cremona transformation.\\

\noindent Observe that the polynomials $F$, $xF_x$, $yF_y$ and $zF_z$ have the same form as $G_{a'b'c'}$. 
\begin{equation*}
\begin{split}
F&=G_{101},\\
xF_x&=G_{0(\f_{k-1})(\f_k)},\\
yF_y&=G_{(\f_{k+1})(-\f_k)0},\\
zF_z&=G_{0(\f_{k-2})(\f_{k-1})}.
\end{split}
\end{equation*}

\noindent The Fibonacci numbers have the property that for any $k \geq 2$, $\f_{k-1}^2-\f_k\f_{k-2}=(-1)^k$ \cite[p.667]{Ore}. Using this property, an inspection of the coefficient matrix $\mathcal{C}$ of the three polynomials $xF_x$, $yF_y$ and $zF_z$ reveals that the determinant is nonzero. Hence, the polynomials are linearly independent.   
\begin{equation*}
\begin{split}
\det(\mathcal{C}) &= \det\begin{bmatrix} 0&\f_{k-1}&\f_k\\\f_{k+1}&-\f_k&0\\0&\f_{k-2}&\f_{k-1}\end{bmatrix}\\
&=-\f_{k+1}(\f_{k-1}^2-\f_{k-2}\f_k)\\
&=(-1)^{k+1}\f_{k+1}\\
&\neq 0 {\text { for all }} k.
\end{split}
\end{equation*}

\noi Furthermore, the matrix $\mathcal{C}$ defines a linear transformation $\tau_p$ which sends $xF_x$, $yF_y$ and $zF_z$ to $G_{100}$, $G_{010}$ and $G_{001}$. Then $\tau_p^{-1}$ composed with the sequence of Cremona transformations described earlier in this section is a Cremona transformation. Hence, $F$ is a toric polar Cremona polynomial.
\begin{flushright}
$\square$
\end{flushright}
\end{pf}

\pagebreak
\begin{rem}\label{rem:node}
Looking at the defining polynomials of the two curves $C_k$ and $C_{k-1}$, we see that they intersect in three points $(1:1:1)$, $(1:0:0)$ and $(0:0:1)$. Hence, the union of curves $C_k \cup C_{k-1}$ has three singularities with these coordinates. Using {\em Maple} we find that there is a nodal singularity in $(1:1:1)$. Furthermore, we have more complex singularities in $(1:0:0)$ and $(0:0:1)$. These two singularities are multiple points with two branches, but they have the same multiplicity and delta invariant as the corresponding cusps of the curve $C_{k+1}$. \\
\end{rem}

\begin{rem}
The Cremona transformation used in the construction of $G_{a'b'c'}$ of degree $\f_{k+1}$ was originally found in order to transform the ordinary Fibonacci curves $C_k$ to lines. In particular, the inverse of the transformation used in the proof will transform any $C_{k+1}$ to the line given by the defining polynomial $L_{111}$. Hence, it is possible to show that all Fibonacci curves are rectifiable. In the process of finding this suitable transformation, only one choice was made. The smooth point $(1:1:1)$ of $C_{k+1}$ was moved to $(0:1:0)$, which is a base point of the first standard Cremona transformation.\\
\end{rem}

\begin{rem}\label{rem:union}
The examples of reducible toric polar Cremona polynomials given by $F=G_{101}$ corresponds to the line given by $x-z$. Performing the described Cremona transformation $\psi$ on this single line breaks down immediately because one only considers the strict transform. By regarding the total transform instead, we get the following remarkable result when applying $\psi$.
$$\V(\psi(x-z))=C_k \cup C_{k-1} \bigcup_{i=0}^{k-3}(C_{k-2-i})^{2^i}.$$
Note that the union of the two curves of highest degree is precisely the reducible toric polar Cremona curve.\\
\end{rem}

\begin{rem}
The above results suggest further remarkable relations between the Fibonacci numbers. These results can also be proved by induction.\\

\noi 1. We have the relation \begin{equation} 2^{k-1}=\f_k+\f_{k-1}+\sum_{i=0}^{k-3}2^i\f_{k-2-i}. \label{1}\end{equation}

\begin{pf} This result comes from comparing the degree of the total transform in Remark \ref{rem:union} with the predicted degree after $k-1$ quadratic Cremona transformations. For $k=2$ this obviously holds, $$2^1=1+1=\f_2+\f_1.$$ The same is true for $k=3$, $$2^2=2+1+1=\f_3+\f_2+\f_1.$$ Now assume that (\ref{1}) holds for $k-1$. We then have $$2^{k-2}=\f_{k-1}+\f_{k-2}+\sum_{i=0}^{k-4}2^{i}\f_{k-3-i}.$$ Multiplying with $2$ on either side of the equation gives 
\begin{equation*}
\begin{split}
2^{k-1}&=2(\f_{k-1}+\f_{k-2})+2\sum_{i=0}^{k-4}2^i\f_{k-3-i}\\
&=\f_k+\f_{k-1}+\f_{k-2}+\sum_{i=0}^{k-4}2^{i+1}\f_{k-3-i}\\
&=\f_k+\f_{k-1}+\sum_{i=0}^{k-3}2^i\f_{k-2-i}.
\end{split}
\end{equation*}
\begin{flushright}
$\square$
\end{flushright}
\end{pf}

\noi 2. We have the relation \begin{equation} \f_k\f_{k-1}=\sum_{i=1}^{k-1}\f_i^2. \label{2}\end{equation}

\begin{pf} By B\'{e}zout's theorem and Remark \ref{rem:node} above, we are led to search for a connection between Fibonacci numbers on the form $$\f_k\f_{k-1}=\sum_{i=1}^3(C_k \cdot C_{k-1})_{p_i},$$ where $(C_k \cdot C_{k-1})_{p_i}$ denotes the intersection multiplicity of $C_k$ and $C_{k-1}$ in the three respective intersection points $p_i$.\\

\noi Since the intersection point $(1:1:1)$ is a node, we can conclude that $$(C_k \cdot C_{k-1})_{(1:1:1)}=1.$$ The other two intersection multiplicities can be found by direct calculation.

{\small
\begin{equation*}
\begin{split}
(C_k \cdot C_{k-1})_{(0:0:1)}&=\mathrm{I}(y^{\f_k}-x^{\f_{k-1}},y^{\f_{k-1}}-x^{\f_{k-2}})\\
&=\mathrm{I}(y^{\f_k}-x^{\f_{k-1}}-y^{\f_{k-2}}(y^{\f_{k-1}}-x^{\f_{k-2}}),y^{\f_{k-1}}-x^{\f_{k-2}})\\
&=\mathrm{I}(x^{\f_{k-2}}, y^{\f_{k-1}}-x^{\f_{k-2}})+\mathrm{I}(y^{\f_{k-1}}-x^{\f_{k-2}},y^{\f_{k-2}}-x^{\f_{k-3}})\\
&=\ldots\\
&=\f_{k-1}\f_{k-2}+\f_{k-2}\f_{k-3}+...+\f_2\f_1\\
&=\sum_{i=2}^{k-1}\f_i\f_{i-1}.
\end{split}
\end{equation*}

\begin{equation*}
\begin{split}
(C_k \cdot C_{k-1})_{(1:0:0)}&=\mathrm{I}(y^{\f_k}-z^{\f_{k-2}},y^{\f_{k-1}}-z^{\f_{k-3}})\\
&=\mathrm{I}(y^{\f_k}-z^{\f_{k-2}}-y^{\f_{k-2}}(y^{\f_{k-1}}-z^{\f_{k-3}}),y^{\f_{k-1}}-z^{\f_{k-3}})\\
&=\mathrm{I}(z^{\f_{k-3}}, y^{\f_{k-1}}-z^{\f_{k-3}})+\mathrm{I}(y^{\f_{k-1}}-z^{\f_{k-3}},y^{\f_{k-2}}-x^{\f_{k-4}})\\
&=\ldots\\
&=\f_{k-1}\f_{k-3}+\f_{k-2}\f_{k-4}+...+\f_3\f_1\\
&=\sum_{i=2}^{k-1}\f_i\f_{i-2}.
\end{split}
\end{equation*}
}

\noi We get the equality
\begin{equation*}
\begin{split}
\f_k\f_{k-1}&=1+ \sum_{i=2}^{k-1}\f_i(\f_{i-1}+\f_{i-2})\\
&=1+\sum_{i=2}^{k-1}\f_i^2\\
&=\sum_{i=1}^{k-1}\f_i^2.
\end{split}
\end{equation*}

\noi This equality is easily proved by induction. It obviously holds for $k=2$, \begin{equation*} \begin{split}\f_2\f_1&=\f_1^2 \\ 1 &=1. \end{split} \end{equation*} Assume that (\ref{2}) holds for $k-1$, $$\f_{k-1}\f_{k-2}=\sum_{i=1}^{k-2}\f_i^2.$$ Then adding $\f_{k-1}^2$ to each side of the equation gives \begin{equation*}
\begin{split}
\f_{k-1}^2+\f_{k-1}\f_{k-2}&=\f_{k-1}^2+\sum_{i=1}^{k-2}\f_i^2\\
\f_k\f_{k-1}&=\sum_{i=1}^{k-1}\f_i^2.
\end{split}
\end{equation*}
\begin{flushright}
$\square$
\end{flushright}
\end{pf}
\end{rem}

\pagebreak
\begin{rem} The toric polar Cremona curves are not strongly linked to cuspidal curves. According to Ranestad, however, there is at least one cuspidal curve which is also a toric polar Cremona curve. The tricuspidal quartic curve given by $C=\V(F)$ is a toric polar Cremona curve, $$F=x^2y^2+x^2z^2+y^2z^2-2xyz(x+y+z).$$ We have the Cremona transformation $$\psi:(x:y:z)\longmapsto (xF_x:yF_y:zF_z).$$
\end{rem}

\appendix

\chapter{Calculations and code}\label{calculationsandcode}

Using the programs {\em Maple} \cite{Maple} and {\em Singular} \cite{Singular}, a lot of properties of a curve can be found. In this appendix we will show some useful codes.\\ 

\noi Note that {\em Maple} and {\em Singular} require that we load packages before we execute the commands.
\begin{itemize}
\item[--] Packages in {\em Maple}.

{\scriptsize
\begin{verbatim}
with(algcurves):
with(LinearAlgebra):
with(VectorCalculus):
\end{verbatim}
}
\item[--] Packages in {\em Singular}.

{\scriptsize
\begin{verbatim}
LIB "all.lib";
\end{verbatim}
}
\end{itemize}

\section{General calculations}

\subsubsection{Multiplicity sequence}
The multiplicity sequence of a cusp $p$ in $(0:0:1)$ can be found with {\em Singular}.

{\scriptsize
\begin{verbatim}
ring r=0,(x,y),dp;
poly f=f(x,y);
displayMultsequence(f);
\end{verbatim}
}

\subsubsection{Intersection multiplicity}
Calculating the intersection multiplicity of two curves $C=\V(F)$ and $D=\V(G)$ at a point $p$ can be done directly as described in chapter \ref{TB}. We can also use the parametrization method described by Fischer in \cite[pp.147--169]{Fisc:01} or the polynomial tangent comparison algorithm given in Fulton \cite[pp.74--75]{Fulton}. Anyhow, the calculations quickly turn messy. The following code in {\em Singular} gives the intersection multiplicity of two curves in $(0:0:1)$.

{\scriptsize
\begin{verbatim}
ring r=0, (x,y), ls;
poly f=f(x,y);
poly g=g(x,y);
ideal I=f,g;
vdim(std(I));
\end{verbatim}
}

\subsubsection{Defining polynomial}
Given a parametrization of a curve $C=\V(F)$ on the form $$(x(s,t):y(s,t):z(s,t)),$$ we find the defining polynomial $F$ of $C$ by eliminating $s$ and $t$. This can easily be done by feeding {\em Singular}

{\scriptsize
\begin{verbatim}
ring R=0, (x,y,z,s,t), dp;
ideal C'=x-x(s,t),y-y(s,t),z-z(s,t);
eliminate(std(C'),st);
\end{verbatim}
}

\begin{ex}
The cuspidal cubic is given by the parametrization $$(s^3:st^2:t^3).$$ The defining polynomial can be found by

{\scriptsize
\begin{verbatim}
ring R=0, (x,y,z,s,t), dp;
ideal C'=x-(s3),y-(st2),z-(t3);
eliminate(std(C'),st);
_[1]=y3-xz2
\end{verbatim}
}
\end{ex}

\subsubsection{Dual curve}
To find the dual curve of a curve $C$ given by a polynomial $F(x,y,z)$, use {\em Singular} and the code

{\scriptsize
\begin{verbatim}
ring R=0,(x,y,z,s,t,u),dp;
poly F=F(x,y,z);
ideal I=F,s-diff(F,x),t-diff(F,y),u-diff(F,z);
short=0;
eliminate(std(I),xyz);
\end{verbatim}
}

\subsubsection{Hessian curve}
To find the Hessian curve of a curve $C$ given by a polynomial $F(x,y,z)$, use {\em Maple} and the code

{\scriptsize
\begin{verbatim}
F := F(x,y,z);
H := Hessian(F, [x, y, z]);
HC := Determinant(H);
\end{verbatim}
}

\subsubsection{Singularities of a curve}
To find the singularities of a curve $C$ given by a polynomial $F(x,y,z)$, use {\em Maple} and the code

{\scriptsize
\begin{verbatim}
F := F(x,y,z);
singularities(F,x,y);
\end{verbatim}
}

\noi For every singularity $p$ of $C$, the output is on the form

$$[[x,y,z],m_p,\delta_p,\# \text{Branches}].$$

\subsubsection{Intersection points of curves}
Intersection points of curves $C=\V(F(x,y,z))$ and $D=\V(G(x,y,z))$ can be found using {\em Maple} and the code

{\scriptsize
\begin{verbatim}
F := F(x,y,z);
G := G(x,y,z);
singularities(F*G,x,y);
\end{verbatim}
}

\noi The output will give the coordinates of the intersections.

\section{Projections}

\subsubsection{The tangent developable}
Defining polynomials for the tangent developable can be found in every degree with the help of {\em Singular}. 

\begin{ex} For degree $d=4$ the tangent developable $T_4$ can be found by feeding {\em Singular}

{\scriptsize
\begin{verbatim}
ring r=0, (x0,x1,x2,x3,x4,a,b,t,s), dp;
ideal T=(x0-4as3,x1-(3as2t+bs3), x2-(2ast2+2bs2t), x3-(at3+3bst2), x4-4bt3);
ideal TD=eliminate(std(T), abst);
std(TD);
_[1]=3*x2^2-4*x1*x3+x0*x4
_[2]=2*x1*x2*x3-3*x0*x3^2-3*x1^2*x4+4*x0*x2*x4
_[3]=8*x1^2*x3^2-9*x0*x2*x3^2-9*x1^2*x2*x4+14*x0*x1*x3*x4-4*x0^2*x4^2
\end{verbatim}
}
\end{ex}

\subsection{Code for analysis of projections}
The analysis of projection centers presented in Section \ref{cuspidalprojections4}, was done with {\em Maple} and {\em Singular}. The calculations are equivalent for all curves, just substitute for the line $L$. 

\subsubsection{\bf Bicuspidal quartic -- $[(2_2),(2)]$} 
The cuspidal quartic curve with two cusps, one $A_4$-cusp and one $A_2$-cusp, and one inflection point of type $1$ is given by the parametrization $$(s^4+s^3t:s^2t^2:t^4).$$ 

\noi The projection center $L$ can be presented as the intersection of three linear hyperplanes simply by reading off the parametrization, $$L=\V(x_0+x_1,x_2,x_4),$$
$$A_L=\left[ \begin {array}{ccccc} 1&1&0&0&0\\\noalign{\medskip}0&0&1&0&0
\\\noalign{\medskip}0&0&0&0&1\end {array} \right].$$
The basis vectors of the kernel of the projection, which are the components of $K_L$, can be found using {\em Maple},

{\scriptsize
\begin{verbatim}
with(linalg):
A := Matrix([[1, 1, 0, 0, 0], [0, 0, 1, 0, 0], [0, 0, 0, 0, 1]]);
K := kernel(A);
{[0 0 0 1 0], [-1 1 0 0 0]}
\end{verbatim}
}

\noi Below we present the code used in {\em Singular} to verify each claim in the analysis.\\

\noi Initial code in {\em Singular} is

{\scriptsize
\begin{verbatim}
LIB "all.lib";
ring r=0,(s,t),dp;
\end{verbatim}
}
\begin{itemize}
\item[--] $L \cap C_4=\emptyset.$

{\scriptsize
\begin{verbatim}
matrix C[3][5]=s4,s3t,s2t2,st3,t4,-1,1,0,0,0,0,0,0,1,0;
ideal I=(minor(C,3));
solve(std(I));
[1]:
   [1]:
      0
   [2]:
      0
\end{verbatim}
}

\item[--] $L \cap T_4=\{p_1,p_2\}$ for $(s:t)=\{(1:0),(0:1)\}$.

{\scriptsize
\begin{verbatim}
matrix T[4][5]=4s3,3s2t,2st2,t3,0,0,s3,2s2t,3st2,4t3,-1,1,0,0,0,0,0,0,1,0;
ideal I=minor(T,4);
ideal Is=I,s-1;
ideal It=I,t-1;
solve(std(Is));
[1]:
   [1]:
      1
   [2]:
      0
solve(std(It));
[1]:
   [1]:
      0
   [2]:
      1
\end{verbatim}
}

\item[--] $L \cap O^2(s,t)=\{p_1,p_2,p_3\}$ for $(s:t)=\{(1:0),(0:1),(1:\textstyle{-\frac{8}{3}})\}$. 
Observe that $(1:-\textstyle{\frac{8}{3}})=(-\textstyle{\frac{3}{8}}:1) \text{ in } \mathbb{P}^1$.
 
{\scriptsize
\begin{verbatim}
matrix O_2[5][5]=6s2,3st,t2,0,0,0,3s2,4st,3t2,0,0,0,s2,3st,6t2,-1,1,0,0,0,0,0,0,1,0;
ideal I=det(O_2);
ideal Is=I,s-1;
ideal It=I,t-1;
solve(std(Is));
[1]:
   [1]:
      1
   [2]:
      -2.66666667
[2]:
   [1]:
      1
   [2]:
      0
solve(std(It));
[1]:
   [1]:
      -0.375
   [2]:
      1
[2]:
   [1]:
      0
   [2]:
      1
\end{verbatim}
}

\item[--] $L \nsubseteq O^2(s,t)$.

{\scriptsize
\begin{verbatim}
matrix O_2[5][5]=6s2,3st,t2,0,0,0,3s2,4st,3t2,0,0,0,s2,3st,6t2,-1,1,0,0,0,0,0,0,1,0;
ideal I=minor(O_2,4);
ideal Is=I,s-1;
ideal It=I,t-1;
solve(std(Is));
 ? ideal not zero-dimensional
solve(std(It));
 ? ideal not zero-dimensional
\end{verbatim}
}

\item[--] $L \subset O^3(s,t)$ for $(s:t)=\{(1:0)\}$.

{\scriptsize
\begin{verbatim}
matrix O_3[6][5]=0,3s,2t,0,0,0,0,2s,3t,0,4s,t,0,0,0,0,0,0,s,4t,-1,1,0,0,0,0,0,0,1,0;
ideal I=minor(O_3,5);
ideal Is=I,s-1;
ideal It=I,t-1;
solve(std(Is));
[1]:
   [1]:
      1
   [2]:
      0
solve(std(It));
  ? ideal not zero-dimensional
\end{verbatim}
}
\end{itemize}

\bibliographystyle{abbrv}
\bibliography{bib}

\end{document}